\documentclass[letterpaper]{article}

\usepackage[margin=1.5in]{geometry}

\usepackage{graphicx,epsf}
\usepackage{bbm,color}
\usepackage{amsmath,amssymb,amsfonts}
\usepackage{multirow}
\usepackage{bm}
\usepackage{algorithmic}
\usepackage[ruled,vlined,linesnumbered]{algorithm2e}
\usepackage{hyperref}
\usepackage{amsthm}
\usepackage{epstopdf}
\usepackage[svgnames]{xcolor}
\usepackage{caption}

\newcommand{\RR}{{\mathbb R}}

\newcommand{\N}{\mathbb{Z}}

\newcommand{\real}{\mathbb{R}}

\newcommand{\E}{\mathbb{E}}

\newcommand{\bbR}{\mathbb{R}}
\newcommand{\bbN}{\mathbb{N}}

\usepackage[ruled,vlined,linesnumbered]{algorithm2e}

\newtheorem{theorem}{Theorem}[section]

\newtheorem{remark}[theorem]{Remark}
\newtheorem{corollary}[theorem]{Corollary}
\newcommand{\GD}[2][ForestGreen]{}

\graphicspath{
{./Figures/}
}

\date{}

\begin{document}

\title{Non-Linear Wavelet Regression and Branch \& Bound Optimization for the Full Identification of Bivariate Operator Fractional Brownian Motion
}

\author{Jordan Frecon\thanks{Jordan Frecon (Corresponding author), Nelly Pustelnik and Patrice Abry are with Univ Lyon, Ens de Lyon, Univ Claude Bernard, CNRS, Laboratoire de Physique, F-­69342
Lyon, France (e-mail: firstname.lastname@ens-lyon.fr).}, Gustavo Didier\thanks{Gustavo Didier is with Mathematics Department, Tulane University, New Orleans, LA 70118, USA, (e-mail: gdidier@tulane.edu).}, Nelly Pustelnik$^*$, and Patrice Abry$^*$\thanks{Work supported by French ANR grant AMATIS  $\# 112432 $, 2010-2014 and by the prime award no.\ W911NF-14-1-0475 from the Biomathematics subdivision of the Army Research Office, USA. G.D. gratefully acknowledges the support of Ens de Lyon.
}
\\
}

\maketitle

\begin{abstract}
Self-similarity is widely considered the reference framework for modeling the scaling properties of real-world data.
However, most theoretical studies and their practical use have remained univariate. Operator Fractional Brownian Motion (OfBm) was recently proposed as a multivariate model for self-similarity. Yet it has remained seldom used in applications because of serious issues that appear in the joint estimation of its numerous parameters.
While the univariate fractional Brownian motion requires the estimation of two parameters only, its mere bivariate extension already involves 7 parameters which are very different in nature. The present contribution proposes a method for the full identification of bivariate OfBm (i.e., the joint estimation of all parameters)
through an original formulation as a non-linear wavelet regression coupled with a custom-made Branch $\&$ Bound numerical scheme.
The estimation performance (consistency and asymptotic normality) is mathematically established and numerically assessed by means of Monte Carlo experiments.
The impact of the parameters defining OfBm on the estimation performance as well as the associated computational costs are also thoroughly investigated.
\end{abstract}

\clearpage
\section{Introduction}
\label{sec.intro}

\noindent {\bf Scale invariance and self-similarity.} Scale invariance, or scaling, is now recognized as an ubiquitous property in a variety of real-world applications which are very different in nature (cf. e.g., \cite{WENDT:2007:E} and references therein for reviews). The so-named scale invariance paradigm is based on the assumption that temporal dynamics in data are not driven by one, or a few, representative time scales, but by a large continuum of them.
Self-similar stochastic processes provide the basal mathematical framework for the modeling of scaling phenomena. In essence, self-similarity states that a signal $X$ cannot be distinguished from any of its dilated copies (cf.\ e.g., \cite{Samorodnitsky1994}):
\begin{equation}
\label{eq:oss}
 \{X(t)\}_{ t \in \real}  \stackrel{fdd}{=} \{a^H X(t/a)\}_{t \in \real}, \forall a > 0,
 \end{equation}
 where $\stackrel{fdd}{=} $ stands for the equality of finite dimensional distributions. The key information on scale-free dynamics is summed up under a single parameter $0 <H <1$, called the Hurst exponent, whose estimation is the main goal in scaling analysis. Amongst the numerous estimators of $H$ proposed in the literature (cf. e.g., \cite{Bardet2003a} for a review), one popular methodology draws upon the computation of the sample variance of a set of multiscale quantities (e.g., wavelet coefficients) $T_X(a,t)$ that behave like a power law with respect to the scale $a$:
 \begin{equation}
 \label{equ:wavcoef}
 \sum_t T_X^2(a,t) \simeq a^{2H+1}.
 \end{equation}
 In view of the relation \eqref{equ:wavcoef}, $H$ can be estimated by means of a linear regression in log--log coordinates (cf. e.g., \cite{Veitch1999}).
 
 Fractional Brownian motion (fBm) $B_H(t)$ -- i.e., the only Gaussian, self-similar, stationary-increment process -- has massively been used as a reference process in the modeling of scaling properties in univariate real-world signals. \\
 
\noindent {\bf Multivariate scaling.}  Notwithstanding its theoretical and practical importance, fBm falls short of providing an encompassing modeling framework for scaling because most modern contexts of application involve the recording of multivariate time series
that hence need to be jointly analyzed.
The construction of a comprehensive multivariate estimation paradigm is still an open problem in the literature.
The so-named Operator fractional Brownian motion (OfBm), henceforth denoted by $\underline{B}_{W, \underline{H}}(t) $, is a natural extension of fBm.
It was recently been defined and studied in \cite{Amblard_P-0_2011_j-ieee-tsp_imfbm,didier:pipiras:2012,Didier_G_2011_bernouilli_irpofbm,Coeurjolly_J-F_2013_ESAIM_wamfbm} as the only Gaussian, multivariate self-similar process with stationary increments. Multivariate self-similarity translates into the relation:
\begin{equation}\label{equ:ssmulti}
 \{\underline{B}_{W, \underline{H}}(t)\}_{ t \in \real}  \stackrel{fdd}{=} \{a^{\underline{\underline{H}}} \underline{B}_{W,\underline{H}}(t/a)\}_{t \in \real}, \forall a > 0,
 \end{equation}
 where the scaling exponent consists of a Hurst matrix $\underline{\underline{H}} = W \makebox{diag} \underline{H} W^{-1}$. In the latter expression, $W$ represents a $ P \times P $ invertible matrix, and $ \underline{H} $ is a $P$-dimensional vector of Hurst eigenvalues, where $ a^{\underline{\underline{H}}} := \sum_{k=0}^{+\infty} \log^k(a) \underline{\underline{H}}^k/k !$. The full parametrization of OfBm further requires a $ P \times P $ point-covariance matrix $\Sigma$. OfBm remains so far rarely used in applications, mostly because its actual use requires, in a general setting, the estimation of $P + P^2 + P(P-1)/2$ parameters which are very different in nature (cf. Section~\ref{sec.OFBM}). In particular, Eq.~\eqref{equ:wavcoef} above results in a mixture of power laws (cf. Eqs.~\eqref{eq:E11}-\eqref{eq:E22} in Section~\ref{sec.OFBM}).
However, the identification of OfBm has been thoroughly studied in the entry-wise scaling case (corresponding to a diagonal mixing matrix $W$) \cite{Amblard_P-0_2011_j-ieee-tsp_imfbm} 
and often used in applications (cf.\ e.g., \cite{Achard2008,Ciuciu14}).
Identification has also been recently achieved under a non-diagonal mixing matrix $W$, yet with more restrictive assumptions on $\Sigma$ \cite{Helgason:Didier:Abry:2015}.
Even more recently, \cite{AbryDidier2015} proposed a general estimator for the vector of Hurst eigenvalues $ \underline{H} $ in the bivariate setting, yet requiring additional assumptions for the estimation of the extra parameters $W$ and $\Sigma$.
The full identification of OfBm without parametric assumptions has remained an open issue. \\

\noindent {\bf Goals, contributions and outline.}  In this work, our contribution is two-fold.
First, the full identification of bivariate ($P=2$) OfBm (Biv-OfBm) is formulated as a non-linear wavelet domain regression.
Second, an algorithmic solution for the associated optimization problem is devised by means of a Branch \& Bound procedure, which is essential in view of the highly non-convex nature of the latter.
To this end, definitions and properties of Biv-OfBm are recapped in Section~\ref{sec.OFBM}.
A parsimonious parametrization of the process is also proposed, which prevents the potential parametric under-determination of Biv-OfBm \cite{didier:pipiras:2012}.
In Section~\ref{sec.WOFBM}, the properties of the wavelet coefficients and of the wavelet spectrum of Biv-OfBm are explicitly laid out and computed.
This provides a mathematical framework for the proposed estimation method.
The full identification of Biv-OfBm is formulated as a minimization problem whose solution is developed based on a Branch \& Bound strategy (cf. Section~\ref{sec.BB}). The consistency and asymptotic normality of the proposed estimator is mathematically established in the general multivariate setting $P \geq 2$ (cf.\ Section~\ref{sec.Co}), and numerically assessed in the bivariate setting $P=2$ by means of Monte Carlo experiments conducted on large numbers of synthetic Biv-OfBm paths (cf.\ Section~\ref{sec.sim}). Comparisons with the Hurst eigenvalues estimators proposed in \cite{AbryDidier2015} are also reported. The routines for the identification and synthesis of OfBm will be made publicly available at the time of publication.

\section{Bivariate Operator fractional Brownian motion}
\label{sec.OFBM}
\subsection{Definitions}

\subsubsection{Preamble}
The most general definitions of OfBm were formulated in \cite{Amblard_P-0_2011_j-ieee-tsp_imfbm,didier:pipiras:2012,Didier_G_2011_bernouilli_irpofbm,Coeurjolly_J-F_2013_ESAIM_wamfbm} as the only multivariate Gaussian, self-similar (i.e., satisfying Eq.~\eqref{equ:ssmulti}) process with stationary increments.
Targeting real-world data and applications, the present contribution is restricted to the (slightly) narrower class of \emph{time reversible} OfBm (cf. \cite{Amblard_P-0_2011_j-ieee-tsp_imfbm,didier:pipiras:2012,Didier_G_2011_bernouilli_irpofbm,Coeurjolly_J-F_2013_ESAIM_wamfbm}) whose scaling exponent matrix $ \underline{\underline{H}} $ can be diagonalized as $\underline{\underline{H}} = W \makebox{diag} \underline{H} W^{-1}$, where $W$ is an invertible matrix.
The definitions and properties of OfBm are stated only in the bivariate setting.

\subsubsection{Entry-wise scaling OfBm}

The entry-wise scaling, time-reversible OfBm $\{X(t)\}_{t\in\mathbb{R}} \equiv \underline{B}_{\makebox{Id},\underline{H}}(t)$ is defined by the condition $W = $ Id. Hence, the Hurst exponent has the form $\underline{\underline{H}}=\mbox{diag}(h_1,h_2)$, $0<h_1\leq h_2<1$.
Let $\Sigma_X \equiv \E X(1) X(1)^* $ denote the point covariance matrix of $X$ with entries $ \sigma_{{\rm x}_m}\sigma_{{\rm x}_n} \rho_{{\rm x}_m,{\rm x}_n} $, where $ \sigma^2_{{\rm x}_m} $ is the variance of component $m$ and $ \rho_{{\rm x}_m,{\rm x}_n} $ is the correlation between components $m$ and $n$.
It was shown in \cite{Ambard_P-O_2013_bsmfbm}, \cite{Didier_G_2011_bernouilli_irpofbm} that the process $X$ is well-defined (i.e., that its covariance matrix $ \E X(t) X(s)^* $ is always positive definite) if and only if (with $ \rho_{\rm x} \equiv  \rho_{{\rm x}_1,{\rm x}_2} $):
\begin{multline}
\label{eq:condPositiveDefinite}
g(h_1, h_2,\rho_{\rm x})  \equiv \Gamma(2h_1+1)\Gamma(2h_2+1)\sin(\pi h_1)\sin(\pi h_2) \\- \rho_{\rm x}^2\Gamma(h_1+h_2+1)^2\sin^2(\pi(h_1+h_2)/2) >0.
\end{multline}
For entry-wise scaling OfBm, self-similarity in Eq.~\eqref{equ:ssmulti} simplifies to:
\begin{equation}
\{X_1(at),X_2(at)\}_{t\in\mathbb{R}}  \stackrel{fdd}{=} \{ a^{h_1} X_1(t), a^{h_2}X_2(t)\}_{t\in\mathbb{R}},\; \forall a>0.
\label{eq:ossew}
\end{equation}
The estimation of the parameters $(h_1,h_2,\rho_{\rm x},\sigma_{\rm x_1},\sigma_{\rm x_2})$, which fully characterize the process, can thus be conducted by following univariate-type strategies, i.e., by making use of extensions of Eq.~\eqref{equ:wavcoef} to all auto- and cross-components (cf.\ \cite{Amblard_P-0_2011_j-ieee-tsp_imfbm,Coeurjolly_J-F_2013_ESAIM_wamfbm} for a theoretical study estimation performance, or \cite{Ciuciu14} for wavelet-based estimation on real-world data).

\subsubsection{Mixing}

Let $W$ denote a $2 \times 2$ invertible matrix, hereinafter called the mixing matrix. OfBm is defined as $\{\underline{B}_{W, \underline{H}}(t) \equiv Y(t)\}_{t\in\mathbb{R}} =\{W \underline{B}_{\makebox{Id},\underline{H}}(t)  \equiv W X(t)\}_{t\in\mathbb{R}}$. Following \cite{Amblard_P-0_2011_j-ieee-tsp_imfbm,Didier_G_2011_bernouilli_irpofbm}, it is straightforward to show that $Y$ is self-similar as in Eq.~\eqref{equ:ssmulti}, with $\underline{\underline{H}} = W \makebox{diag} \underline{H} W^{-1}$. When $W$ is not diagonal, OfBm is no longer entry-wise scaling. Instead, the entry-wise scaling behavior of OfBm consists of mixtures of univariate power laws (cf.\ Eqs.~(\ref{eq:E11})-(\ref{eq:E22})). For this reason, the construction of estimators in the bivariate setting cannot rely on a direct extension of a univariate procedure.

 \subsection{Properties}

\subsubsection{Under-determination}
\label{sec:underD}
Because $W$ is invertible, one can show that, for $\Sigma_{\rm y}(t) \equiv \mathbb{E} Y(t)Y(t)^{*} $,
\begin{equation}
\Sigma_{\rm y}(t) = W  \mathbb{E} X(t)X(t)^{*} W^{*} \equiv  W \Sigma_{\rm x}(t) W^{*},
\label{eq:generalCovX}
\end{equation}
which reveals three forms of under-determination in the parametrization
of OfBm:\\
\indent i) Writing $T_X = \mbox{diag}(\sigma_{{\rm x}_1},  \sigma_{{\rm x}_2})$ and $ \Sigma_X = T_X C_X T^*_X$, where $C_X \equiv \{ 1 \,  \rho_{\rm x}  ; \rho_{\rm x}   \, 1\}$ is the correlation matrix of $X$, one cannot discriminate between $\underline{Y} = W \underline{X} $ and $\underline{Y} = W' \underline{X}' $, where $ W' = W T_X $ and $\underline{X}' =  T^{-1}_X \underline{X}$. \\ 
\indent ii) Let $\Pi $ denote a $2 \times 2 $ permutation matrix, i.e., there is only one non-zero entry (equal to $1$) per column or line.
Then, $\underline{Y} = W \underline{X} = W' \underline{X'} $, where $W' = W \Pi $ and $X' = \Pi^T X$.\\
\indent iii) Let $S$ be a diagonal matrix with entries $\pm 1$ and $X' = S X$, then $\underline{Y} = W \underline{X} = W' \underline{X'} $, where $W' \equiv W S $.\\

\subsubsection{Parametrization}
\label{sec.param}

To fix the parametric under-determination of OfBm, we adopt the following conventions:
i) the columns of $W$ are normalized to 1; ii) $h_1\leq h_2$; iii) the diagonal entries of $W$ are positive. This leads us to propose the following generic 7-dimensional parametrization $\Theta=(h_1,h_2,\rho_{\rm x},\sigma_{\rm x_1},\sigma_{\rm x_2},\beta,\gamma)$ of Biv-OfBm $\{Y(t)\}_{t\in\mathbb{R}}$:
\begin{equation}
W=\left(\begin{matrix}
 \frac{1}{\sqrt{1+\gamma^2}} & \frac{\beta}{\sqrt{1+\beta^2}} \\
  \frac{-\gamma}{\sqrt{1+\gamma^2}} & \frac{1}{\sqrt{1+\beta^2}}
 \end{matrix}\right),\;
\Sigma_X =
\left(\begin{matrix}
\sigma^2_{\rm x_1} &  \sigma_{\rm x_1}\sigma_{\rm x_2} \rho_{\rm x} \\
  \sigma_{\rm x_1}\sigma_{\rm x_2} \rho_{\rm x} & \sigma^2_{\rm x_2}
 \end{matrix}\right).
 \label{eq:mixtureMatrix}
\end{equation}

\section{Wavelet Analysis of OfBm}
\label{sec.WOFBM}
 \subsection{Multivariate discrete wavelet transform (DWT)}

Let $\psi_0$ be a mother wavelet, namely, $\psi_0 \in L^2(\bbR)$ and $ \int_\real t^k \psi_0(t) dt \equiv 0$, $k = 0,1,\hdots,N_{\psi}-1$.
Let $\{ \psi_{j,k}(t) = 2^{-j/2} \psi_0(2^{-j}t-k) \}_{(j,k) \in \N^2}$ denote the collection of dilated and translated templates of $\psi_0$ that forms an orthonormal basis of $L^2(\bbR)$. The multivariate DWT coefficients of $\{Y(t)\}_{t\in\mathbb{R}}$ are defined as $ (D_{\rm y}(j,k)) \equiv (D_{\rm y_1}(j,k), D_{\rm y_2}(j,k)) $, where
\begin{equation}
\label{eq:multquant}
\displaystyle D_{\rm y_m}(j,k) = \int_\mathbb{R} 2^{-j/2} \psi_0(2^{-j}t-k) Y_m(t)\mathrm{d}t, \quad m=1,2.
\end{equation}
For a detailed introduction to wavelet transforms, interested readers are referred to, e.g., \cite{Mallat_S_2008_book_awtspsw}.

 \subsection{Wavelet spectrum}

The properties of the wavelet coefficients of OfBm in a $P$-variate setting were studied in detail in \cite{AbryDidier2015}. Here, we only recall basic properties and expand on what is needed for actual full identification  (i.e., the estimation of all parameters entering its definition) of Biv-OfBm.

\subsubsection{Mixture of power laws}

From Eq.~\eqref{eq:ossew} and $Y(t)=WX(t)$, it can be shown that the wavelet spectrum reads:
\begin{align}
\mathbb{E} D_{\rm y}(j,k) D_{\rm y}(j,k)^*
& = W 2^{j (\underline{H}+ \makebox{Id} /2)} E_0 2^{j (\underline{H}^*+ \makebox{Id} /2)} W^*
\label{eq:covDyGeneral}
\end{align}
with $E_0  \equiv \mathbb{E} D_{\rm x}(0,k) D_{\rm x}(0,k)^*=$
 \begin{equation}   \left(\begin{matrix}
  \sigma_{{\rm x}_1}^2\eta_{h_1} & \rho_{\rm x}\sigma_{{\rm x}_1}\sigma_{{\rm x}_2} \eta_{\frac{h_1+h_2}{2}} \\
  \rho_{\rm x}\sigma_{{\rm x}_1}\sigma_{{\rm x}_2} \eta_{\frac{h_1+h_2}{2}} & \sigma_{{\rm x}_2}^2\eta_{h_2}
 \end{matrix}\right)
\end{equation}
\begin{equation}
\makebox{and } \eta_h = -\frac{1}{2} \int_\mathbb{R} |u|^{2h} \mathrm{d}u \int_\mathbb{R} \psi_0(v) \psi_0(v-u)^*\mathrm{d}v > 0.
\label{eq:etah}
\end{equation}

The OfBm parametrization proposed above yields the following explicit form of the wavelet spectrum:
\begin{equation}
\label{eq:ensave}
 \mathbb{E} D_{\rm y}(j,k)D_{\rm y}(j,k)^*
 \equiv E(2^j,\Theta)= \left( \begin{matrix} E_{11}(2^j,\Theta) & E_{12}(2^j,\Theta) \\  E_{12}(2^j,\Theta)& E_{22}(2^j,\Theta) \end{matrix} \right),
 \end{equation}
 \begin{multline}\label{eq:E11}
\makebox{with }  E_{11}(2^j,\Theta) = (1+\gamma^2)^{-1}\sigma_{{\rm x}_1}^2\eta_{h_1} 2^{j(2h_1+1)}\\+ 2\beta (1+\beta^2)^{-1/2} (1+\gamma^2)^{-1/2} \rho_{\rm x}\sigma_{{\rm x}_1}\sigma_{{\rm x}_2}\eta_{\frac{h_1+h_2}{2}} 2^{j(h_1+h_2+1)} \\+   \beta^2 (1+\beta^2)^{-1} \sigma_{{\rm x}_2}^2\eta_{h_2}2^{j(2h_2+1)},
 \end{multline}
 \begin{multline}\label{eq:E12}
 E_{12}(2^j,\Theta) = - \gamma (1+\gamma^2)^{-1} \sigma_{{\rm x}_1}^2\eta_{h_1}2^{j(2h_1+1)}\\+ (1-\beta\gamma)(1+\beta^2)^{-1/2} (1+\gamma^2)^{-1/2}  \rho_{\rm x}\sigma_{{\rm x}_1}\sigma_{{\rm x}_2}\eta_{\frac{h_1+h_2}{2}}2^{j(h_1+h_2+1)}\\+   \beta(1+\beta^2)^{-1} \sigma_{{\rm x}_2}^2\eta_{h_2}2^{j(2h_2+1)},
 \end{multline}
 \begin{multline}\label{eq:E22}
   E_{22}(2^j,\Theta) = \gamma^2(1+\gamma^2)^{-1} \sigma_{{\rm x}_1}^2\eta_{h_1} 2^{j(2h_1+1)} \\ -  2\gamma (1+\beta^2)^{-1/2} (1+\gamma^2)^{-1/2}  \rho_{\rm x}\sigma_{{\rm x}_1}\sigma_{{\rm x}_2}\eta_{\frac{h_1+h_2}{2}}2^{j(h_1+h_2+1)} \\+   (1+\beta^2)^{-1} \sigma_{{\rm x}_2}^2\eta_{h_2}2^{j(2h_2+1)}.
 \end{multline}

\subsubsection{Further under-determination}
\label{sec:funderD}

Eqs.~\eqref{eq:E11}, \eqref{eq:E12} and \eqref{eq:E22} reveal that $E_{11}(2^j,\Theta)$, $|E_{12}(2^j,\Theta)|$ and $E_{22}(2^j,\Theta)$ are invariant under the transformation $(\beta,\gamma,\rho_{\rm x}) \rightarrow -(\beta,\gamma,\rho_{\rm x})$. Therefore, the definition of $\rho_{\rm x} $ can be restricted to $\rho_{\rm x}\geq 0$. 

\subsection{Empirical wavelet spectrum}

\begin{figure}[!t]
\centerline{
\includegraphics[width=0.22\linewidth]{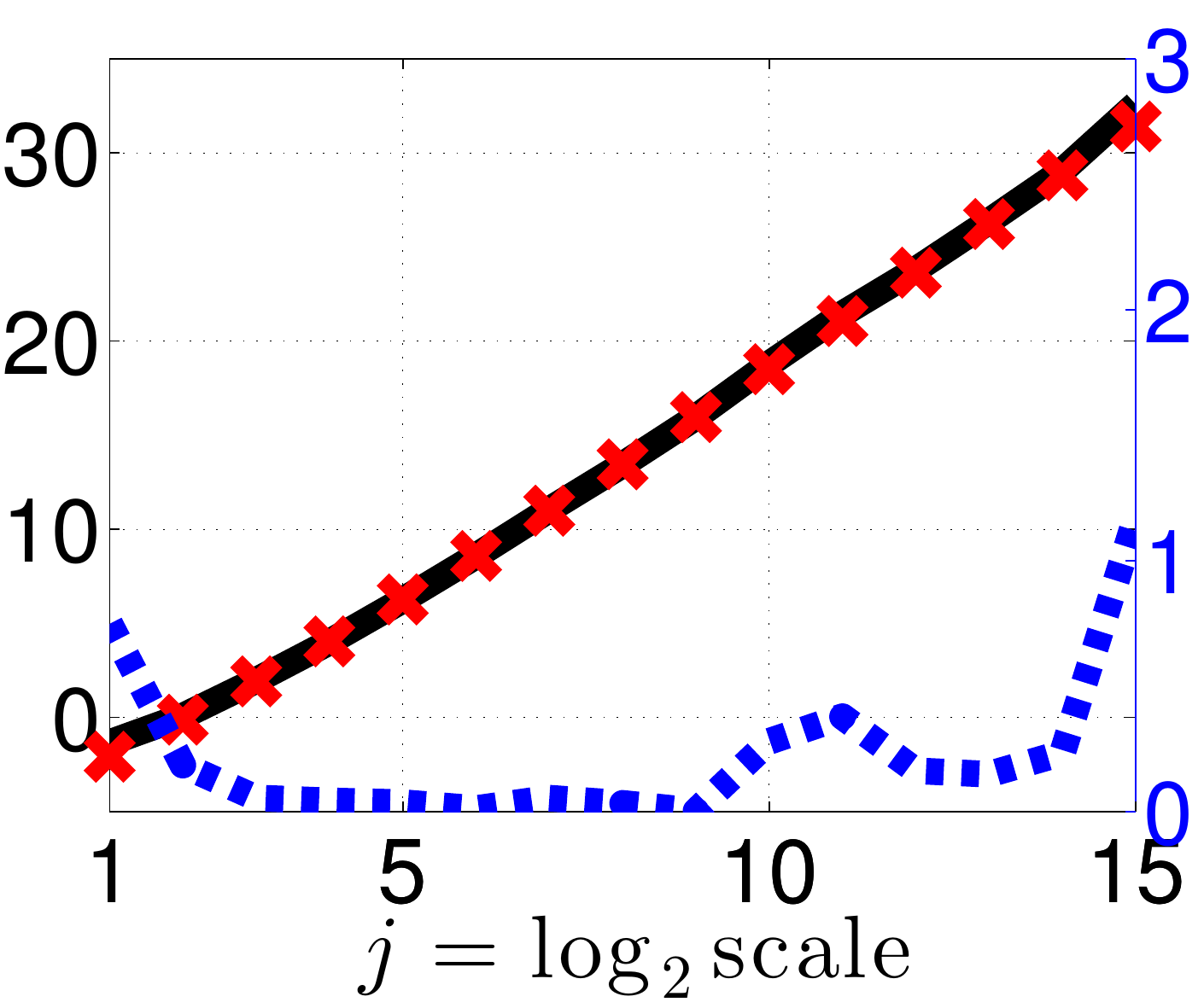}\qquad
\includegraphics[width=0.22\linewidth]{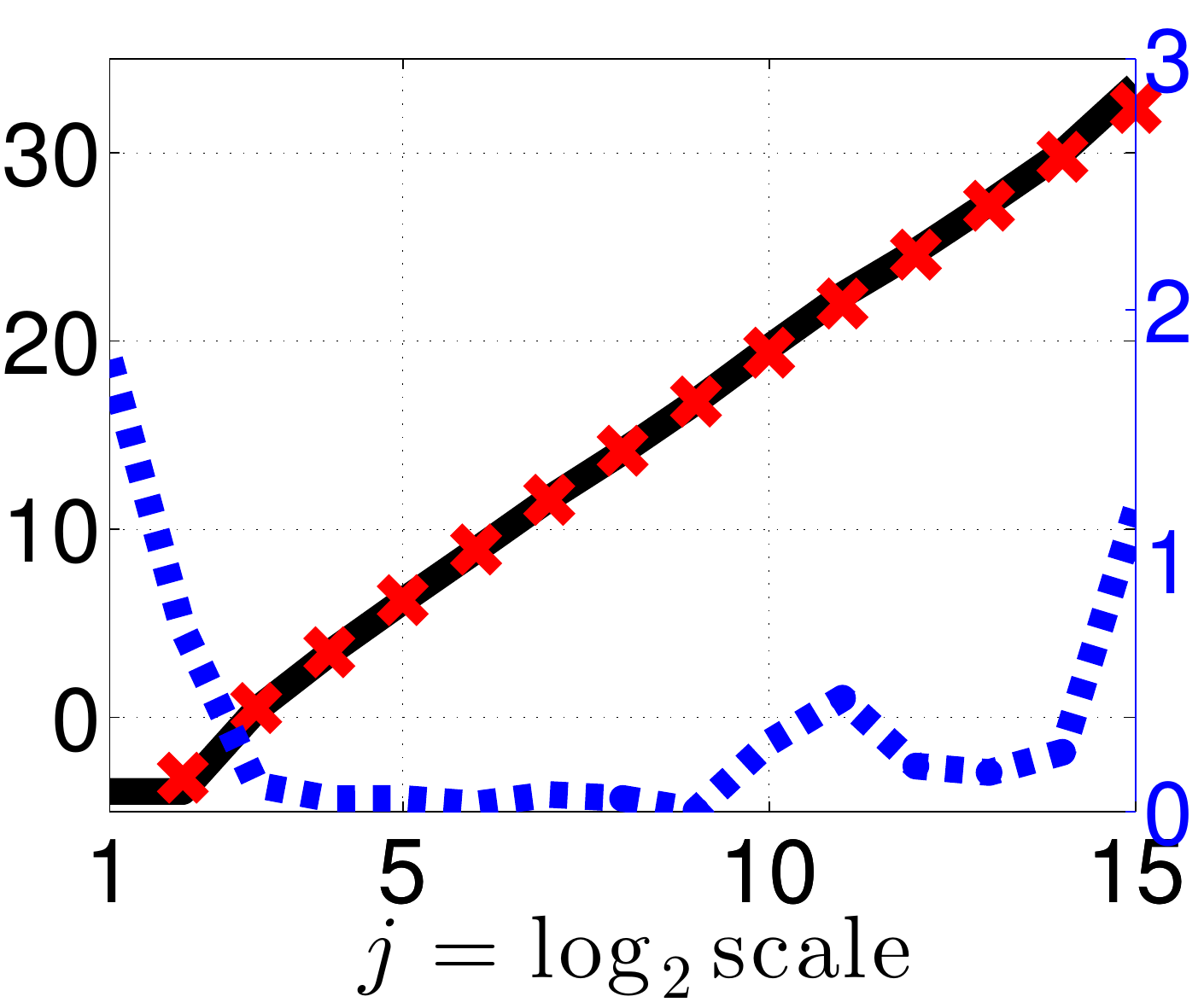}\qquad
\includegraphics[width=0.22\linewidth]{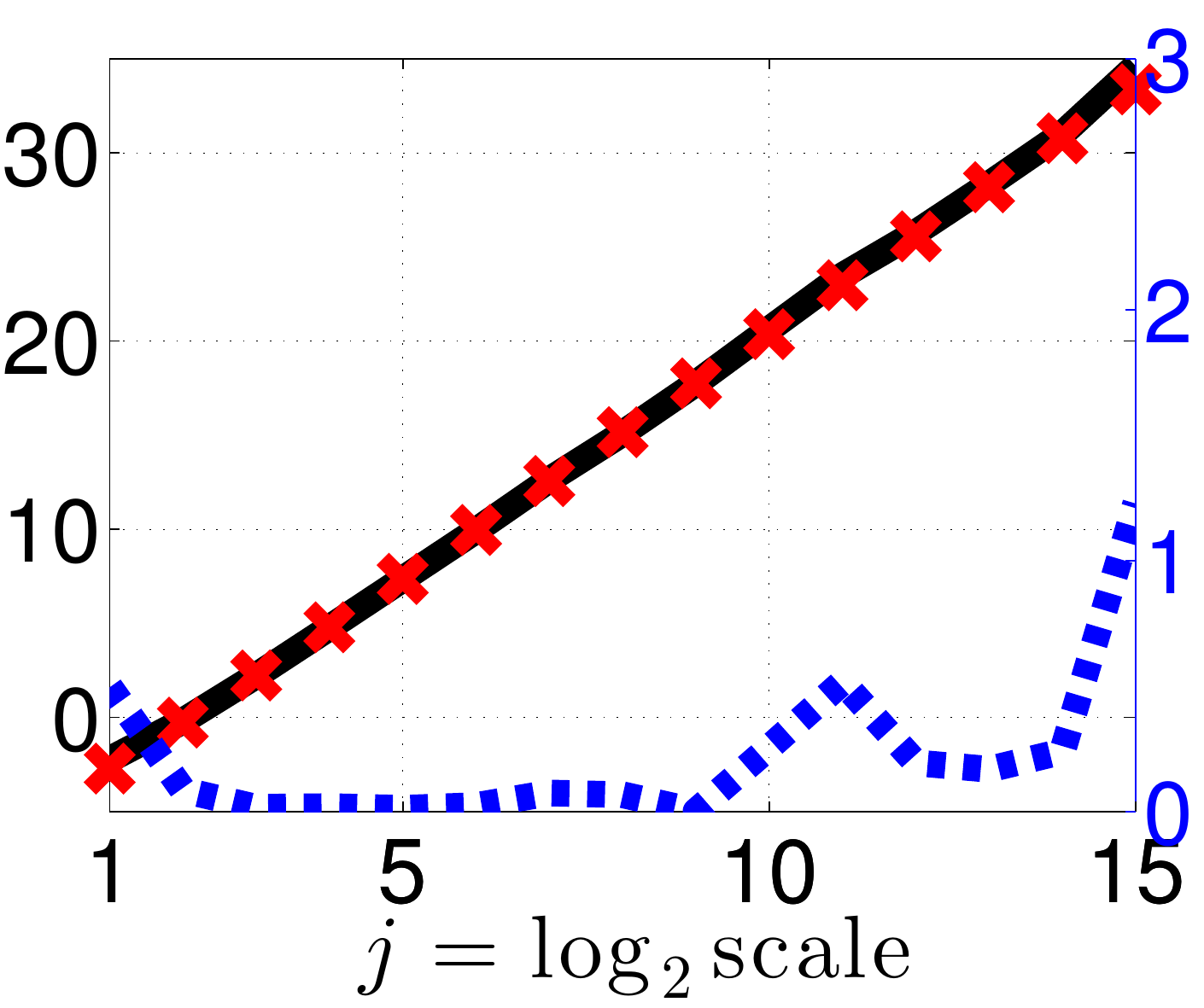}
}
\caption{\textbf{Wavelet Spectrum.} Superimposition of $\log_2 | E_{p,p'}(2^j,\Theta)|$ (red '+') and $\log_2 |S_{p,p'}(2^j)|$  (solid black line)  (with $(p,p') = (1, 1), (1,2), (2,2) $ from left to right) for a single realization of Biv-ofBm with $\Theta=(h_1=0.4, h_2=0.8, \rho_{\rm x}=0.1, \sigma_{\rm x_1}=1, \sigma_{\rm x_2}=1, \beta=0.5, \gamma=0.5)$ (absolute difference between data and model is shown in dashed blue). \label{fig:illustrationMixedOFBM}}
\end{figure}

The goal is to estimate the Biv-OfBm parameters $\Theta=(h_1,h_2,\rho_{\rm x},\sigma_{\rm x_1},\sigma_{\rm x_2},\beta,\gamma)$ starting from the wavelet spectrum $ \mathbb{E} D_{\rm y}(j,\cdot)D_{\rm y}(j,\cdot)^*$. The plug-in estimator of the ensemble variance $ \mathbb{E} D_{\rm y}(j,\cdot)D_{\rm y}(j,\cdot)^*$ is the sample variance
$$
S(2^j) = \frac{1}{K_j}\sum^{K_j}_{k=1} D_{\rm y}(j,k) D_{\rm y}(j,k)^*, \quad K_j = \frac{N}{2^j},
$$
where  $N$ denotes the sample size.
Fig.~\ref{fig:illustrationMixedOFBM} illustrates the fact that $S(2^j)$ is a satisfactory estimator for $ \mathbb{E} D_{\rm y}(j,\cdot)D_{\rm y}(j,\cdot)^*$.

\section{Non-linear regression based estimation and Branch and Bound algorithm}
\label{sec.BB}
\subsection{Identification procedure as a minimization problem}

The estimation of the parameter vector $\Theta $ of Biv-OfBm is challenging because its entry-wise wavelet (or Fourier) spectrum is a mixture of power laws (cf.\ Eqs.~(\ref{eq:E11})-(\ref{eq:E22})). This precludes the direct extension of classical univariate techniques, based on the scalar relation Eq.~\eqref{equ:wavcoef} \cite{Veitch1999}. For this reason, we formulate the full identification of Biv-OfBm (i.e., the estimation of $\Theta$) as a minimization problem\footnote{The superscript $\cdot^M$ has been added in order to refer to the $M$-estimator whose theoretical details are given in Section~\ref{sec.Co}.}:
\begin{equation}
\label{eq:globMinProb}
\hat{\Theta}^M_N =  \underset{\Theta\in\mathcal{Q}_0}{\mathrm{argmin}}\; C_N(\Theta), \makebox{ where }
\end{equation}
\begin{equation}
C_N(\Theta) \equiv{\hskip -.2cm} \sum_{\substack{i_1, i_2 =1\\ i_1\leq i_2}}^{P=2} \sum_{j= j_1}^{j_2}{\hskip -.2cm} \left( \log_2 |S_{i_1,i_2}(2^j)| - \log_2 | E_{i_1,i_2}(2^j,{\Theta})| \right)^2.
\label{eq:globMinProbb}
\end{equation}
The use of $ \log_2$ ensures that the scales $2^j$, $j = j_1,\ldots,j_2$ contribute equally to $C_N$. The search space incorporates prior information in the shape of constraints:
Sections \ref{sec:underD} and \ref{sec:funderD} impose $h_1\leq h_2$ and $\rho_{\rm x}\in[0,1]$~;
feasible solutions must satisfy the constraint $g(h_1,h_2,\rho_{\rm x})>0$ (cf. Eq.~\eqref{eq:condPositiveDefinite}) and $(\beta, \gamma) \in [-1,1]^2$.
For the sake of feasibility, we further restrict $(\sigma_{\rm x_1}, \sigma_{\rm x_2}) \in[0,\sigma_{\max}]^2$, with $\sigma_{\max}=\sqrt{\widehat{\sigma}^2_{\rm y_1} +\widehat{\sigma}^2_{\rm y_2}}$, where $\widehat{\sigma}^2_{\rm y_m} $ denotes the sample variance estimates of the increments of $Y_m$. We arrive at the parameter space
\begin{align}\label{e:Q0}
\mathcal{Q}_0 = \Big\{ \Theta=(h_1,h_2,&\rho_{\rm x},\sigma_{\rm x_1},\sigma_{\rm x_2},\beta,\gamma) \in \RR^7 \,\vert \,\Theta\in [0,1]^3 \times [0,\sigma_{\max}]^2 \times [-1,1]^2, \nonumber\\& g(h_1,h_2,\rho_{\rm x})>0, h_1 \leq h_2 \Big\}.
\end{align}

The minimization of $C_N(\Theta) $ is an intricate task for two reasons. First, because it involves disentangling a mixture of power laws, which yields a highly non-convex function.
Second, because the parameters to be estimated in $\Theta$ (scaling exponents, mixing coefficients, variances and correlation) are very different in nature.
The present contribution makes the original proposition of searching the global minimum of Eq.~(\ref{eq:globMinProb}) by means of a Branch \& Bound procedure detailed in the next section.
\subsection{Global minimization via a Branch \& Bound strategy}

Branch \& Bound algorithms consist of smart enumeration methods, which were shown to solve a variety of constrained global non-convex optimization problems \cite{Ichida_K_1979_computing_aiamgo,Kearfott_RB_1992_j-global-opt_aibbabcop,Moore_R_1992_book_rmgo,Ratschek_H_1988_book_ncmgo}.
In the context of the estimation problem \eqref{eq:globMinProb}, it amounts to partitioning (branching) the search space $\mathcal{Q}_0$ into smaller and smaller subregions, bounding the range of the objective function $C_N$ in each subregion, and then identifying the region containing the global minimum. This can be rephrased as 4 steps which are repeated until a stopping criterion is reached:\\
\noindent \textit{- Selecting:} Choose any region $\mathcal{R}$ from the search space and relax it into a closed convex set, i.e. an \emph{interval}, as illustrated by the dashed line in Fig.~\ref{fig:innerConvex} (left plot).\\
\noindent\textit{- Partitioning:} Divide $\mathcal{R}$ into two smaller regions $\mathcal{R}_a$ and $\mathcal{R}_b$.\\
\noindent\textit{- Bounding:} Compute lower and upper bounds of $C_N$ on $\mathcal{R}_a$ and $\mathcal{R}_b$.
Upper bounds can be obtained by evaluating $C_N$ anywhere in the region at hand.
Lower bounds  are computed by resorting to \emph{interval arithmetic} techniques (cf. Appendix~\ref{app:IA} and \cite{Moore_RE_1966_book_ia,Jaulin_J_2001_book_aia, Moore_R_2009_book_iia}), which combine elementary operations to produce rough lower bounds for the range of a given function, here $C_N$, on any interval. \\
\noindent\textit{- Pruning: } Pruning is driven by three mechanisms: discard regions that do not satisfy constraints (\textit{infeasibility})~;
discard regions whose lower bound is larger than the smallest upper bound as they cannot contain the global minimum (\textit{bound})~;
discard regions whose size (for all parameters) has reached the targeted precision (\textit{size}).

\subsection{Branch \& Bound procedure for Biv-OfBm identification}

\subsubsection{Convex relaxation}

\begin{figure}[t]
\centerline{
\includegraphics[scale=.50]{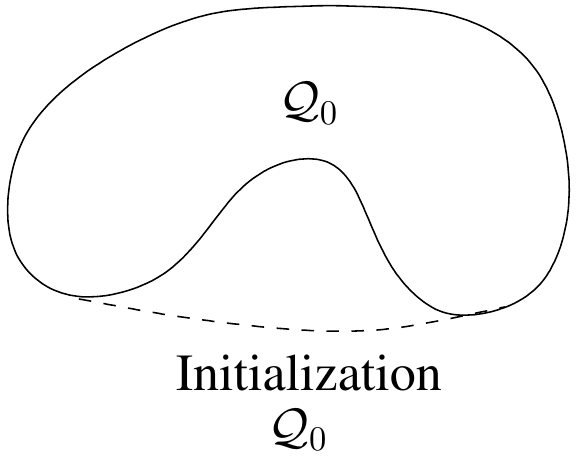} \hspace{2mm}
\includegraphics[scale=.50]{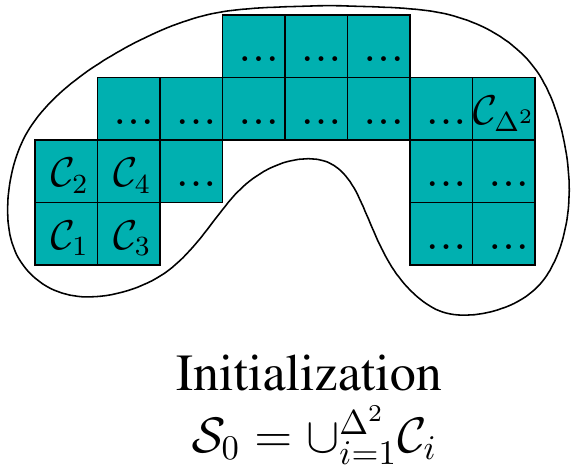}
}
\caption{\textbf{Schematic view of the proposed approximation $\mathcal{S}_0\subseteq\mathcal{Q}_0$.}\label{fig:innerConvex}}
\end{figure}

By nature, \emph{interval arithmetic} techniques apply only to \emph{intervals}, i.e., to convex sets.
Therefore, in most Branch \& Bound procedures, a \emph{convex relaxation} of the search space $\mathcal{Q}_0$ is required at initialization, as sketched in Fig.~\ref{fig:innerConvex} (left). However, in the present case, a convex relaxation of $\mathcal{Q}_0$ is not feasible because of the constraint  $g(h_1,h_2,\rho_{\rm x})>0$.

Instead, we propose to approximate $\mathcal{Q}_0$ by an inner convex relaxation $\mathcal{S}_0$, consisting of the union of $\Delta^2$ separable convex sets $\mathcal{C}_i$, i.e., $ \left({\mathcal{S}}_0=\cup_{i=1}^{\Delta^2} \mathcal{C}_i\right)\subset\mathcal{Q}_0$, as illustrated in Fig.~\ref{fig:innerConvex} (right).
In practice, $(h_1,h_2,\rho_{\rm x})\in[0,1]^3$ is approximated by a union of $\Delta^2$ non-overlapping parallelepipedic sets $\{\mathcal{P}_i\}_{1\leq i\leq \Delta^2}$, denoted $\mathcal{P} \subset [0,1]^3$.
They are obtained by dividing $(h_1,h_2)\in[0,1]^2$ into squares $\mathcal{T}_i$ with a discretization step $\Delta^{-1}$ and defining $\mathcal{P}_i=\mathcal{T}_i\times[0,\rho_i]$ where $\rho_i$ is the largest value such that $ (\forall(h_1,h_2)\in\mathcal{T}_i), \quad  g(h_1,h_2,\rho_i)>0$:
\begin{equation}
\mathcal{C}_i = \Big\{ \Theta=(h_1,h_2,\rho_{\rm x},\sigma_{\rm x_1},\sigma_{\rm x_2},\beta,\gamma) \in \RR^7 \,\vert \, \Theta\in \mathcal{P}_i \times [0,\sigma_{\max}]^2 \times [-1,1]^2, h_1 \leq h_2 \Big\}.
\label{eq:mathcalC_i}
\end{equation}
In this inner convex relaxation strategy, the constraint $g(h_1,h_2,\rho_{\rm x})>0$ is necessarily satisfied, a major practical benefit as infeasible regions need not be explored.

\subsubsection{Algorithm}
\label{sec:algoBB}

The full identification of Biv-OfBm is achieved via the following proposed  sequence of operations:

\noindent\textbf{Inputs.} \\
\noindent  - From data, compute the wavelet spectrum $S(2^j)$, $j = j_1, \hdots,j_2$~;\\
\noindent  - Pick the Biv-OfBm model (i.e., Eqs.~(\ref{eq:E11})-(\ref{eq:E22}))~;\\
\noindent  - Set the precision $\delta$ for each parameter~;\\
\noindent  - Set the inner convex relaxation parameter $\Delta$ in order to approximate the set $\mathcal{Q}_0$ by ${\mathcal{S}}_0=\cup_{i=1}^{\Delta^2}\mathcal{C}_i$~;\\
\noindent\textbf{Initialization.} \\
\noindent  - Set  $\widehat{\mathcal{S}}=\emptyset$ and $k=0$.\\
\noindent  - Compute lower bounds $l_i$ of $C_N$ on $\mathcal{C}_i$ ($\forall i=1,\ldots,\Delta^2$)~;\\
\noindent  - Compute upper bounds $u_i$ of $C_N$ on $\mathcal{C}_i$ ($\forall i=1,\ldots,\Delta^2$)~; \\
\noindent  -  Set $U=\min(u_1,\ldots,u_{\Delta^2})$ and $L=\min(l_1,\ldots,l_{\Delta^2})$.\\
\noindent\textbf{Iteration.} Let $\mathcal{S}_k$ denote the partitioning at the step $k$.\\
\noindent  - {Selecting:} Select region $\mathcal{R} \subset \mathcal{S}_k$ with lowest lower bound $L$ (best-first-search strategy)~\footnote{At the step $k=0$, this amounts in choosing one $\mathcal{C}_i$. More generally, it consists in selecting one element of the partitioning $\mathcal{S}_k$.}.\\
\noindent  - {Cutting:} Divide $\mathcal{R}$ into $\mathcal{R}_a$ and $\mathcal{R}_b$ such that $\mathcal{R} = \mathcal{R}_a\cup\mathcal{R}_b$ and $\mathcal{R}_a\cap\mathcal{R}_b = \emptyset$, along its longest edge, in half, where
the length of an edge is defined relatively to the maximum
accuracy  $\delta$ prescribed by the practitioner. \\
\noindent  - {Lower bound:} Compute lower bounds of $\mathcal{R}_a$ and $\mathcal{R}_b$, using interval arithmetic.  \\
\noindent  - {Upper bound:} Compute upper bound of $\mathcal{R}_a$ (resp. $\mathcal{R}_b$), by evaluating $C_N(\Theta)$ for $\Theta$ chosen at the center of $\mathcal{R}_a$  \mbox{(resp. $\mathcal{R}_b$).} \\
\noindent  - {Branching:} Update the partitioning ${\mathcal{S}}_{k+1} = ({\mathcal{S}}_k\backslash\mathcal{R})\cup\mathcal{R}_a\cup\mathcal{R}_b,$
and update $U$ and $L$ on this new partition.\\
\noindent  - {Pruning:} Discard regions $\mathcal{R^*}$ of ${\mathcal{S}}_{k+1} $ either by \textit{bound}, \textit{infeasibility} or \textit{size}, i.e.
${\mathcal{S}}_{k+1} \leftarrow {\mathcal{S}}_{k+1} \backslash\mathcal{R}^*.$
Append to $\widehat{\mathcal{S}}$ the regions of ${\mathcal{S}}_{k+1}$ discarded by \textit{size}~; Discard regions in $\widehat{\mathcal{S}}$ for \textit{bound}. \\
\noindent  - Set $k\leftarrow k + 1$\\
\noindent\textbf{Stop and Output.} Stop iterations when ${\mathcal{S}}_k$ is empty.
Output $\widehat{\mathcal{S}}$ as the list of potential solutions at targeted precision $\delta$.
Output region in $\widehat{\mathcal{S}}$ with lowest upper bound and corresponding  best estimate at the targeted precision
\begin{equation}\label{e:Theta-hat^(M,BB)}
\hat \Theta^{M, BB}_N.
\end{equation}

\section{Asymptotic theory in the multivariate setting}
\label{sec.Co}

In this section, the asymptotic properties of the exact solution $\widehat{\Theta}^{M}_N$ (see \eqref{eq:globMinProb}) are studied theoretically for a general multivariate OfBm \cite{Amblard_P-0_2011_j-ieee-tsp_imfbm,didier:pipiras:2012,Didier_G_2011_bernouilli_irpofbm}. In other words, the results encompass, but are not restricted to, the bivariate framework of Section \ref{sec.OFBM}.
Regarding the OfBm $\underline{B}_{W,{\underline{H}}} $, it is assumed that:
(OFBM1) $\underline{B}_{W,{\underline{H}}}$ is an $\bbR^P$-valued OfBm with Hurst parameter $\underline{\underline{H}}$, not necessarily diagonalizable, where the eigenvalues of $\underline{\underline{H}}$ satisfy $0 < \Re(h_k) < 1$, $k = 1,\hdots,n$~;
(OFBM2) $\E \underline{B}_{W,\underline{H}}(t)\underline{B}_{W,\underline{H}}(t)^*$, $t \neq 0$, is a full rank matrix (properness)~;
(OFBM3) $\underline{B}_{W,{\underline{H}}}$ is a time reversible stochastic process.
Regarding $\psi_0 \in L^1(\bbR)$, it is assumed that:
(W1) $N_{\psi}\geq 2$~;
(W2) $\textnormal{supp}(\psi_0)$ is a compact interval~;
(W3) $\sup_{t \in \bbR} |{\psi}_0(t)| (1 + |t|)^{\alpha} < \infty$ for $\alpha > 1$.
The summation range in the objective function \eqref{eq:globMinProbb} is generalized to  $i_1,i_2 = 1,\hdots,P\geq 2$.
The proofs of the statements in this section can be found in Appendix \ref{app:AA}.

\subsection{The asymptotic normality of the wavelet spectrum}
\label{s:asymptotic_normality}

The asymptotic behavior of the estimator $\widehat{\Theta}^{M}_N$ draws upon the asymptotic normality of the wavelet variance for fixed scales. Under the assumptions (OFBM1-3) and (W1-3), the latter property can be established by an argument that is almost identical to that in \cite[Theorem 3.1]{AbryDidier2015}. For the reader's convenience, we reproduce the claim here.

\begin{theorem}\label{t:asymptotic_normality_wavecoef_fixed_scales}
Suppose the assumptions (OFBM1-3) and (W1-3) hold. Let $F \in {\mathcal S}(\frac{P(P+1)}{2}m,\bbR)$\footnote{$\mathcal{S}(n,\mathbb{R})$ is the space of real symmetric matrices of size $n \times n$.} be the asymptotic covariance matrix described in \cite[Proposition 3.3]{AbryDidier2015}. Then,
\begin{equation}\label{e:asymptotic_normality_wavecoef_fixed_scales}
\Big(\sqrt{K_j}(\textnormal{vec}_{{\mathcal S}} S(2^j)- \textnormal{vec}_{{\mathcal S}} E(2^j,\Theta) ) \Big)_{j = j_1 , \hdots, j_2} \stackrel{d}\rightarrow {\mathbf Z},
\end{equation}
as $N \rightarrow \infty$, where $j_1 < \hdots < j_2$ and $m = j_2 - j_1 + 1$ and ${\mathbf Z} \stackrel{d}= {\mathcal N}_{\frac{P(P+1)}{2} \times m}(0,F)$,
where $\textnormal{vec}_{{\mathcal S}} $ defines the operator that vectorizes the upper triangular entries of a symmetric matrix:
\begin{equation*}\label{e:vec_S}
\textnormal{vec}_{{\mathcal S}}(S) = (s_{11}, \hdots, s_{1P}; \hdots;s_{P-1,P-1}, s_{P-1,P};s_{P,P})^*.
\end{equation*}
\end{theorem}

\subsection{Consistency of $\hat \Theta^M_N$}

Let $ \Theta_0 $ be the true parameter value. To prove the consistency of $\hat \Theta^M_N$, the following additional assumptions on the parametrization $\Theta$ are required: \\
\noindent i)  The parameter space $\Xi \subseteq {\mathcal Q}_0$ (see \eqref{e:Q0}) is a finite-dimensional compact set and
\begin{equation}\label{e:parametrization_Theta0_intXi}
\Theta_0 \in \textnormal{int}\hspace{0.5mm}\Xi ;
\end{equation}
\noindent ii) For some $j^*=j_1,\hdots,j_2$,
\begin{equation}\label{e:EW(2j)_is_identifiable}
\Theta \neq \Theta' \Rightarrow  |E_{i^*_1,i^*_2}(2^{j^*},\Theta)| \neq | E_{i^*_1,i^*_2}(2^{j^*},\Theta')|
\end{equation}
for some matrix entry $(i^*_1,i^*_2)$;\\
\noindent iii) $ \forall  i_1,i_2=1,\hdots,P, \quad j=j_1,\hdots,j_2$,
\begin{equation}\label{e:EW(2j,Theta)_entry_not_zero}
E_{i_1,i_2}(2^j,\Theta_0) \neq 0~;
\end{equation}
\begin{equation}\label{e:smoothness}
\makebox{\hskip -5.75cm iv) The mapping }\Theta \mapsto E(2^j,\Theta) \quad \quad \quad \quad \quad \quad\quad\quad\quad\;
\end{equation}
is three times continuously differentiable on $\textnormal{int}\hspace{0.5mm}\Xi$.\\

Under \eqref{e:EW(2j,Theta)_entry_not_zero}, the functions $\log_2 | E_{i_1,i_2} (2^{j},\Theta)|$, $i_1,i_2 = 1,\hdots, P$, are well-defined.
This fact and Theorem \ref{t:asymptotic_normality_wavecoef_fixed_scales} then imply that the functions $\log_2 |S_{i_1,i_2}(2^{j})|$ are well-defined with probability going to 1.
In turn, condition \eqref{e:EW(2j)_is_identifiable} implies that the (entry-wise) absolute value of the target matrix $E(2^{j^*},\Theta)$ is (parametrically) identifiable, namely, there is an injective function $\Xi \ni \Theta \mapsto |E(2^{j^*},\Theta)|$.

The objective function $C_N(\Theta)$ is a function of $N$, and so is $S(2^j)$, $j = j_1,\hdots,j_2$. Since $C_N(\cdot)$ is continuous and $\Xi$ is compact, then for all $N$ a minimum $\widehat{\Theta}_N$ is attained (a.s.), whence we can form one such sequence
\begin{equation}\label{e:minima}
\{\widehat{\Theta}^M_N\}_{N \in \bbN}.
\end{equation}
Any sequence \eqref{e:minima} defines an $M$-estimator of $\Theta_0$, e.g., \cite[chapter 5]{vandervaart:2000}.
The next theorem shows that \eqref{e:minima} is consistent.

\begin{theorem}\label{t:consistency}
Under the assumptions of Theorem \ref{t:asymptotic_normality_wavecoef_fixed_scales}, suppose in addition that the conditions i) to iv) hold. Then, the sequence of minima \eqref{e:minima} is consistent for $\Theta$, namely,
\begin{equation}\label{e:consistency}
\widehat{\Theta}^M_{N} \rightarrow \Theta_0 \quad \textnormal{in probability}.
\end{equation}
\end{theorem}

\begin{remark}
The uniqueness of $\widehat{\Theta}^M_{N}$ for a given $N$ is not ensured by the conditions $i)$ to $iv)$, but  it is not needed in Theorem \ref{t:consistency}.
\end{remark}

\subsection{Asymptotic normality of $\hat \Theta^M_N$}

By comparison to consistency, showing asymptotic normality will require an additional assumption, laid out next.\\
\begin{equation}\label{e:parametrization_Theta0_intXi_deriv}
\noindent \makebox{\hskip 0cm v) } \small\det\Big(\sum^{j_2}_{j=j_1} \sum_{1 \leq i_1 \leq i_2 \leq n}\Lambda_{i_1,i_2}(2^j,\Theta_0)\Lambda_{i_1,i_2}(2^j,\Theta_0)^*\Big)> 0,\qquad\qquad\qquad\qquad\qquad\qquad\qquad\qquad\qquad\quad
\end{equation}
where we define the score-like vector $\Lambda_{i_1,i_2}(2^j,\Theta)^* = \nabla_{\Theta}\log_2  |E_{i_1,i_2}(2^j,\Theta)|$.

\begin{theorem}\label{t:asymptotic_normality}
Under the assumptions of Theorem \ref{t:asymptotic_normality_wavecoef_fixed_scales}, suppose in addition that the condition v) holds.
Let $\{\widehat{\Theta}^{M}_{N}\}_{N \in \bbN}$ be a consistent sequence of minima of $\{C_{N}\}_{N \in \bbN}$, respectively.
Then,
\begin{equation}\label{e:asymptotic_normality}
\sqrt{N}(\widehat{\Theta}^M_{N} - \Theta_0) \stackrel{d}\rightarrow {\mathbf W} , \quad N \rightarrow \infty, \makebox{ where }
\end{equation}
$$
{\mathbf W} \stackrel{d}=  \Big(\sum^{j_2}_{j=j_1}\sum_{1 \leq i_1 \leq i_2 \leq P}\Lambda_{i_1,i_2}(2^j,\Theta)\Lambda_{i_1,i_2}(2^j,\Theta)^*\Big)^{-1}\Big( \sum^{j_2}_{j_1=1} \frac{2^{j/2}}{ \log 2} \sum_{1 \leq i_1 \leq i_2 \leq P} Z_{i_1,i_2}(2^j)\frac{\Lambda_{i_1,i_2}(2^j,\Theta)}{E_{i_1,i_2}(2^j,\Theta_0)} \Big),
$$
where ${\mathbf Z} = (Z_{i_1,i_2}(2^j))_{j}$ is a random vector whose distribution is obtained in the weak limit \eqref{e:asymptotic_normality_wavecoef_fixed_scales}.
\end{theorem}

The next result is a corollary to Theorems \ref{t:consistency} and \ref{t:asymptotic_normality}.
\begin{corollary}\label{c:M-BB_is_asympt_normal}
Under the assumptions of Theorem \ref{t:asymptotic_normality}, let $\{\widehat{\Theta}^M_{N}\}_{N \in \bbN}$ be a sequence of minima of the objective function \eqref{eq:globMinProbb}. Also, let $\{\widehat{\Theta}^{M,BB}_{N}\}_{N \in \bbN}$ be an estimator of the form \eqref{e:Theta-hat^(M,BB)} which satisfies
\begin{equation}\label{e:assumption_BB}
\|\widehat{\Theta}^M_{N}- \widehat{\Theta}^{M, BB}_{N}\| \leq \frac{C}{N^{1/2 + \varepsilon}} \quad \textnormal{a.s.}
\end{equation}
for constants $C, \varepsilon > 0$. Then,
$$
\sqrt{N}(\widehat{\Theta}^{M, BB}_{N} - \Theta_0) \stackrel{d}\rightarrow {\mathbf W}, \quad N \rightarrow \infty,
$$
where the random vector ${\mathbf W}$ is given in \eqref{e:asymptotic_normality}.
\end{corollary}

\begin{remark}
The condition \eqref{e:assumption_BB} is easily satisfied in practice, since over a compact set and at a low computational cost a Branch and Bound algorithm is guaranteed to yield a solution which lies at a controllable distance of the true minimum.
\end{remark}

\begin{remark}
\label{rem:technical}
The technical condition \eqref{e:parametrization_Theta0_intXi_deriv} should be satisfied in many cases of interest, as discussed in Appendix \ref{app:AC}.
\end{remark}

\section{Estimation performance: empirical study}
\label{sec.sim}

  \newcommand{\lscale}{.165}
 \newcommand{\trleftbord}{.3}
 \newcommand{\trbotbord}{.2}
 \newcommand{\trleft}{2.5}
 \newcommand{\trbot}{1}
 \newcommand{\trright}{0}
 \newcommand{\trtop}{.2cm}
 \subsection{Numerical simulation setting}

Monte Carlo experiments were performed to empirically quantify the finite size performance of the estimator $\hat \Theta^{M, BB}_N$.
To examine the influence of $\Theta$ on the estimation performance, 9 different values of $\Theta $ were used, obtained essentially by varying the strength of the correlation amongst components ($\rho_{\rm x}=0.1$, $0.45$ and $0.8$) and of the mixing factor (no mixing, $\beta=\gamma=0$~; orthogonal mixing, $\beta=\gamma=0.5$, and $\beta=-\gamma=0.5$, referred to as anti-orthogonal). Three different sample sizes (short, $N=2^{10}$, medium, $N=2^{14}$ and large, $N=2^{18}$) were investigated. Results are reported here for $(h_1, h_2) = (0.4, 0.8)$.
Equivalent conclusions are obtained for other choices of $(h_1, h_2) $. For each set of parameters $\Theta$, the estimation performance were assessed by means of box plots computed from 100 independent copies. The synthesis of OfBm is achieved by using the multivariate toolbox devised in \cite{Helgason_H_2011_j-sp_fessmgtsuce,Helgason_H_2011_j-sp_smsspmdccme}, cf. \href{www.hermir.org}{www.hermir.org}.
The computational loads are also quantified as a percentage of the number of iterations that would be required by a systematic greedy grid search. The wavelet analysis was based on least asymmetric orthonormal Daubechies wavelets \cite{Mallat_S_2008_book_awtspsw}. All available scales were used to compute $C_N$:  $j_1=1 \leq j \leq j_2=\log_2 N - N_\psi-1$.
The results are reported for $N_\psi = 2$; it was found that further increasing $N_\psi$ did not improve the performance. The proposed Branch \& Bound procedure was run with $\mathcal{S}_0=\cup_{i=1}^{\Delta^2}$ for $\delta=50$. The impact of varying the requested precision $\delta $ was also investigated.

\begin{figure}[t]
     \centering
{\scriptsize \hskip .9cm$\rho_{\rm x}=0.1$\hskip 1.7cm$\rho_{\rm x}=0.45$  \hskip 1.4cm$\rho_{\rm x}=0.8$} \\
\rotatebox{90}{\hskip .2cm\scriptsize$\beta=\gamma=0$} \includegraphics[scale=\lscale, clip=true, trim=\trleftbord cm \trbot cm \trright cm \trtop cm]{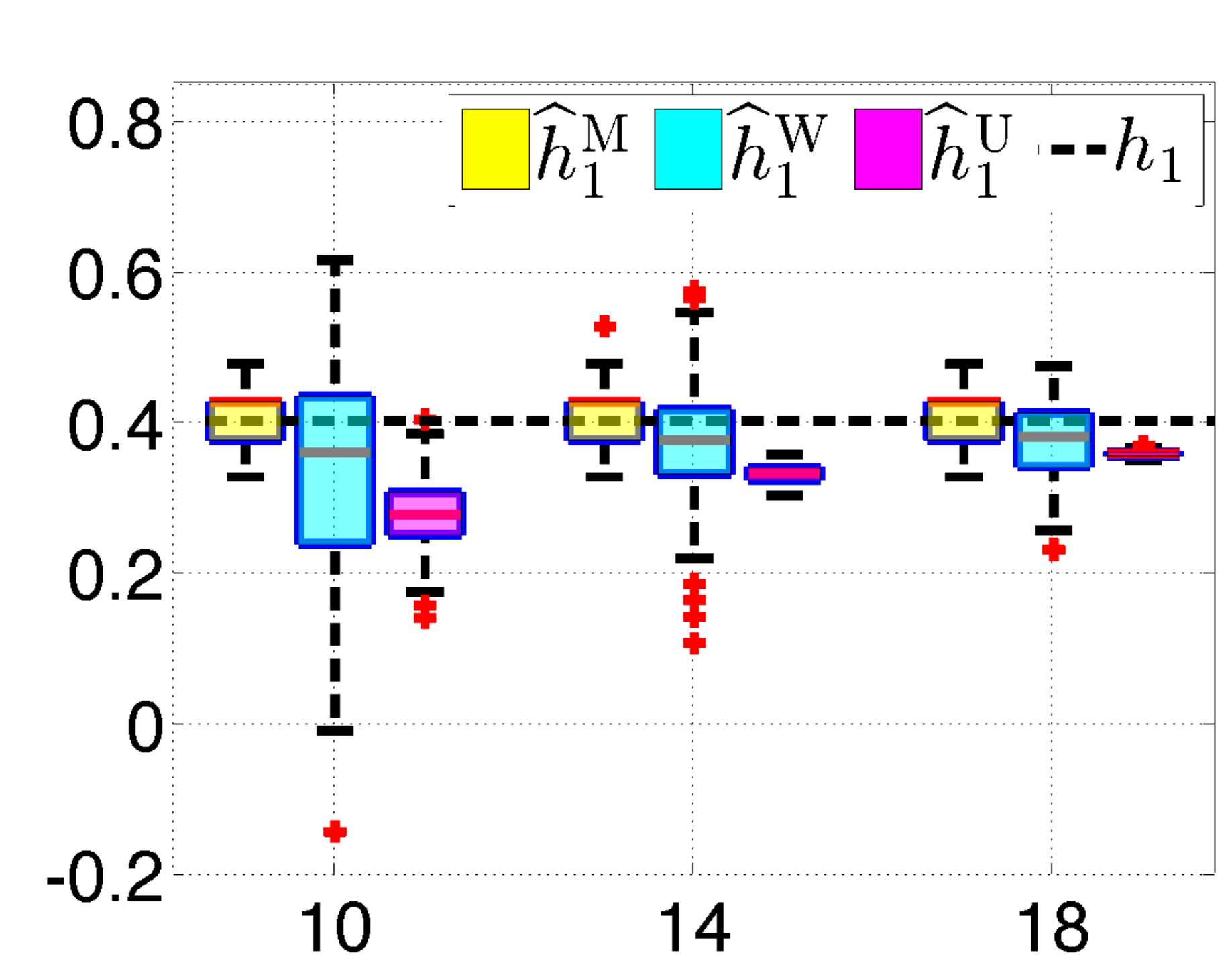}
  \includegraphics[scale=\lscale, clip=true, trim=\trleft cm \trbot cm \trright cm \trtop cm]{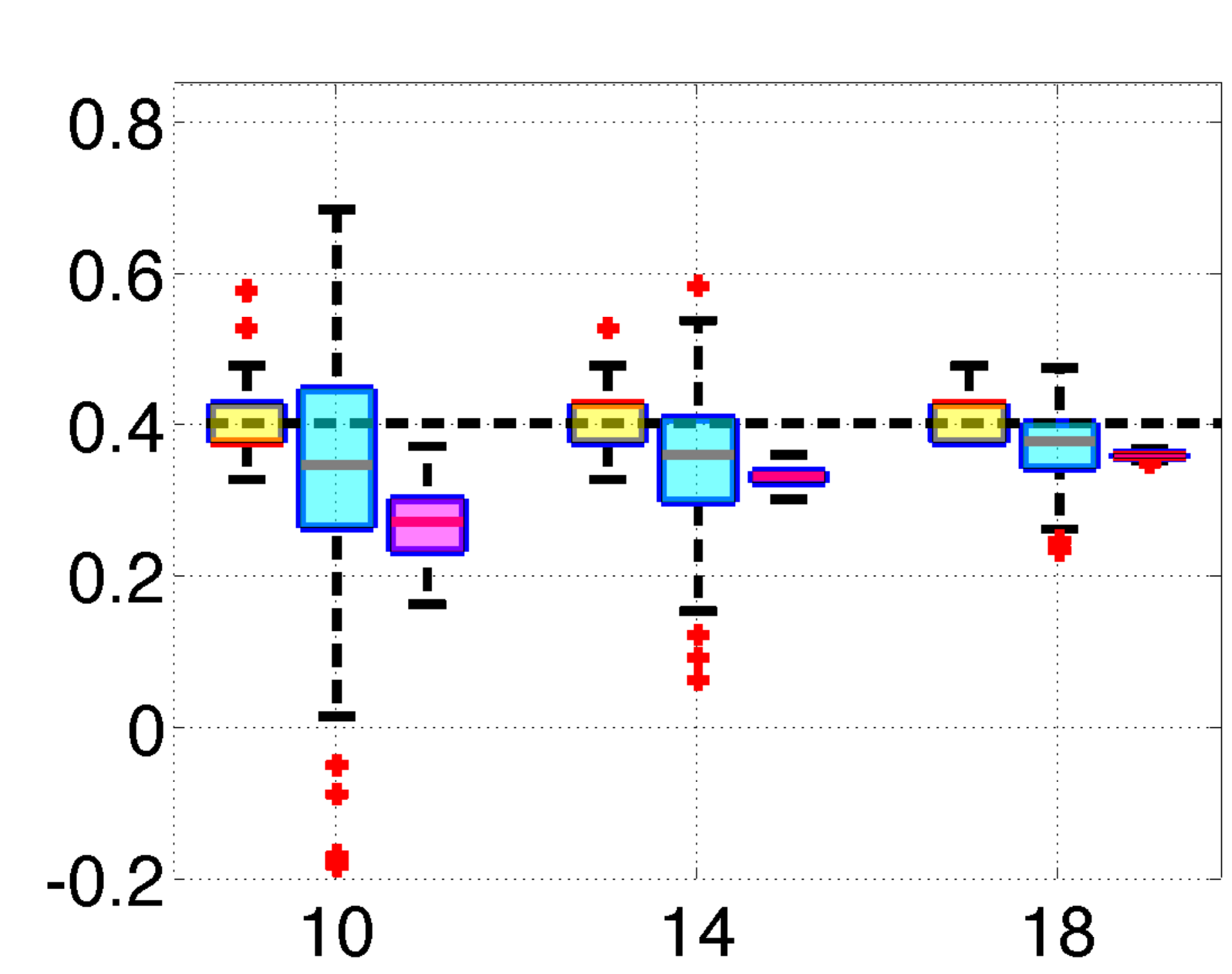}
 \includegraphics[scale=\lscale, clip=true, trim=\trleft cm \trbot cm \trright cm \trtop cm]{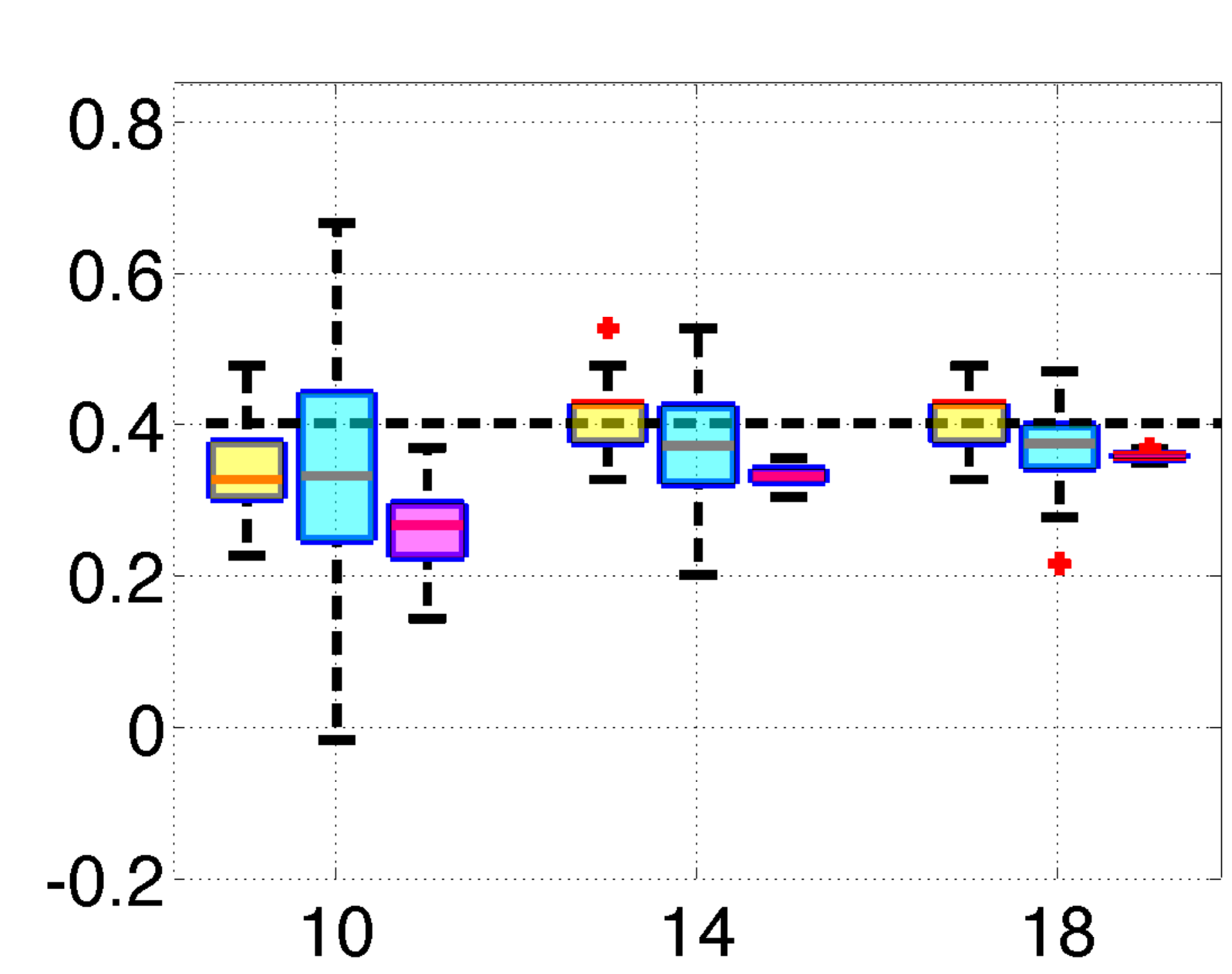}\\ \vskip -.3cm
 \rotatebox{90}{\hskip .2cm\scriptsize$\beta=\gamma=0.5$} \includegraphics[scale=\lscale, clip=true, trim=\trleftbord cm \trbot cm \trright cm \trtop cm]{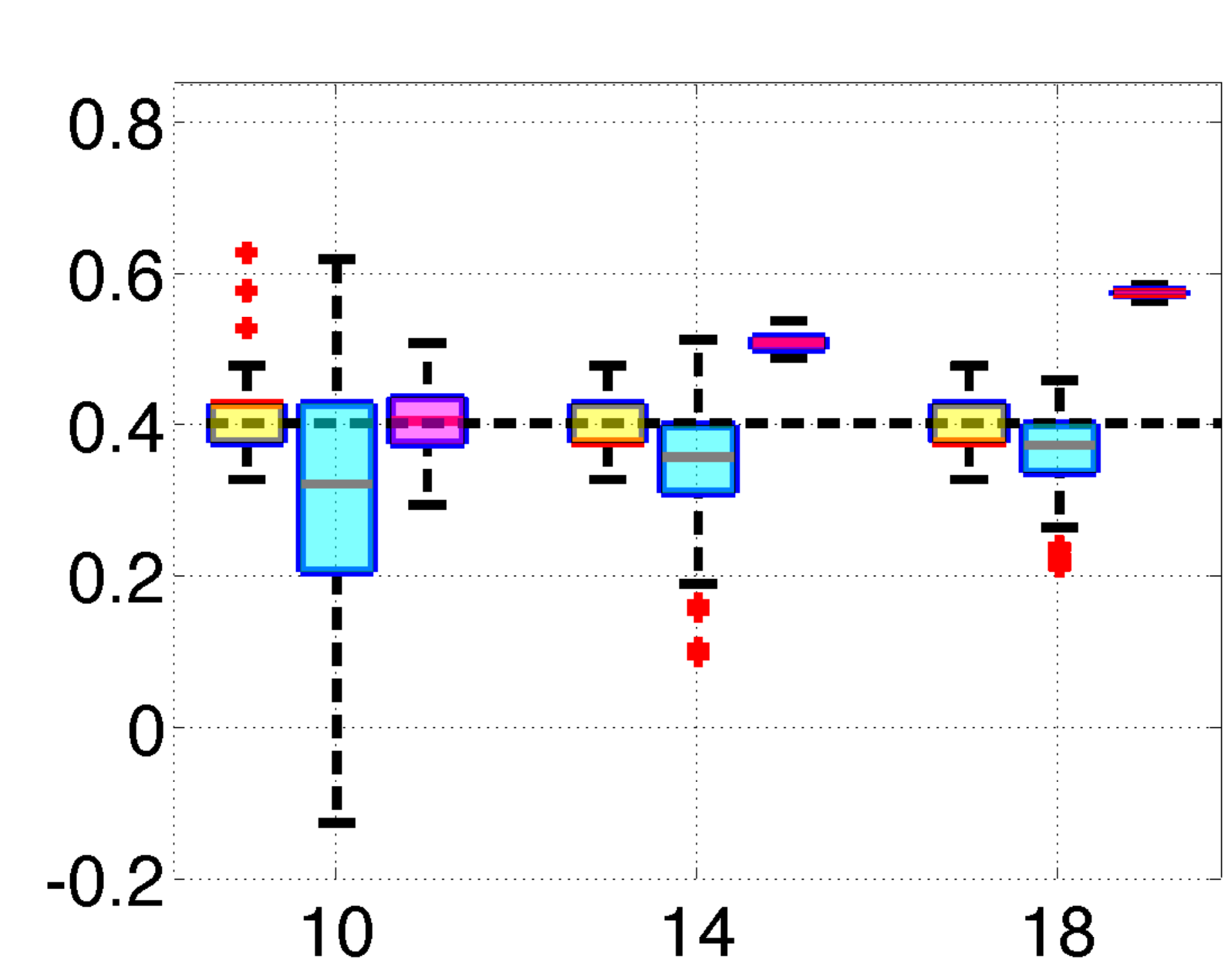}
  \includegraphics[scale=\lscale, clip=true, trim=\trleft cm \trbot cm \trright cm \trtop cm]{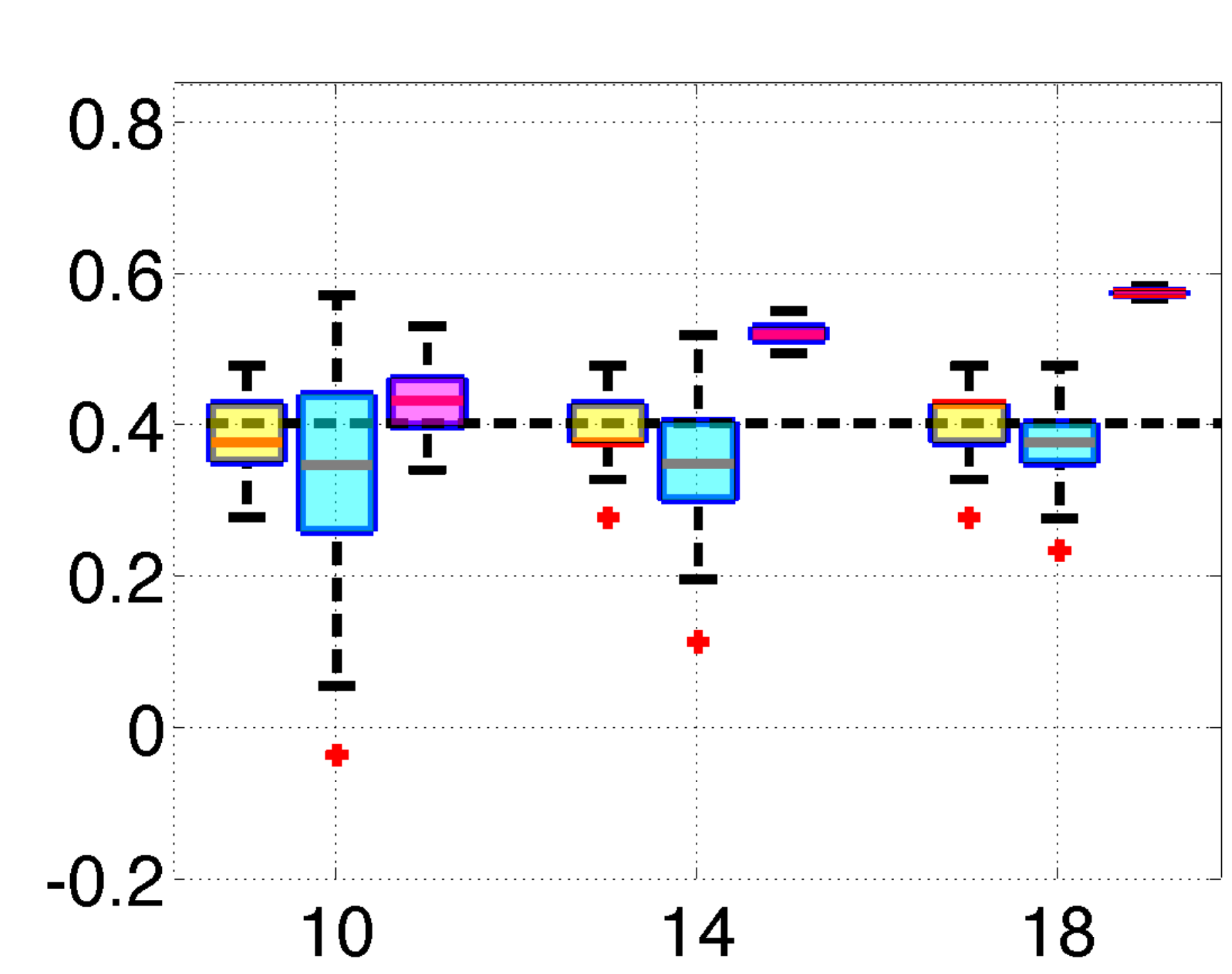}
 \includegraphics[scale=\lscale, clip=true, trim=\trleft cm \trbot cm \trright cm \trtop cm]{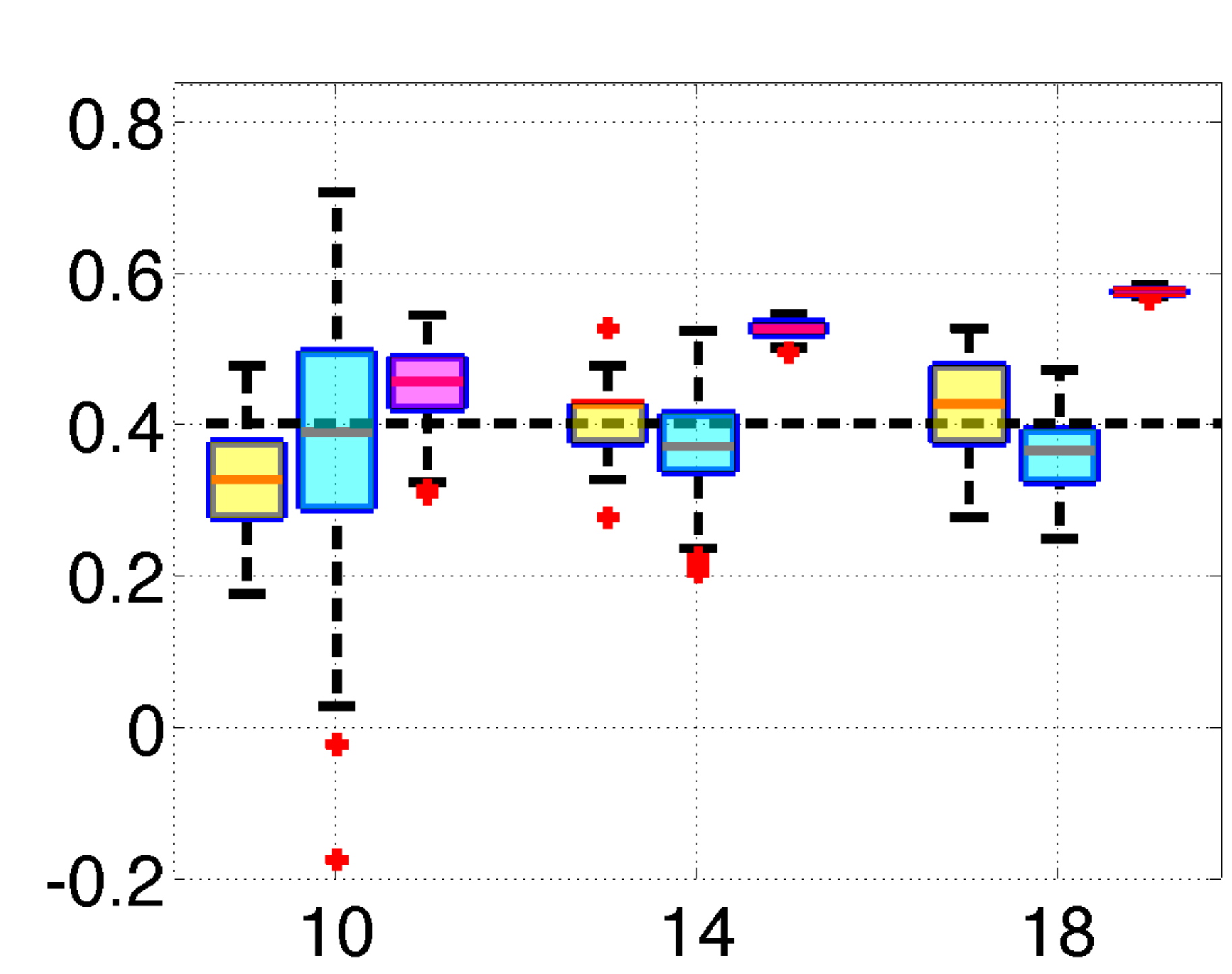}\\ \vskip -.3cm
 \rotatebox{90}{\hskip .1cm\scriptsize$\beta=-\gamma=0.5$} \includegraphics[scale=\lscale, clip=true, trim=\trleftbord cm \trbotbord cm \trright cm \trtop cm]{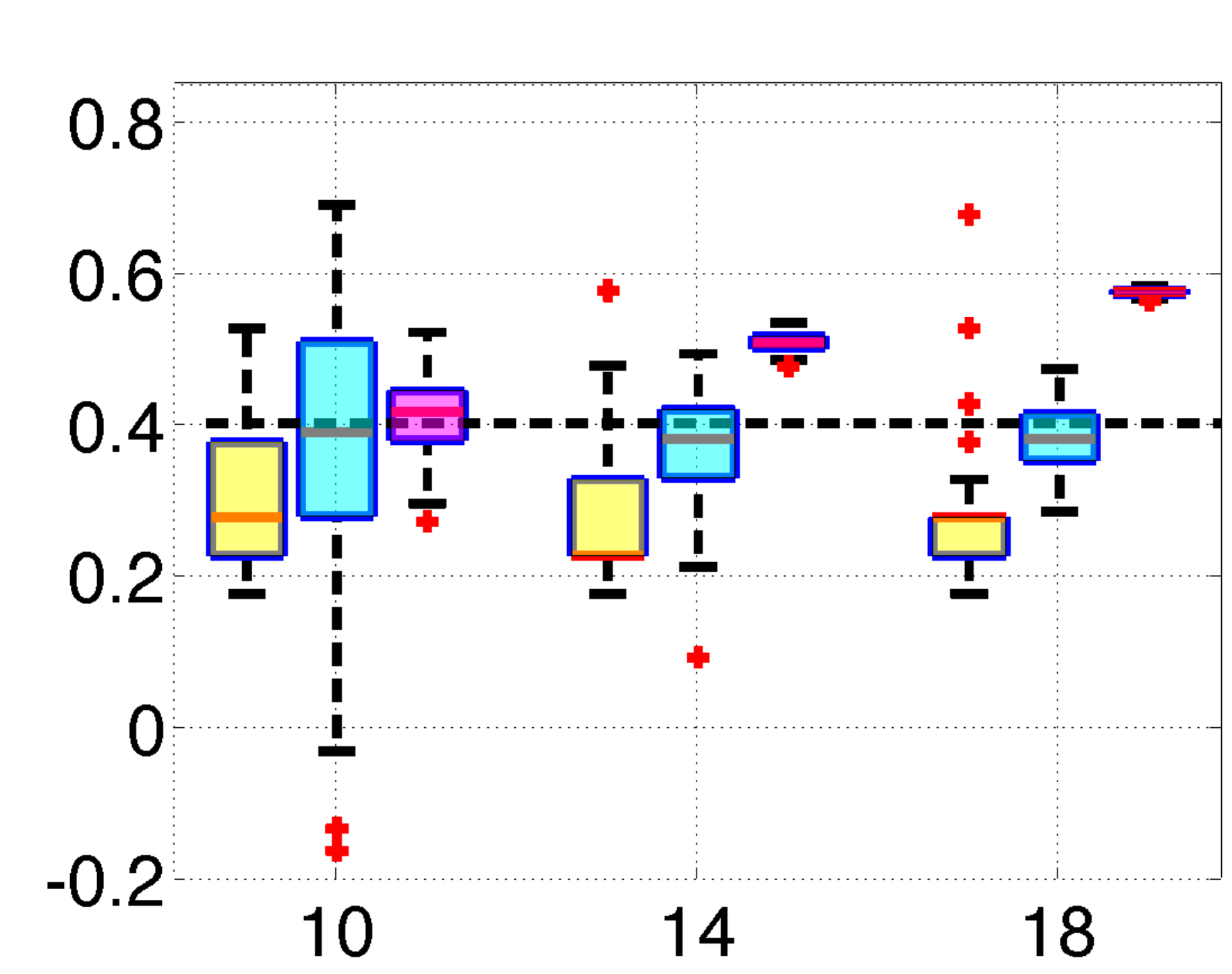}
  \includegraphics[scale=\lscale, clip=true, trim=\trleft cm \trbotbord cm \trright cm \trtop cm]{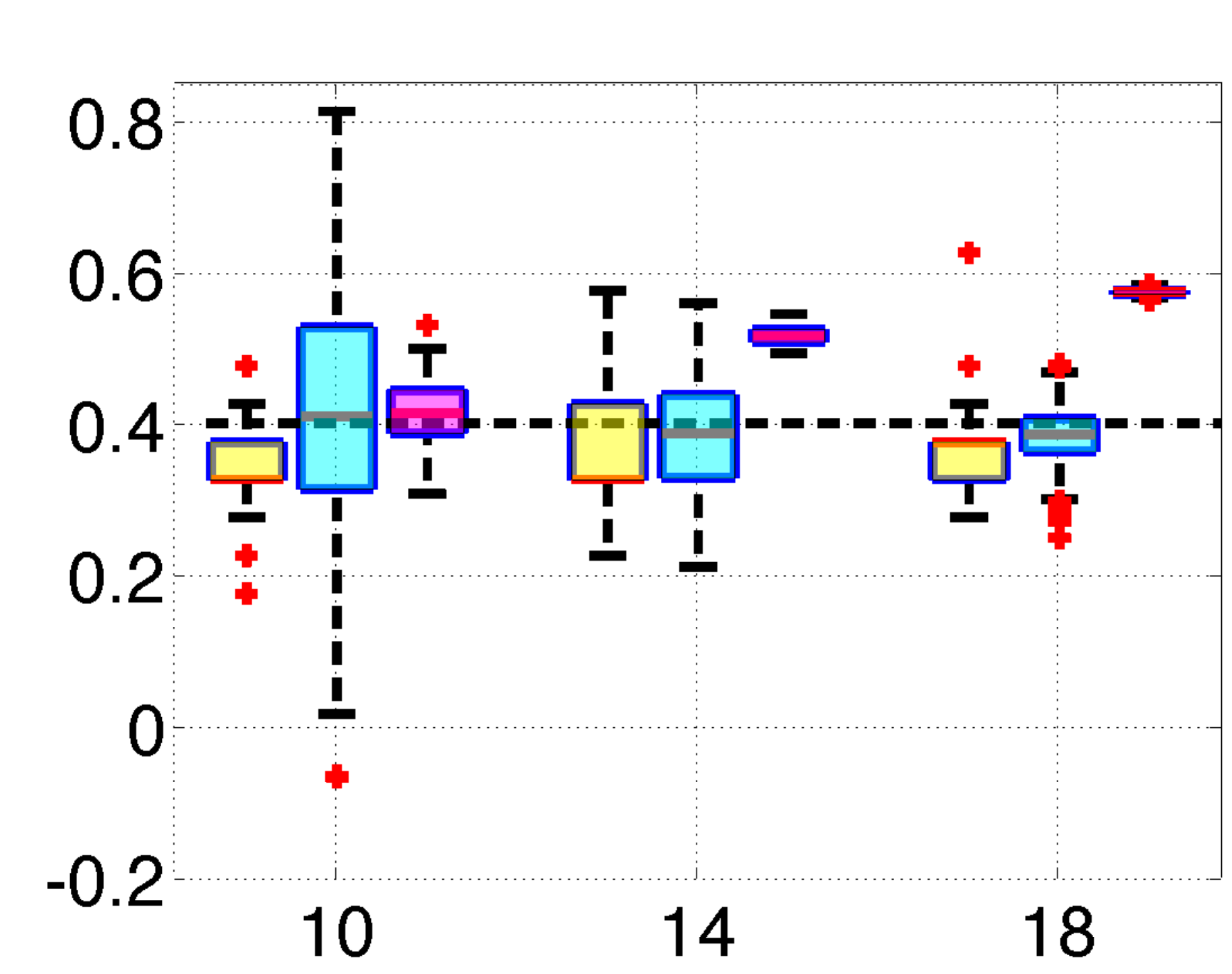}
 \includegraphics[scale=\lscale, clip=true, trim=\trleft cm \trbotbord cm \trright cm \trtop cm]{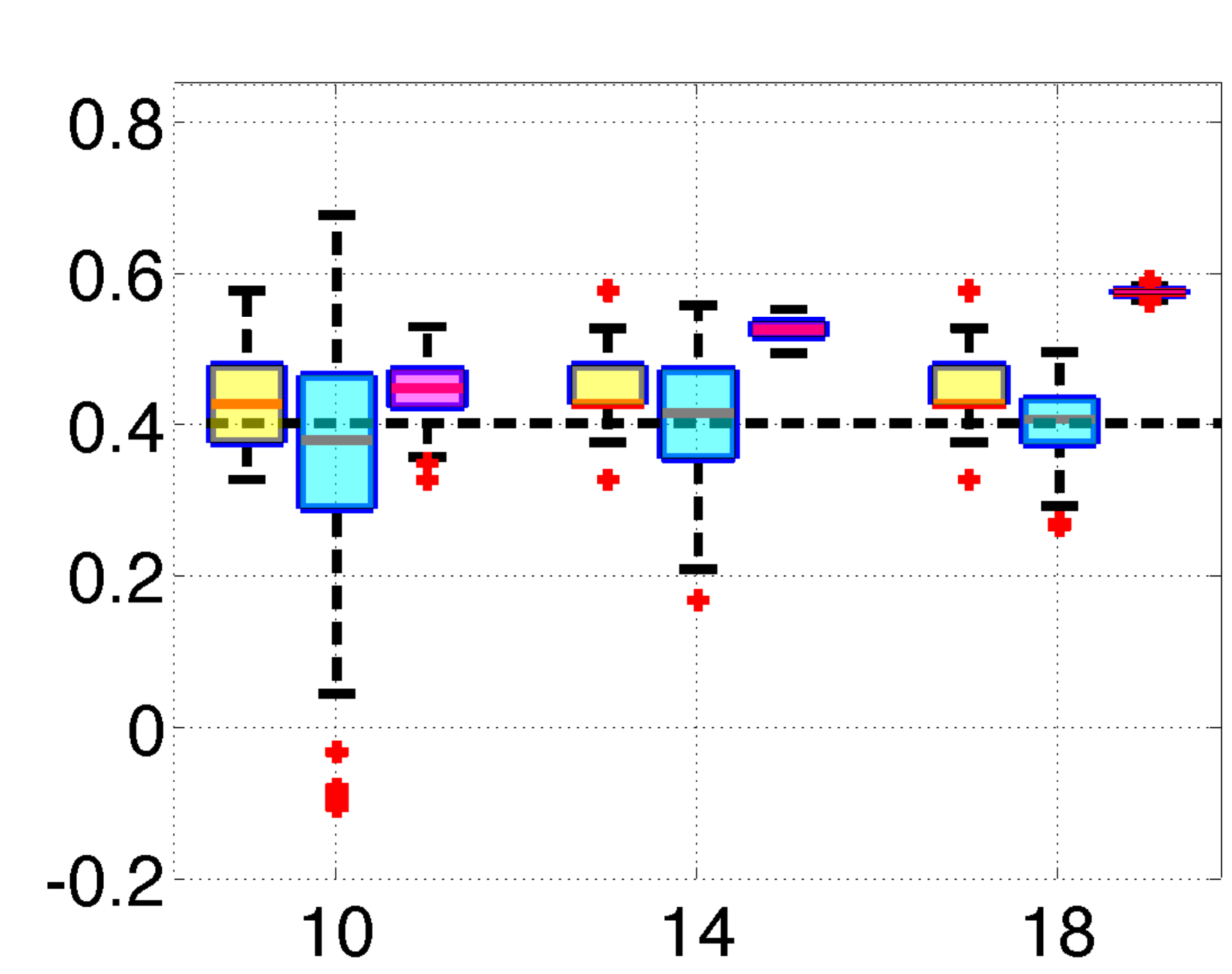}\\
 {\scriptsize $\hphantom{xxxx}\log_2 N$\hskip 1.9cm$\log_2 N$  \hskip 1.9cm$\log_2 N$}
   \vskip-2mm
 \caption{\textbf{Estimation performance of $h_1$} as function of $\log_2 N$.
 \label{fig:perf_h1}}
 \end{figure}

\begin{figure}[t]
    \centering
{\scriptsize \hskip .9cm$\rho_{\rm x}=0.1$\hskip 1.7cm$\rho_{\rm x}=0.45$  \hskip 1.4cm$\rho_{\rm x}=0.8$} \\
\rotatebox{90}{\hskip .2cm\scriptsize$\beta=\gamma=0$} \includegraphics[scale=\lscale, clip=true, trim=\trleftbord cm \trbot cm \trright cm \trtop cm]{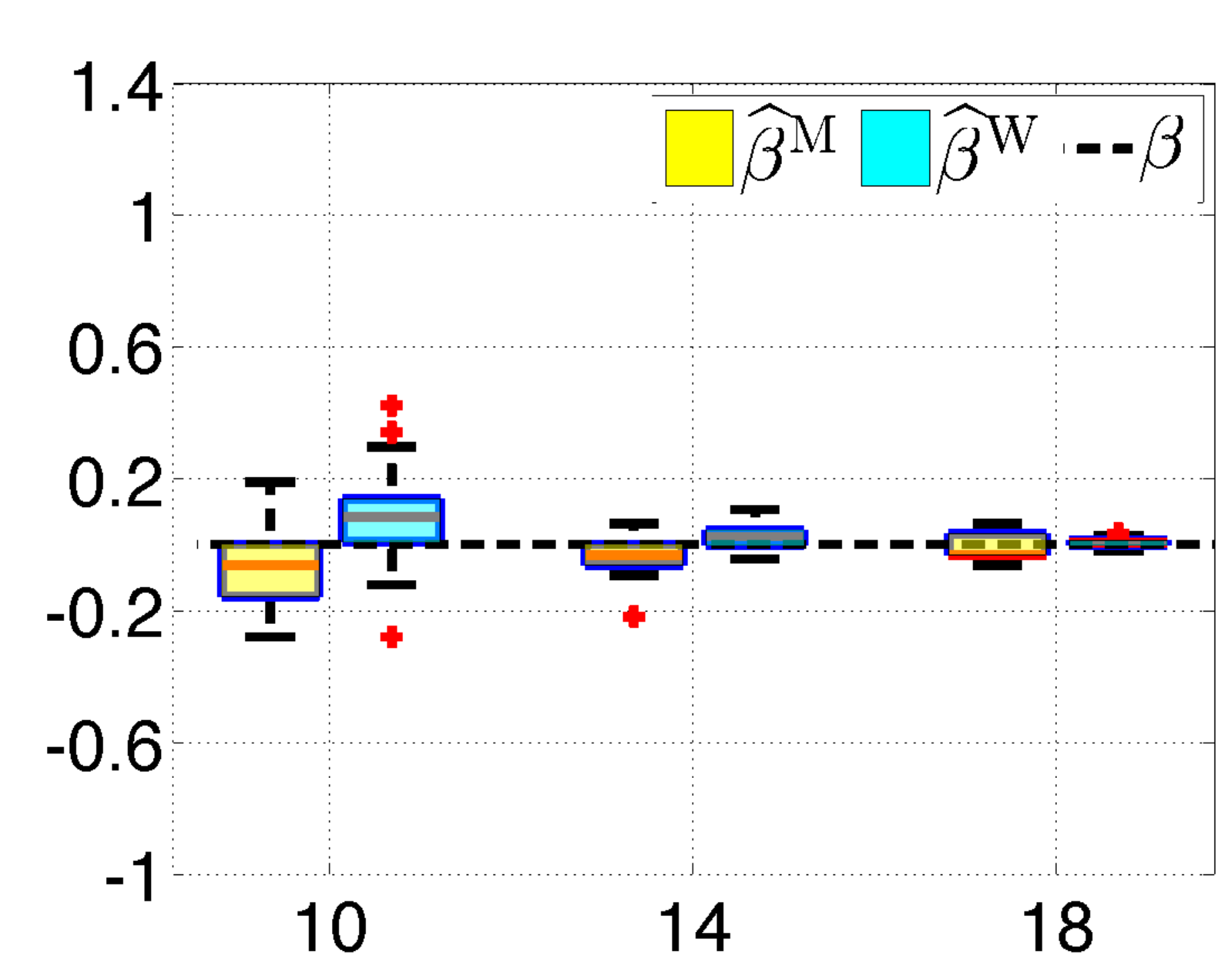}
  \includegraphics[scale=\lscale, clip=true, trim=\trleft cm \trbot cm \trright cm \trtop cm]{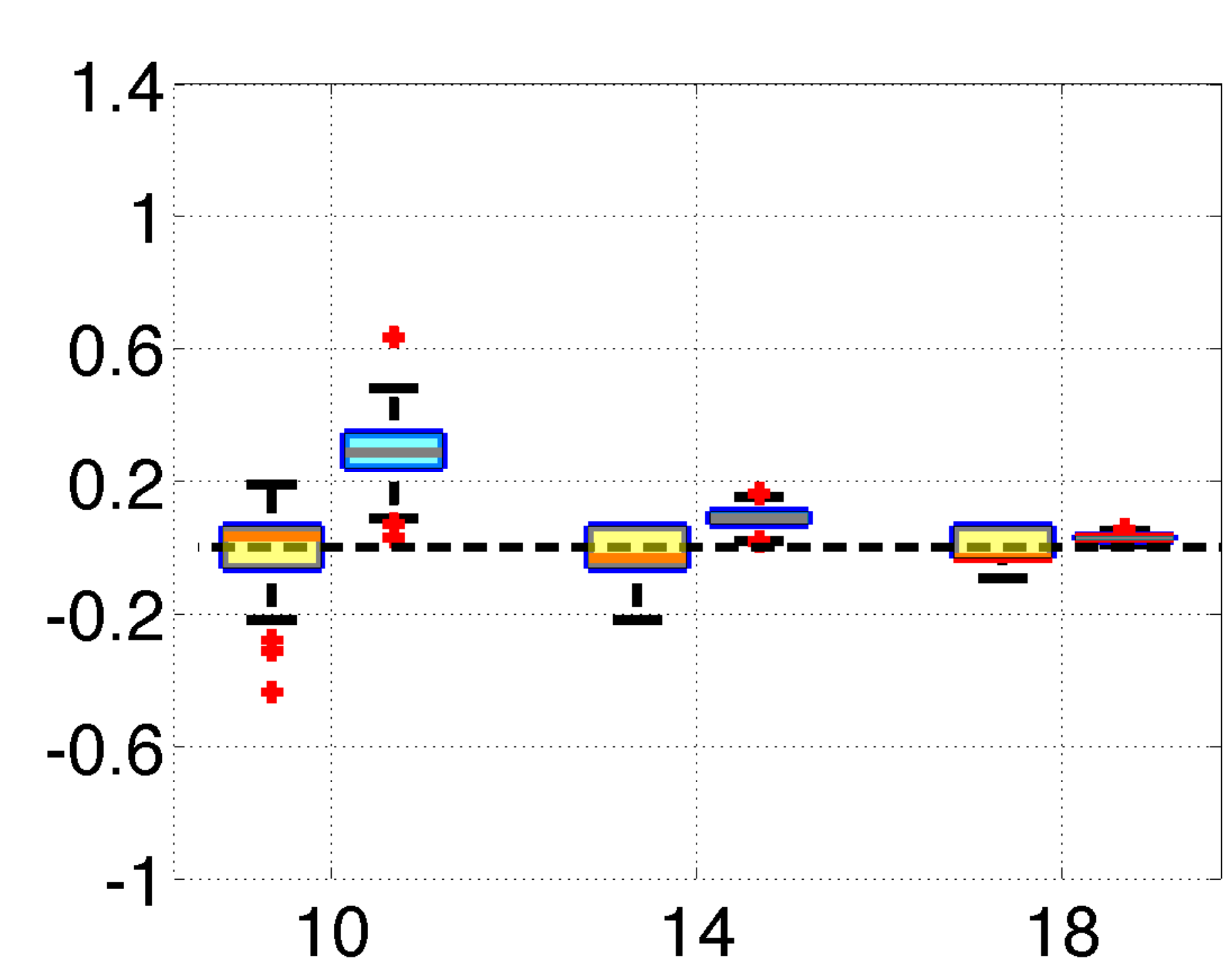}
 \includegraphics[scale=\lscale, clip=true, trim=\trleft cm \trbot cm \trright cm \trtop cm]{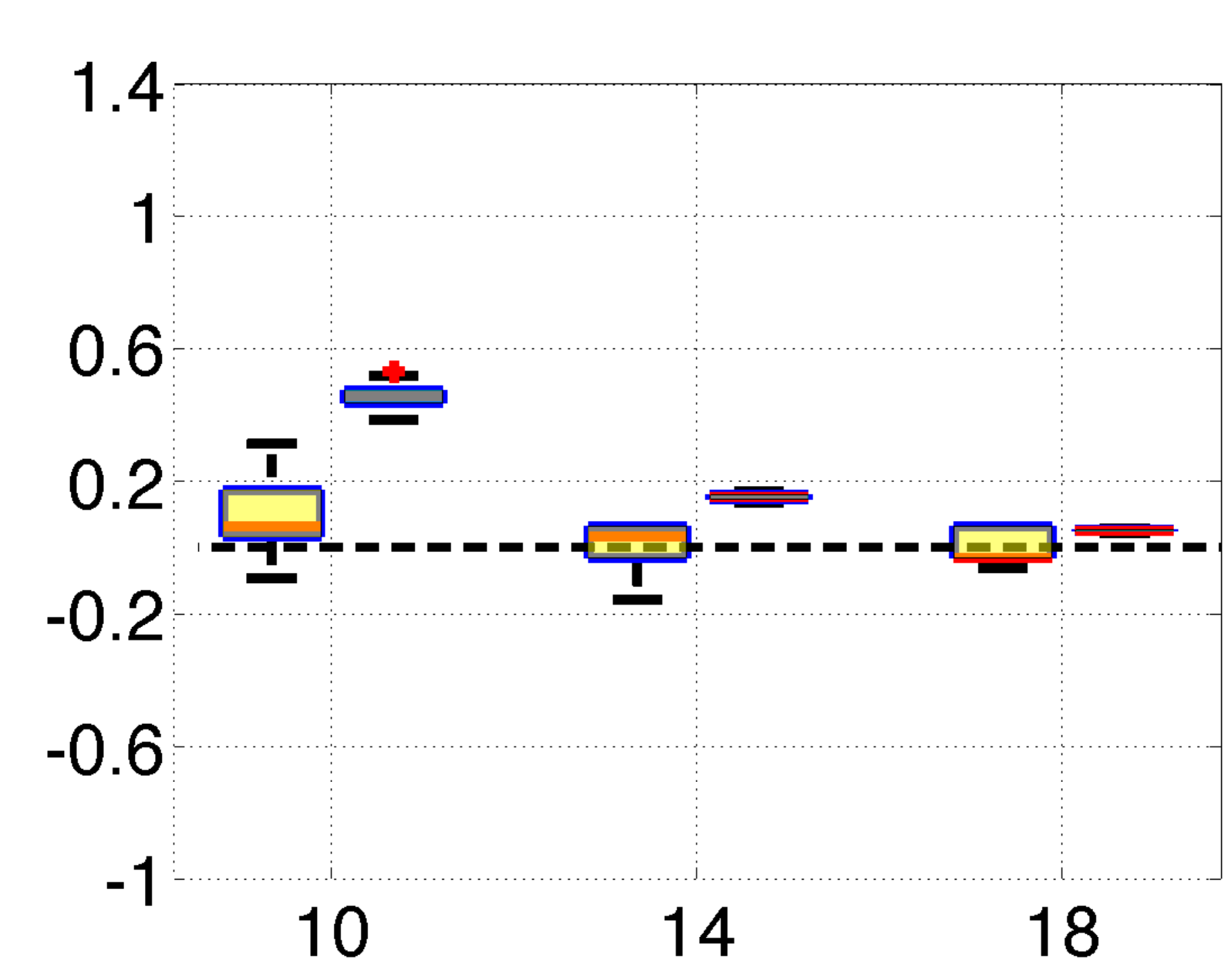}\\ \vskip -.3cm
 \rotatebox{90}{\hskip .2cm\scriptsize$\beta=\gamma=0.5$} \includegraphics[scale=\lscale, clip=true, trim=\trleftbord cm \trbot cm \trright cm \trtop cm]{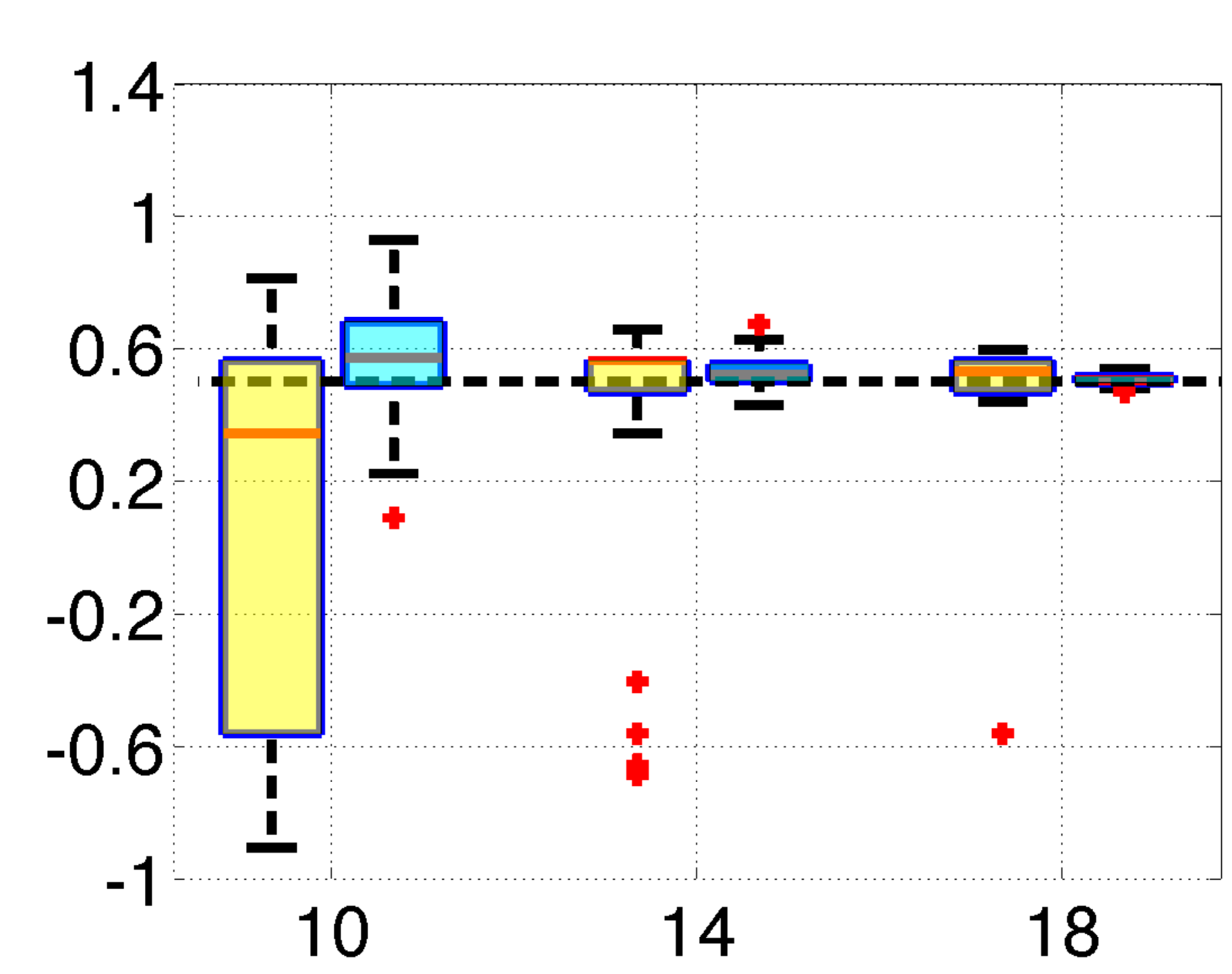}
  \includegraphics[scale=\lscale, clip=true, trim=\trleft cm \trbot cm \trright cm \trtop cm]{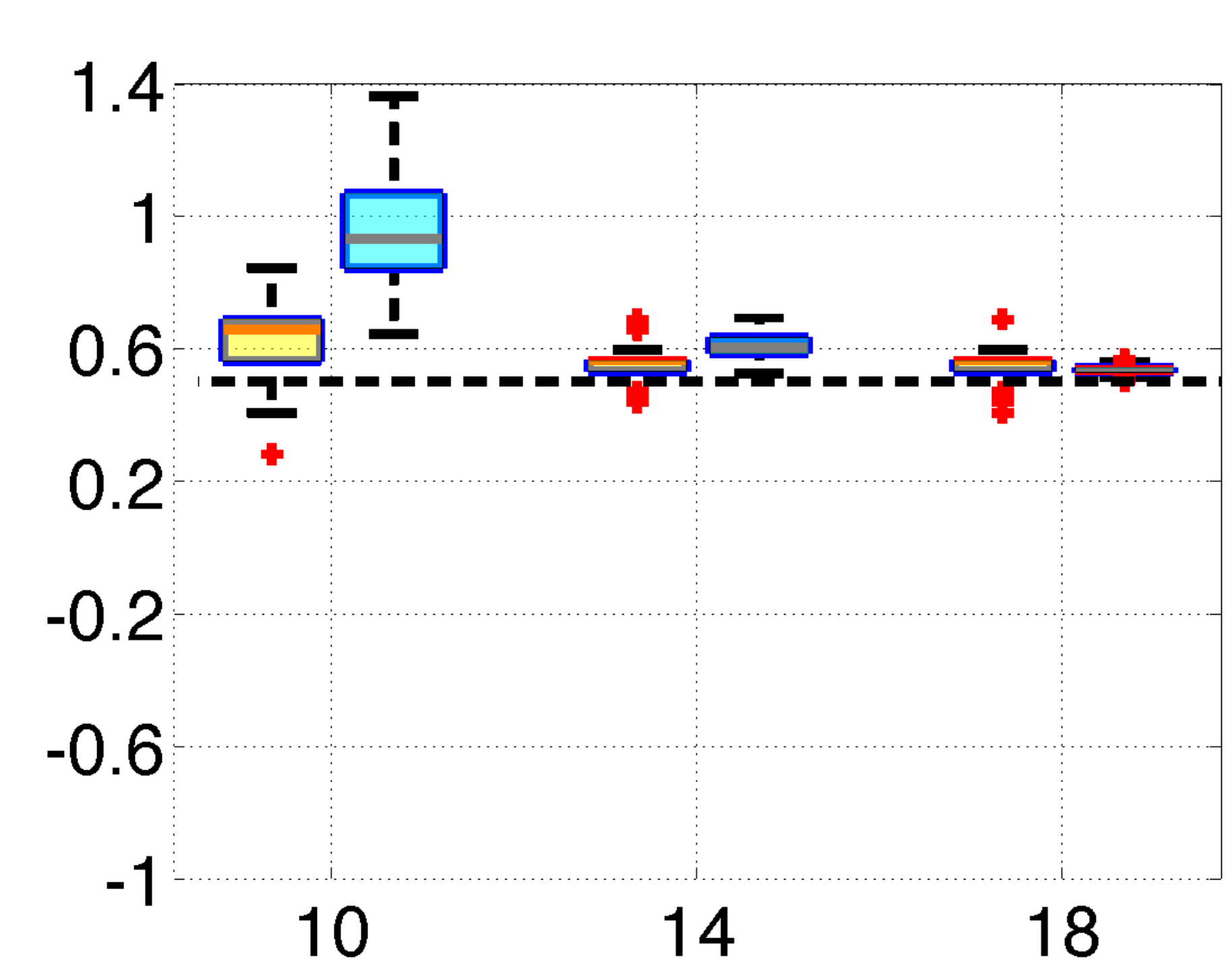}
 \includegraphics[scale=\lscale, clip=true, trim=\trleft cm \trbot cm \trright cm \trtop cm]{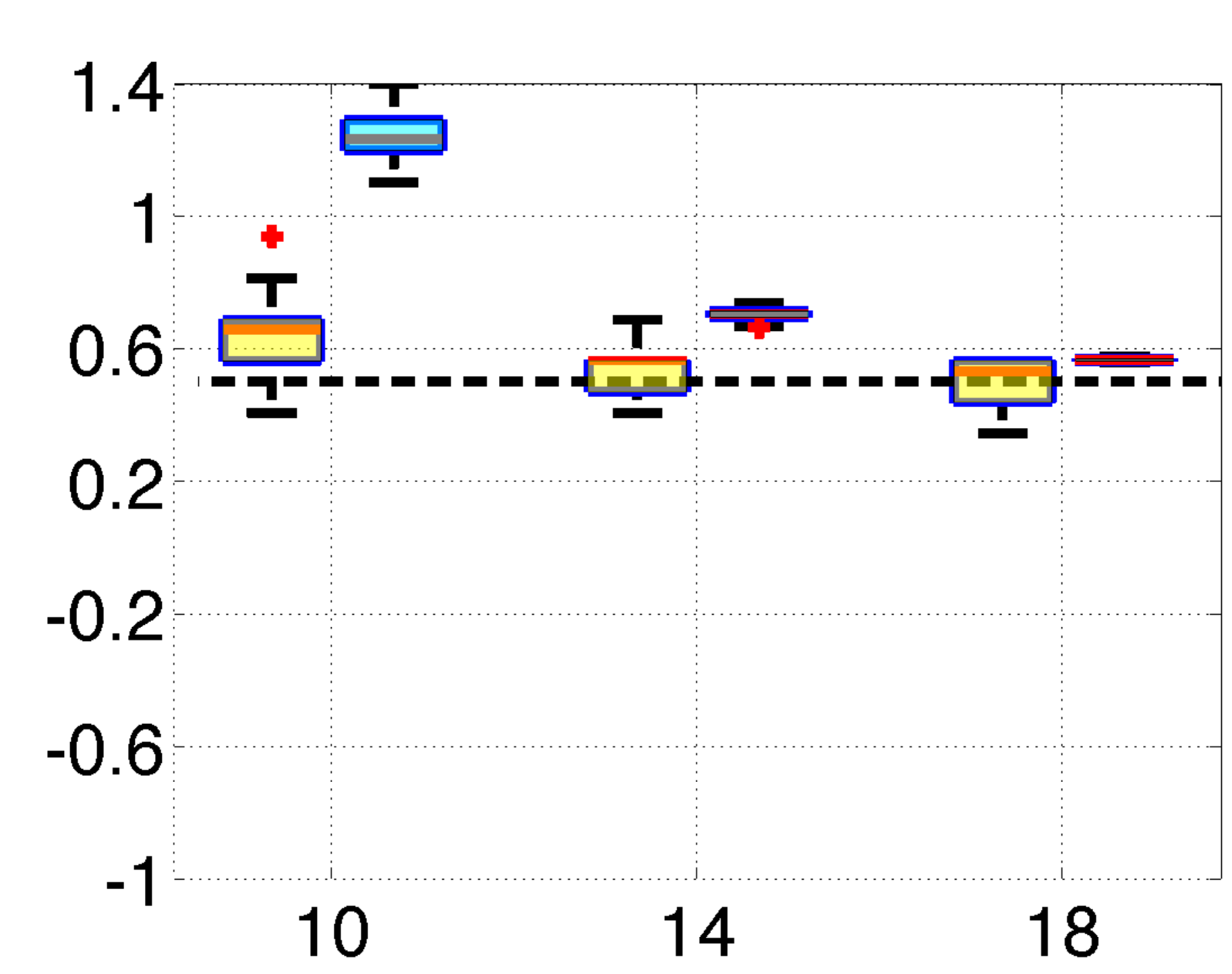}\\ \vskip -.3cm
 \rotatebox{90}{\hskip .1cm\scriptsize$\beta=-\gamma=0.5$} \includegraphics[scale=\lscale, clip=true, trim=\trleftbord cm \trbotbord cm \trright cm \trtop cm]{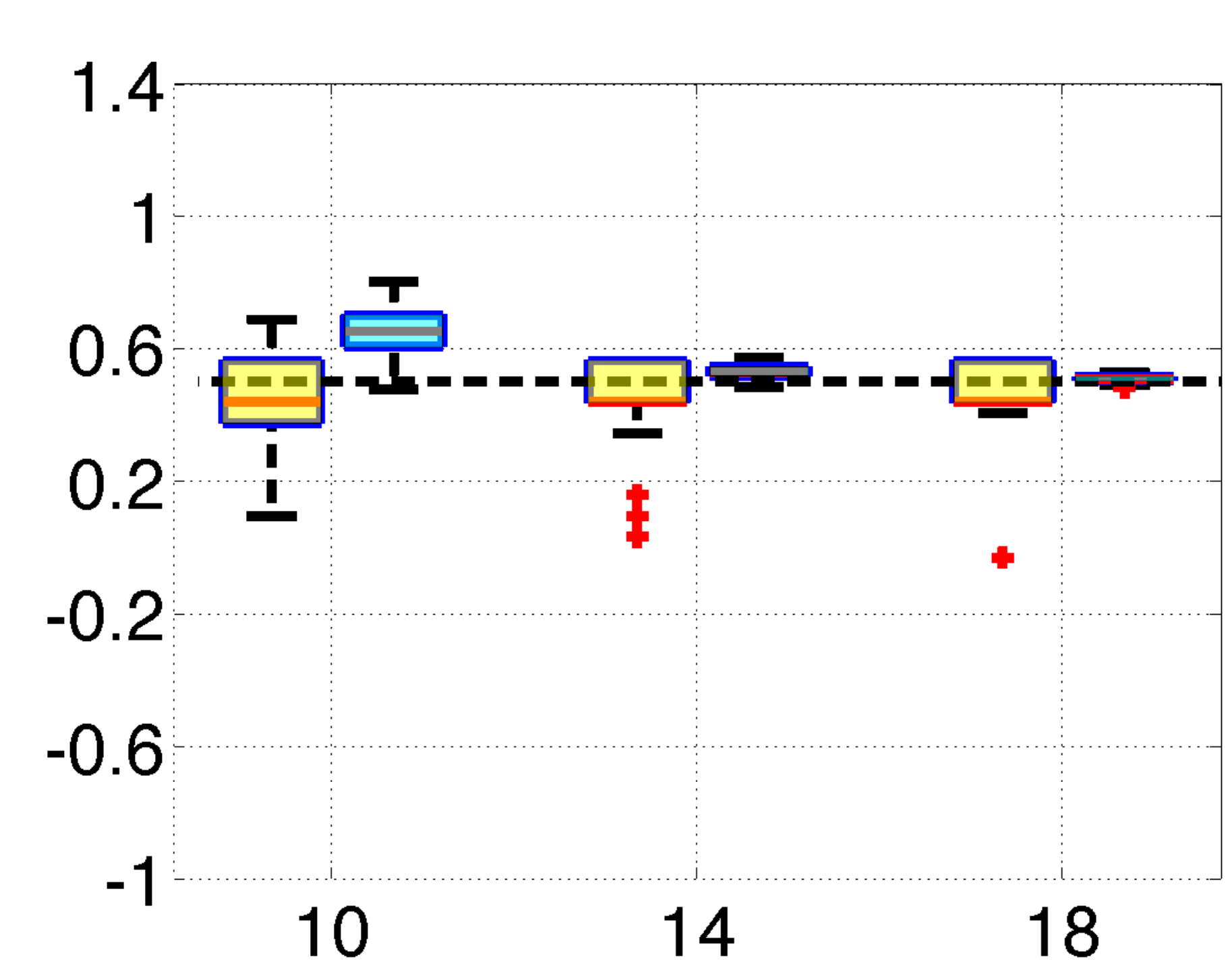}
  \includegraphics[scale=\lscale, clip=true, trim=\trleft cm \trbotbord cm \trright cm \trtop cm]{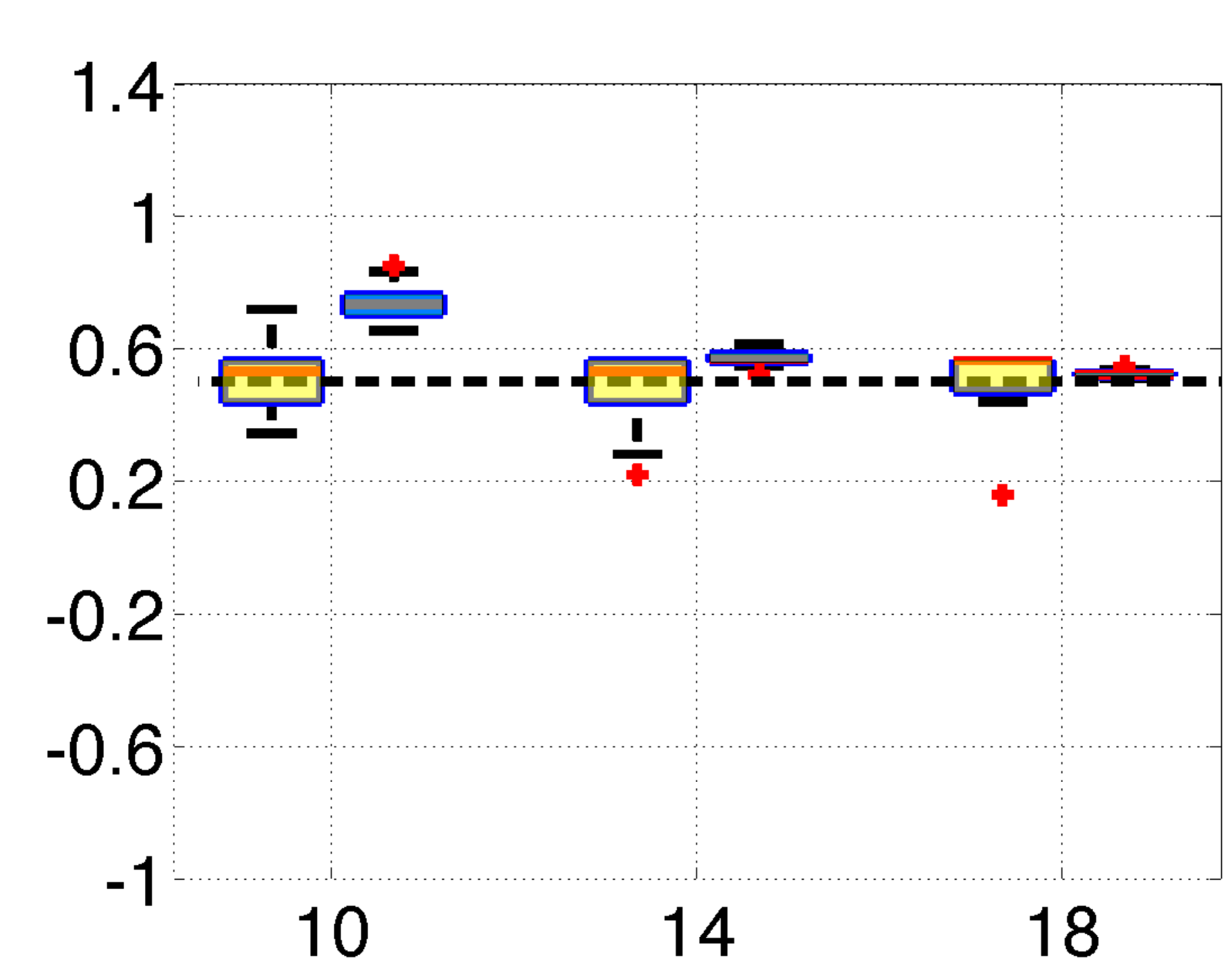}
 \includegraphics[scale=\lscale, clip=true, trim=\trleft cm \trbotbord cm \trright cm \trtop cm]{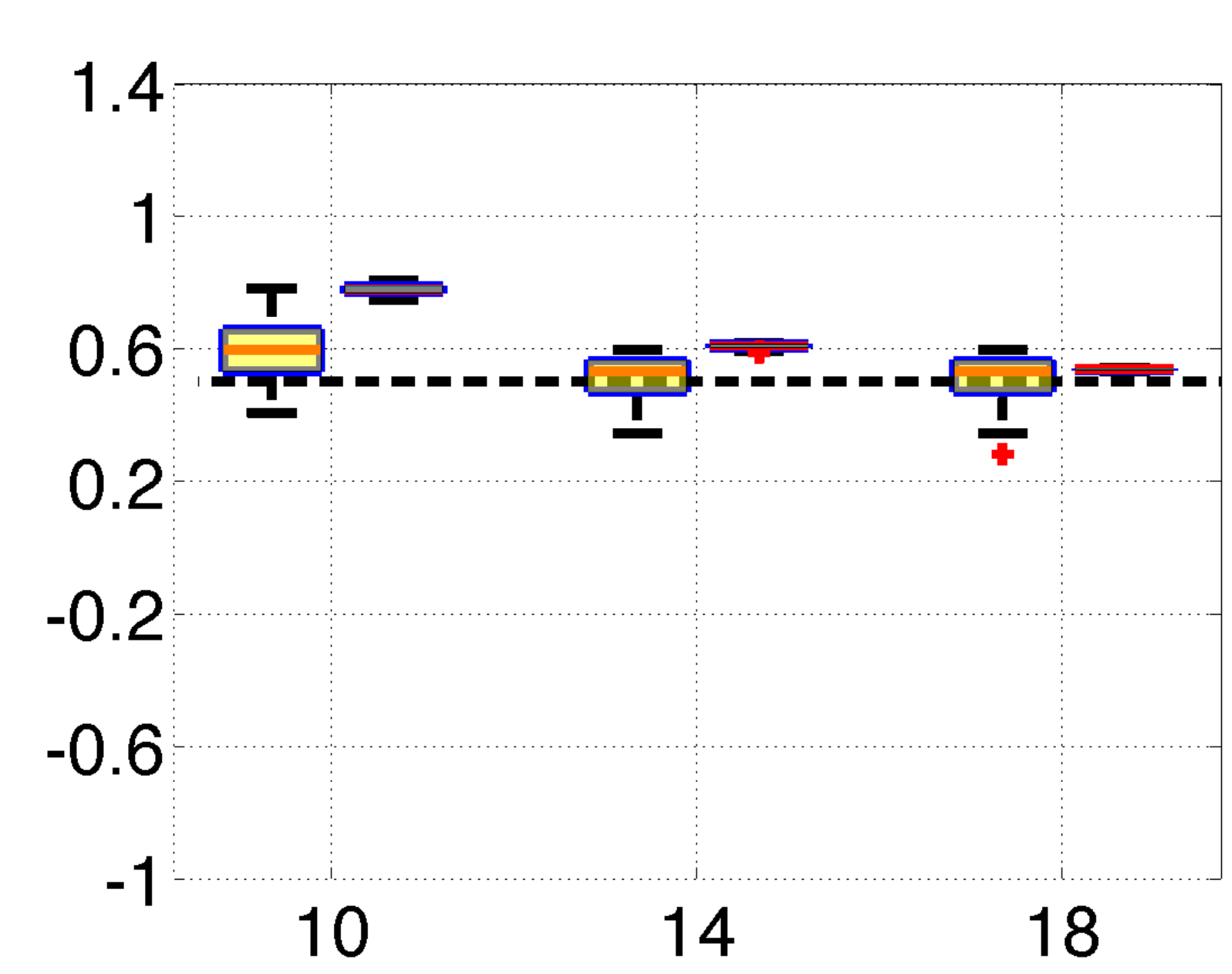}\\
 {\scriptsize $\hphantom{xxxx}\log_2 N$\hskip 1.9cm$\log_2 N$  \hskip 1.9cm$\log_2 N$}
   \vskip-2mm
 \caption{\textbf{Estimation performance of $\beta$} as function of $\log_2 N$.
 \label{fig:perf_beta}}
 \end{figure}

\begin{figure}[t]
      \centering
{\scriptsize \hskip .9cm$\rho_{\rm x}=0.1$\hskip 1.7cm$\rho_{\rm x}=0.45$  \hskip 1.4cm$\rho_{\rm x}=0.8$} \\
\rotatebox{90}{\hskip .2cm\scriptsize$\beta=\gamma=0$} \includegraphics[scale=\lscale, clip=true, trim=\trleftbord cm \trbot cm \trright cm \trtop cm]{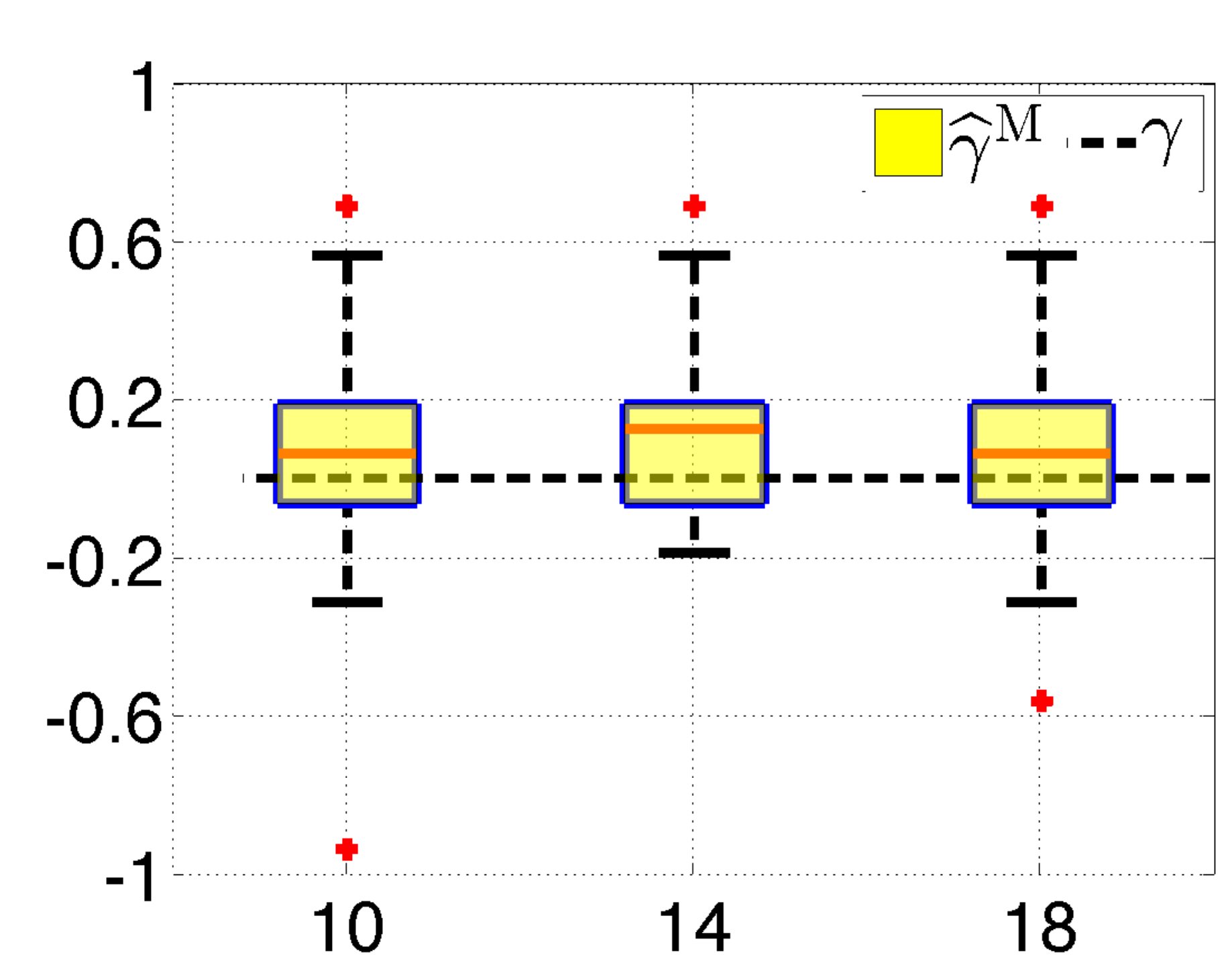}
  \includegraphics[scale=\lscale, clip=true, trim=\trleft cm \trbot cm \trright cm \trtop cm]{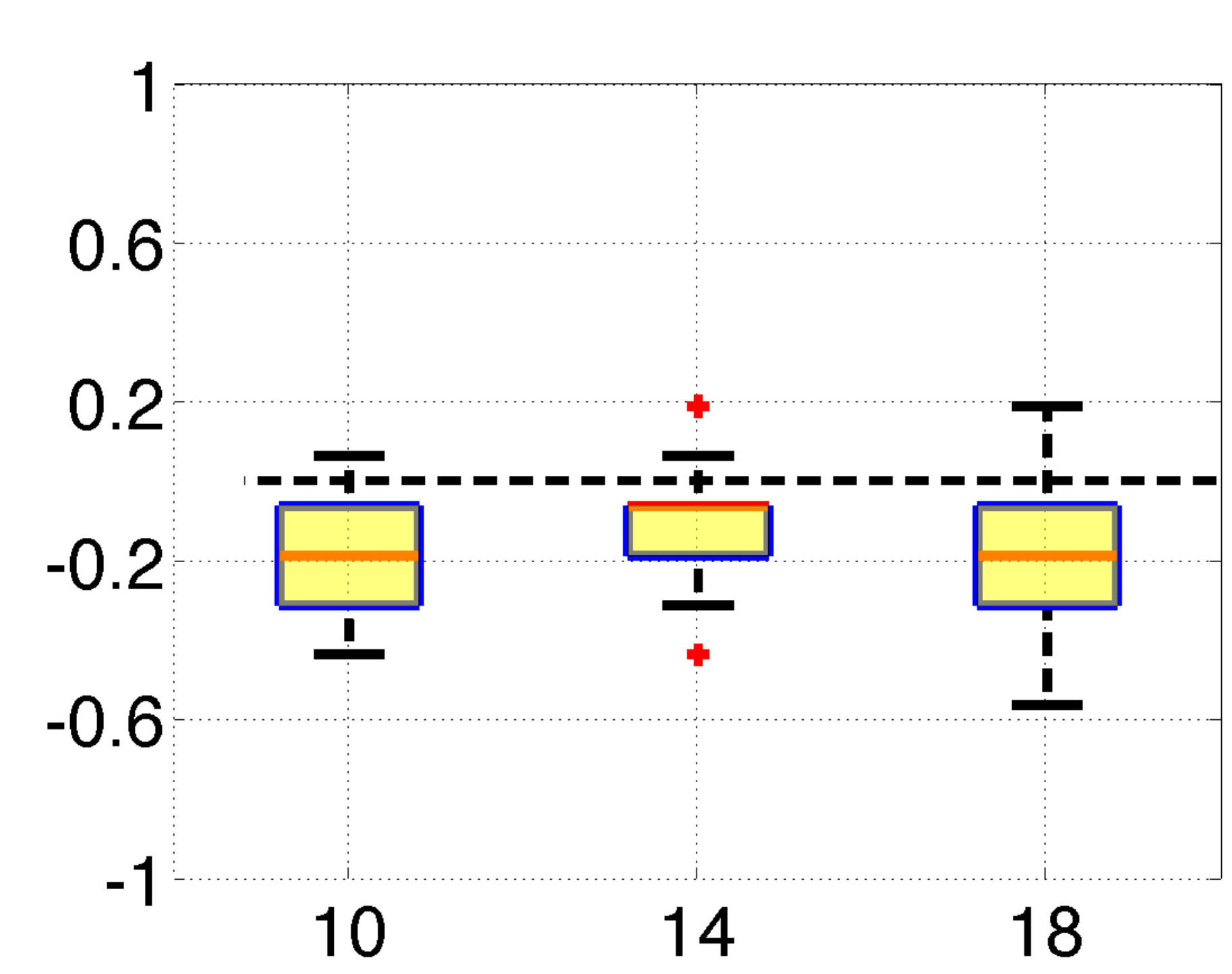}
 \includegraphics[scale=\lscale, clip=true, trim=\trleft cm \trbot cm \trright cm \trtop cm]{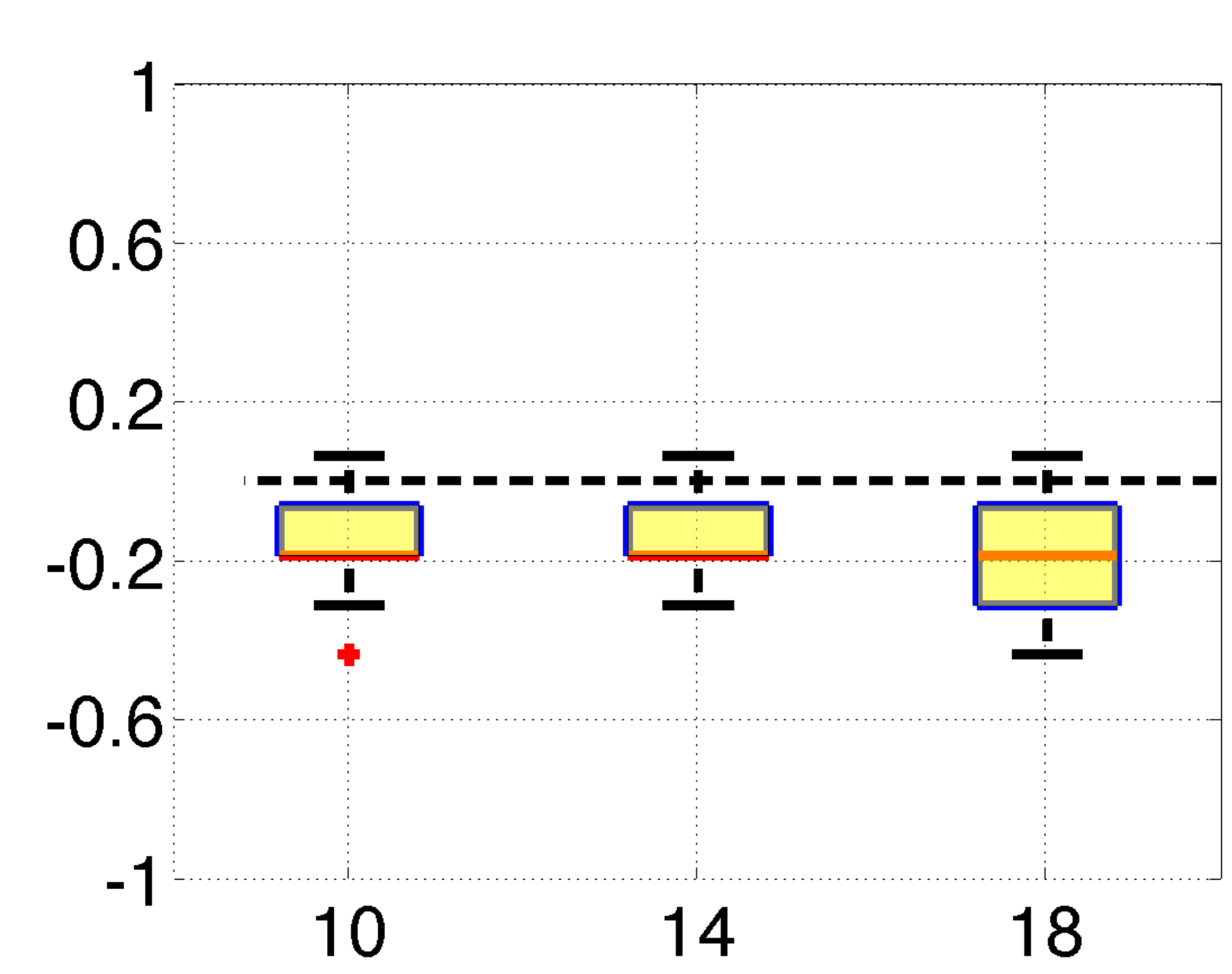}\\ \vskip -.3cm
 \rotatebox{90}{\hskip .2cm\scriptsize$\beta=\gamma=0.5$} \includegraphics[scale=\lscale, clip=true, trim=\trleftbord cm \trbot cm \trright cm \trtop cm]{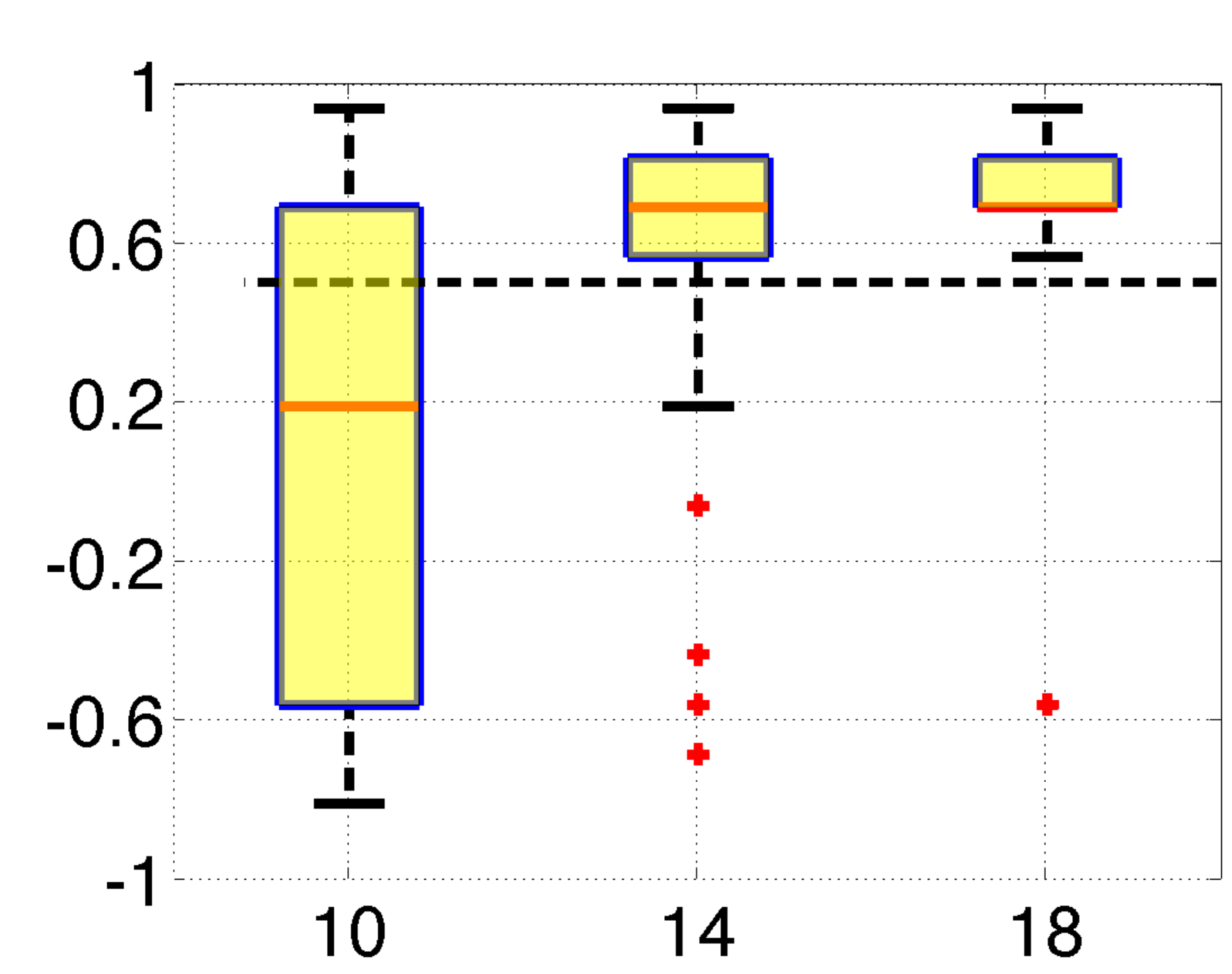}
  \includegraphics[scale=\lscale, clip=true, trim=\trleft cm \trbot cm \trright cm \trtop cm]{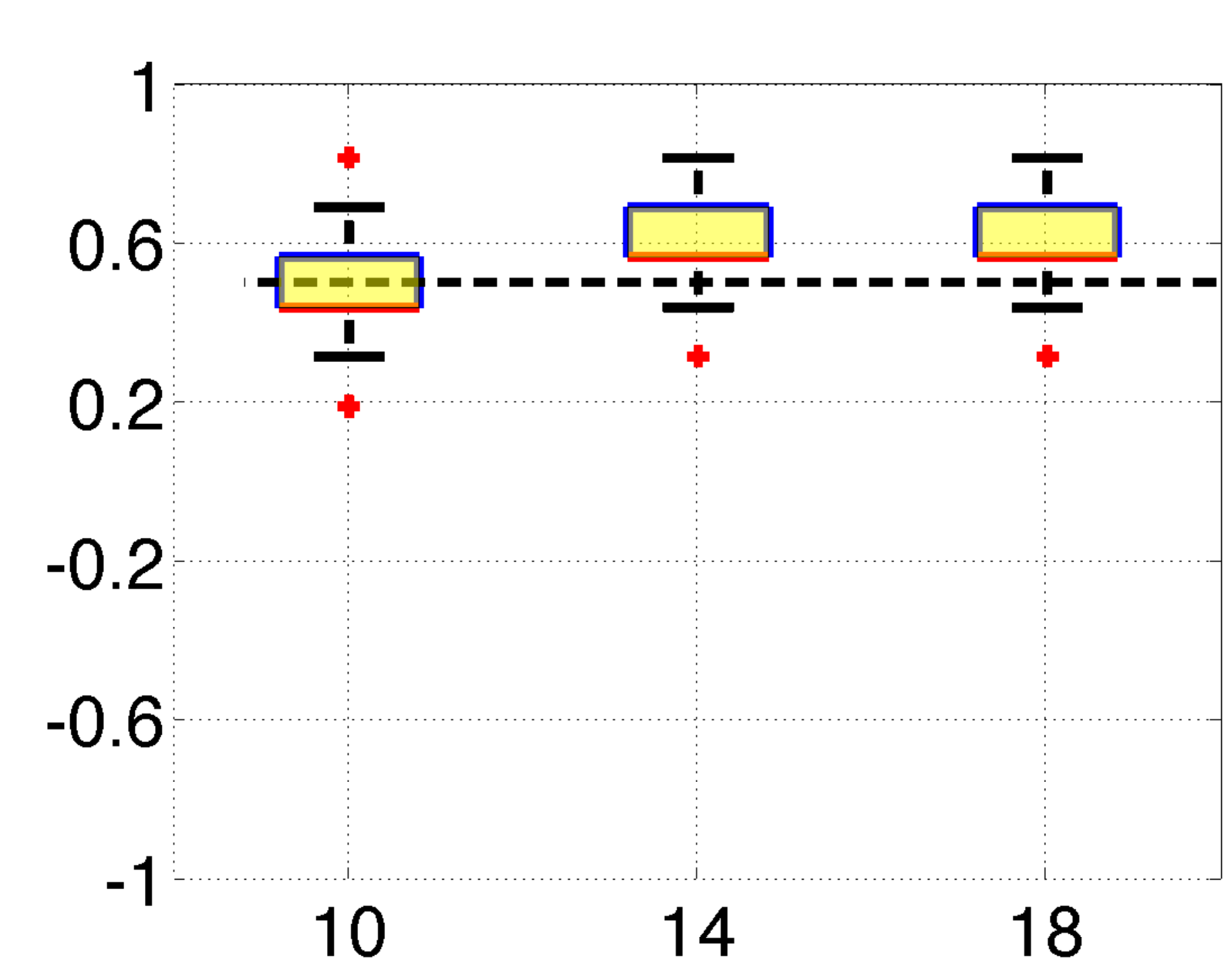}
 \includegraphics[scale=\lscale, clip=true, trim=\trleft cm \trbot cm \trright cm \trtop cm]{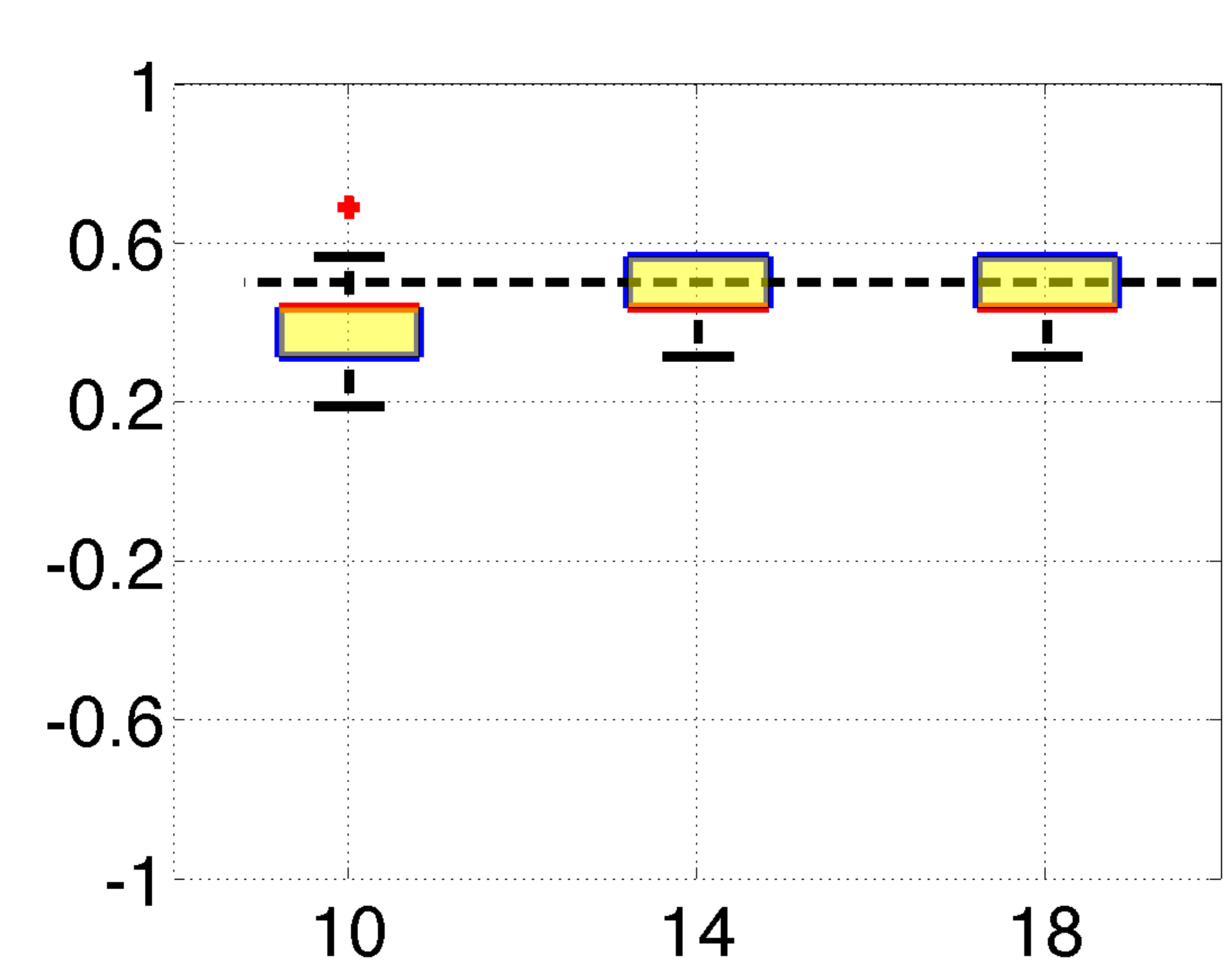}\\ \vskip -.3cm
 \rotatebox{90}{\hskip .1cm\scriptsize$\beta=-\gamma=0.5$} \includegraphics[scale=\lscale, clip=true, trim=\trleftbord cm \trbotbord cm \trright cm \trtop cm]{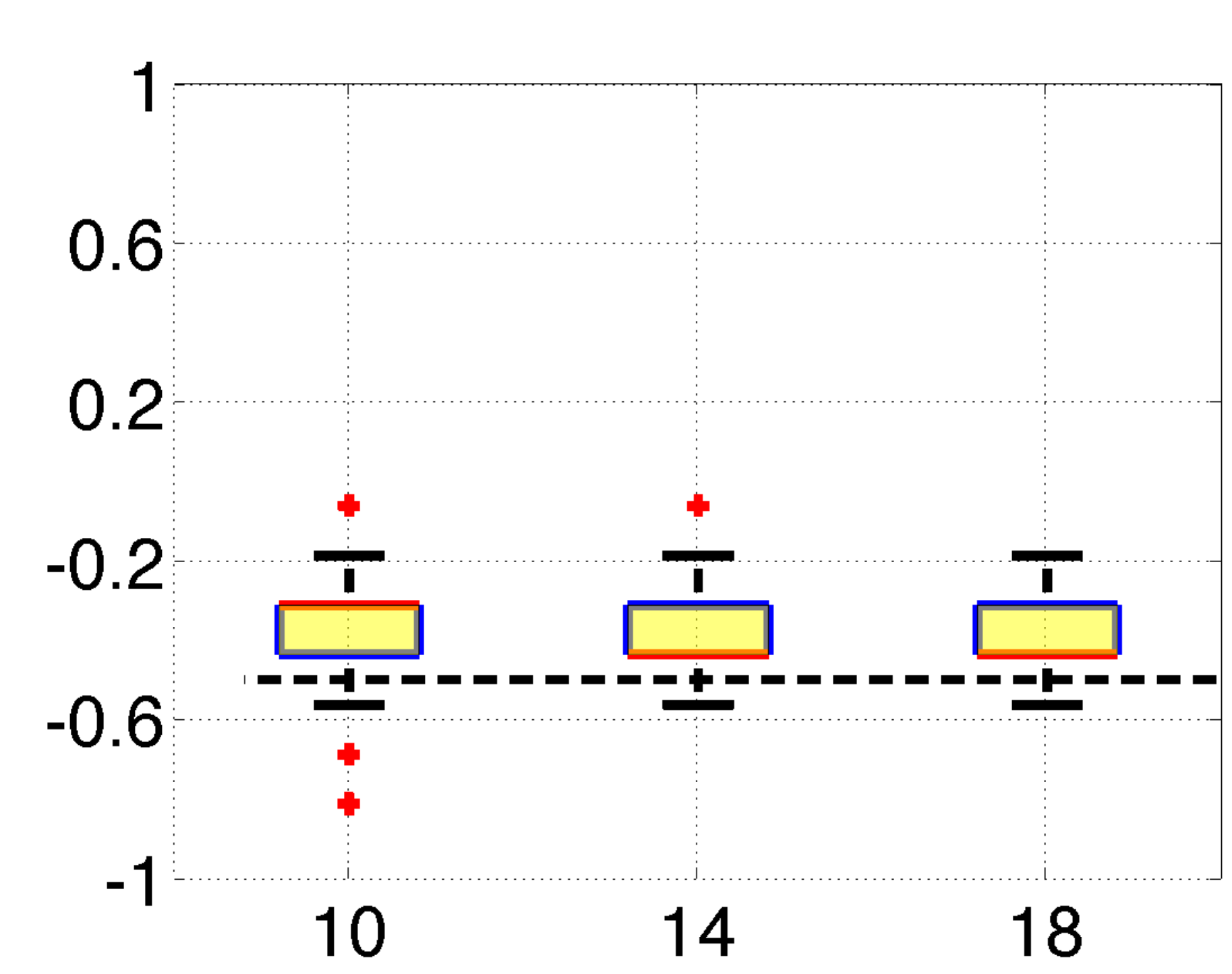}
  \includegraphics[scale=\lscale, clip=true, trim=\trleft cm \trbotbord cm \trright cm \trtop cm]{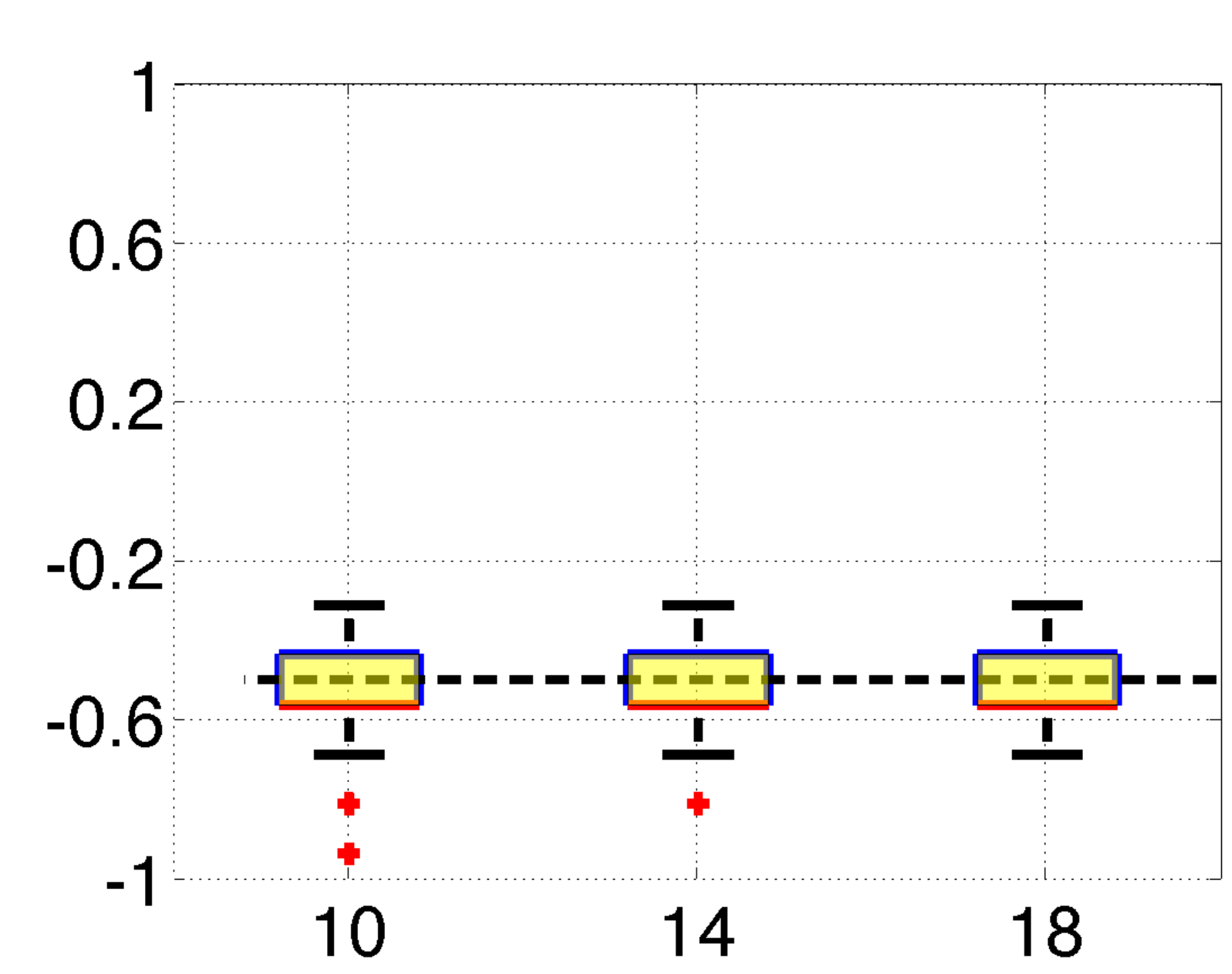}
 \includegraphics[scale=\lscale, clip=true, trim=\trleft cm \trbotbord cm \trright cm \trtop cm]{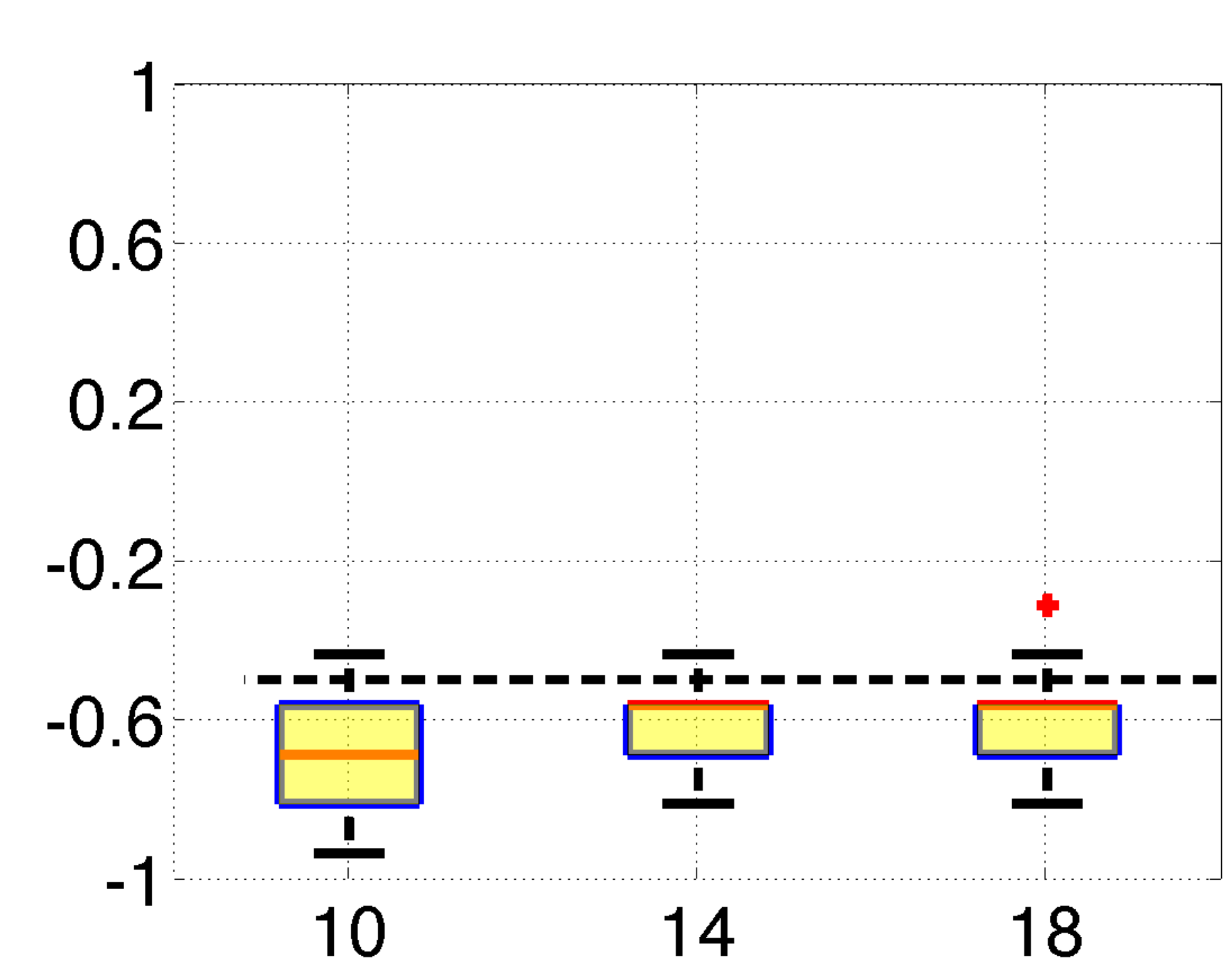}\\
 {\scriptsize $\hphantom{xxxx}\log_2 N$\hskip 1.9cm$\log_2 N$  \hskip 1.9cm$\log_2 N$}
   \vskip-2mm
 \caption{\textbf{Estimation performance of $\gamma$} as function of $\log_2 N$.
 \label{fig:perf_gamma}}
 \end{figure}

\begin{figure}[t]
 \centering
{\scriptsize \hskip .9cm$\rho_{\rm x}=0.1$\hskip 1.7cm$\rho_{\rm x}=0.45$  \hskip 1.4cm$\rho_{\rm x}=0.8$} \\
\rotatebox{90}{\hskip .2cm\scriptsize$\beta=\gamma=0$} \includegraphics[scale=\lscale, clip=true, trim=\trleftbord cm \trbot cm \trright cm \trtop cm]{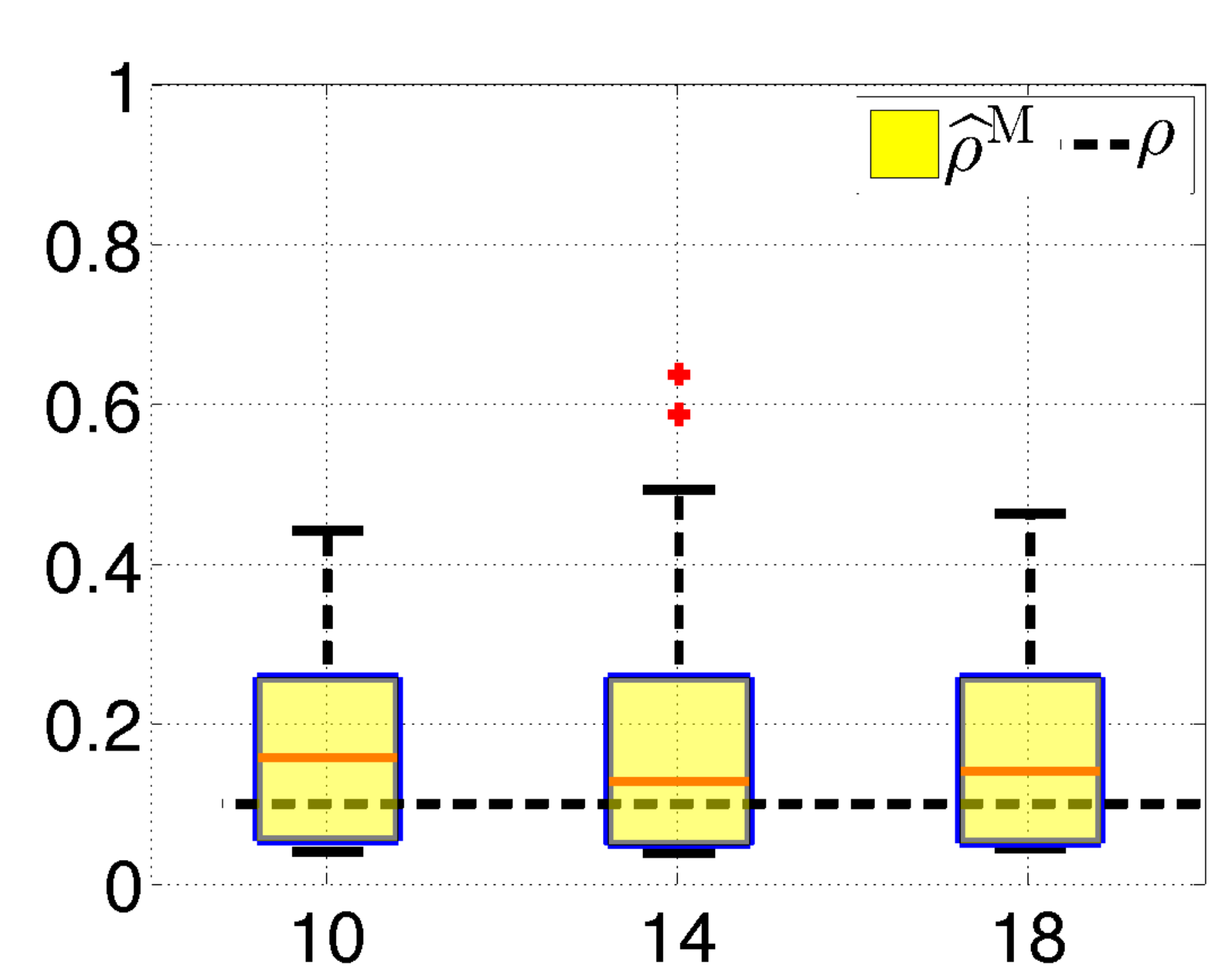}
  \includegraphics[scale=\lscale, clip=true, trim=\trleft cm \trbot cm \trright cm \trtop cm]{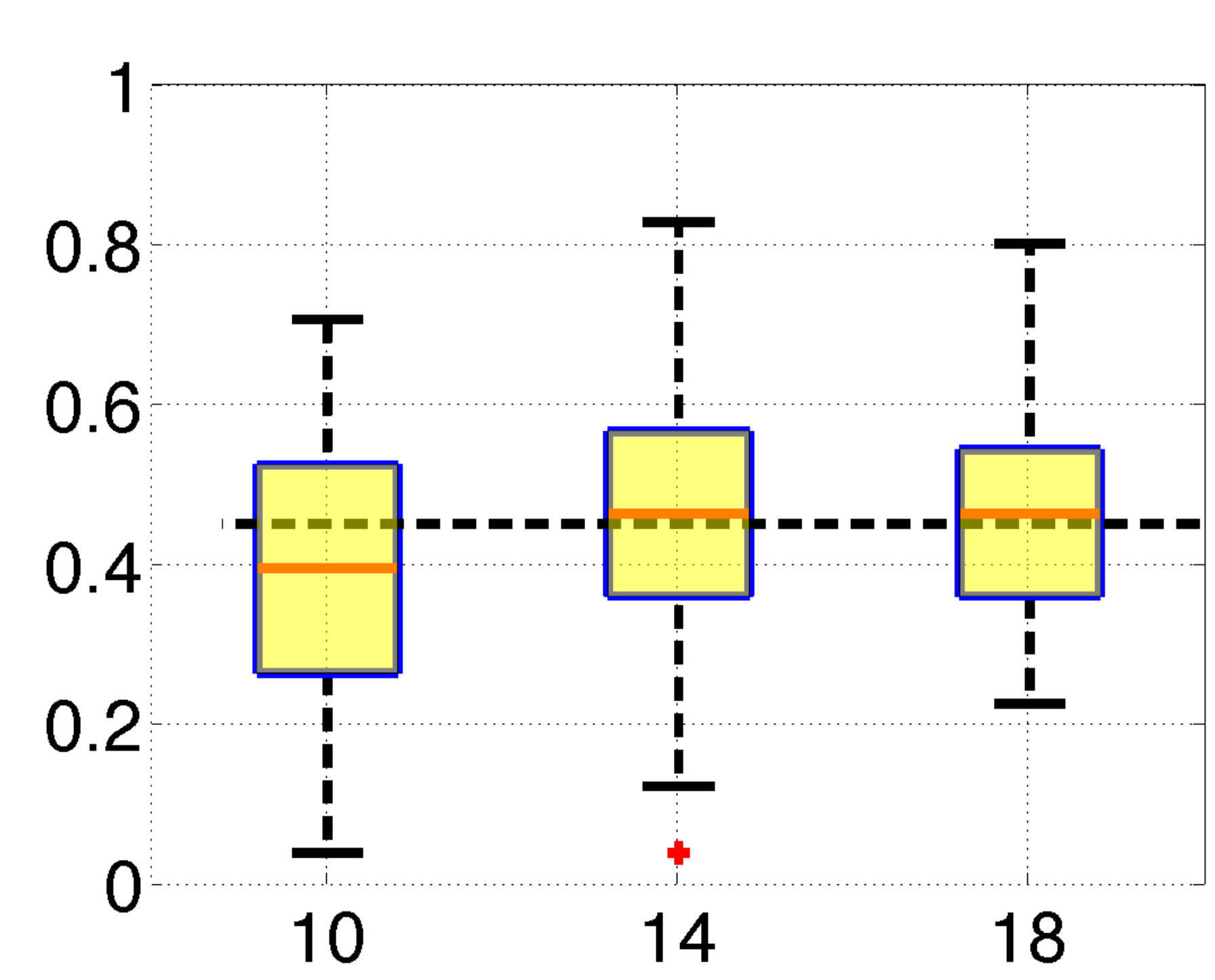}
 \includegraphics[scale=\lscale, clip=true, trim=\trleft cm \trbot cm \trright cm \trtop cm]{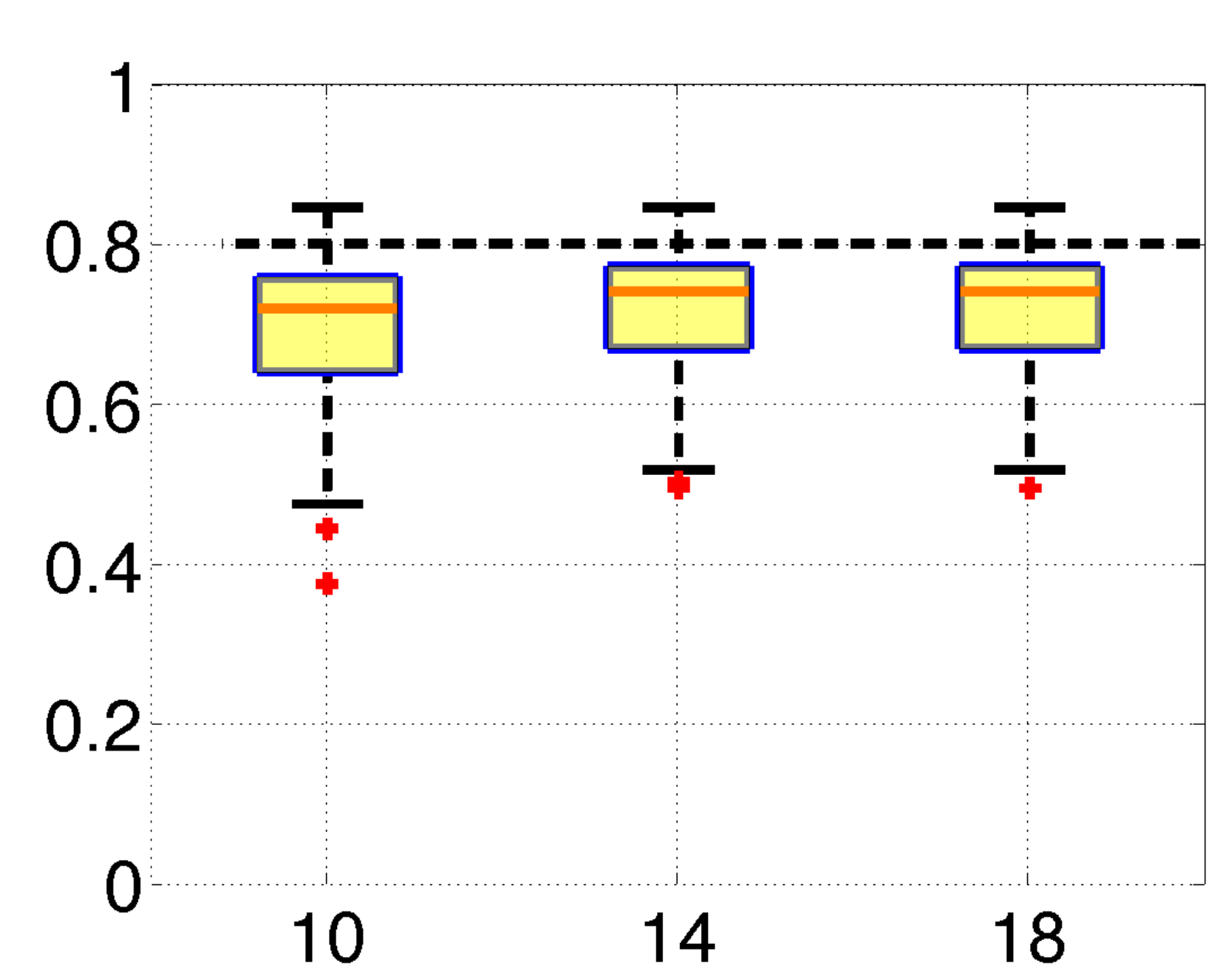}\\ \vskip -.3cm
 \rotatebox{90}{\hskip .2cm\scriptsize$\beta=\gamma=0.5$} \includegraphics[scale=\lscale, clip=true, trim=\trleftbord cm \trbot cm \trright cm \trtop cm]{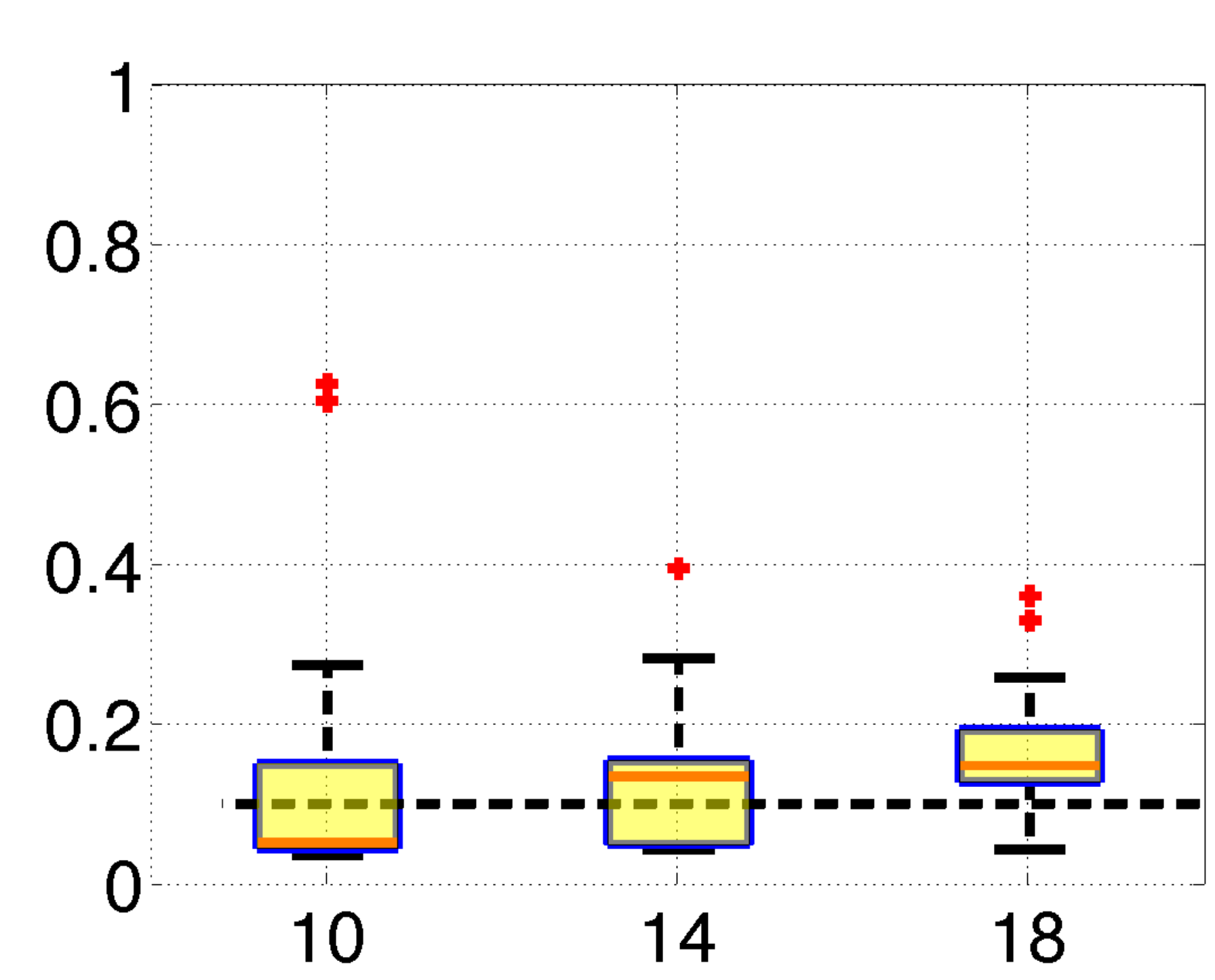}
  \includegraphics[scale=\lscale, clip=true, trim=\trleft cm \trbot cm \trright cm \trtop cm]{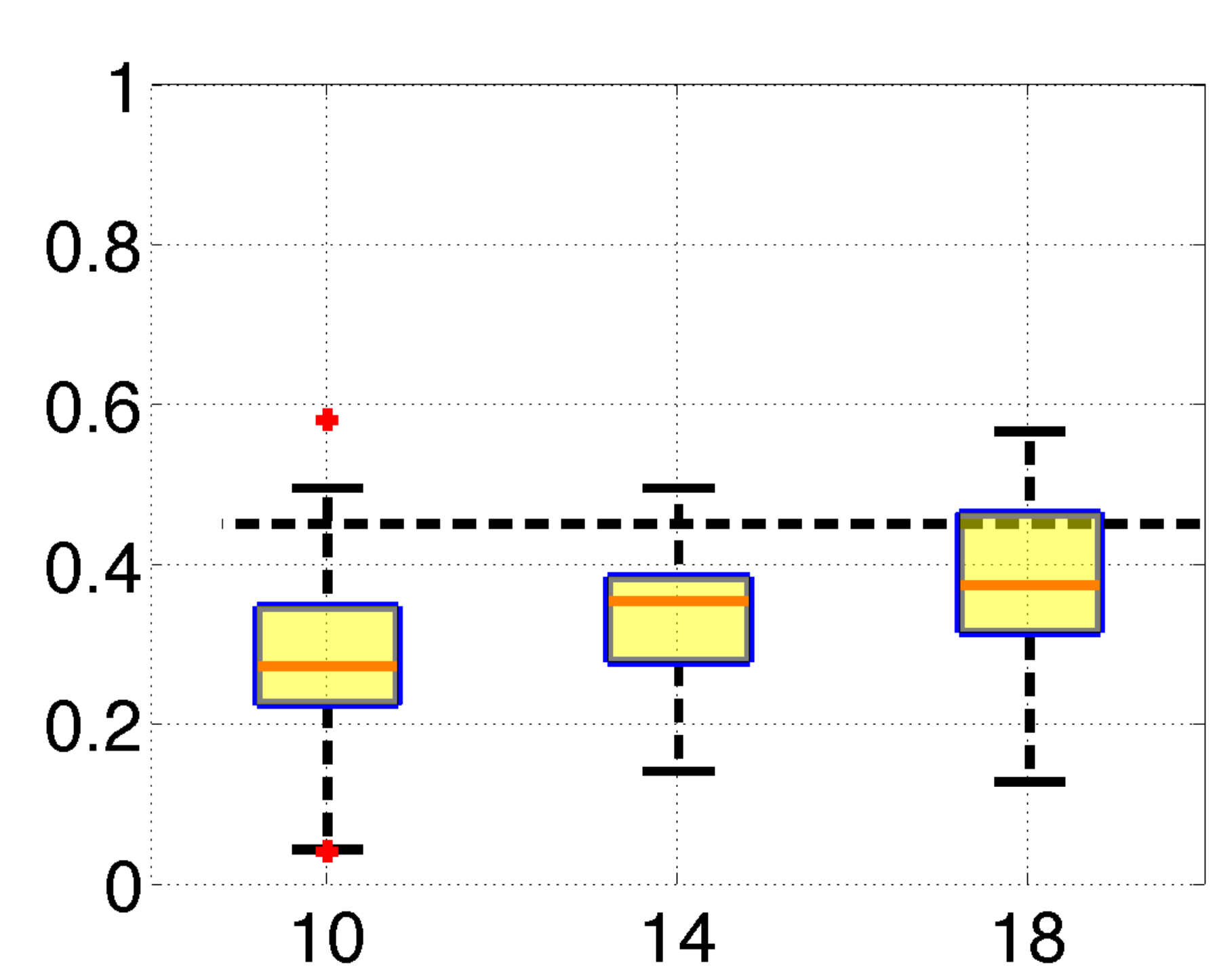}
 \includegraphics[scale=\lscale, clip=true, trim=\trleft cm \trbot cm \trright cm \trtop cm]{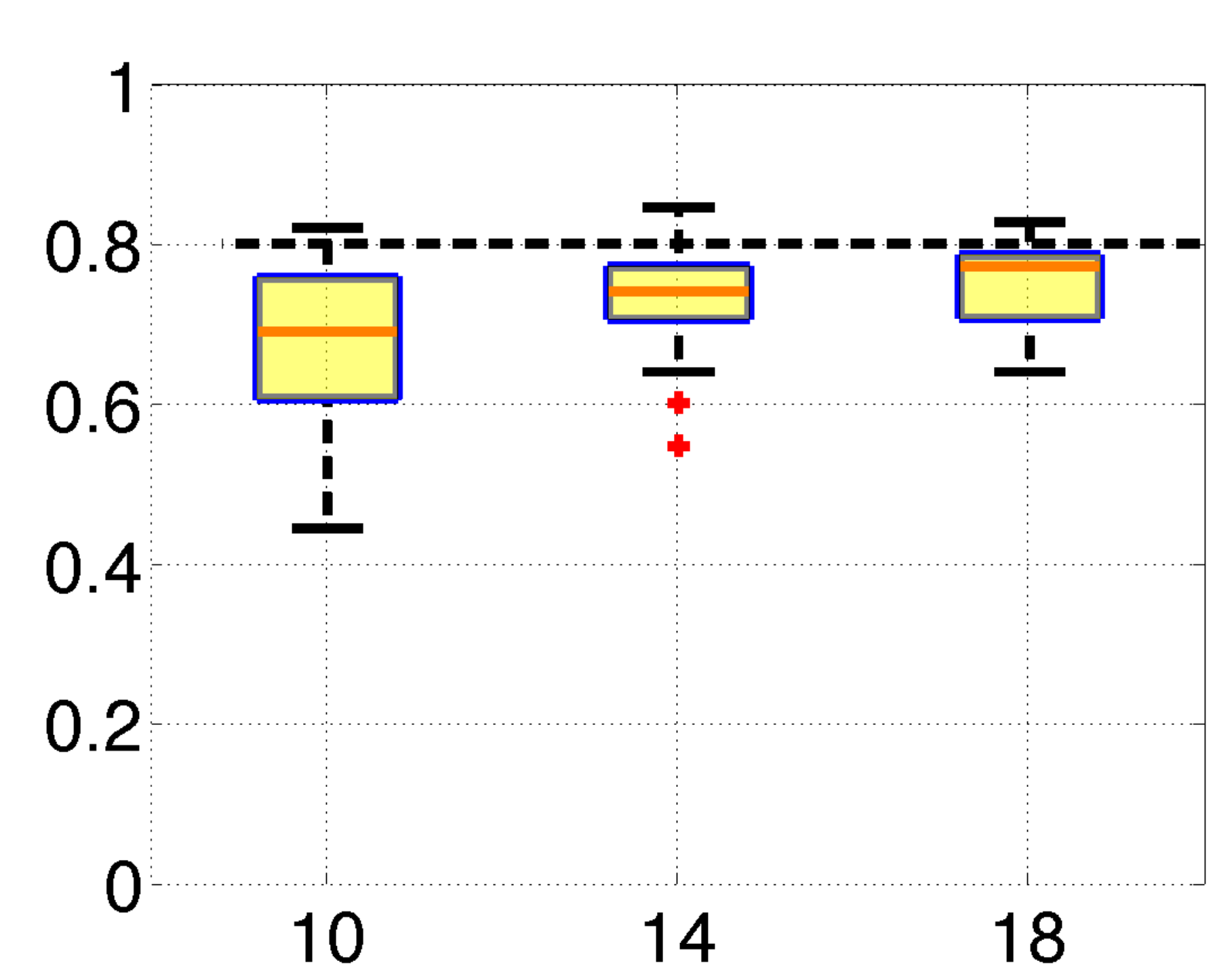}\\ \vskip -.3cm
 \rotatebox{90}{\hskip .1cm\scriptsize$\beta=-\gamma=0.5$} \includegraphics[scale=\lscale, clip=true, trim=\trleftbord cm \trbotbord cm \trright cm \trtop cm]{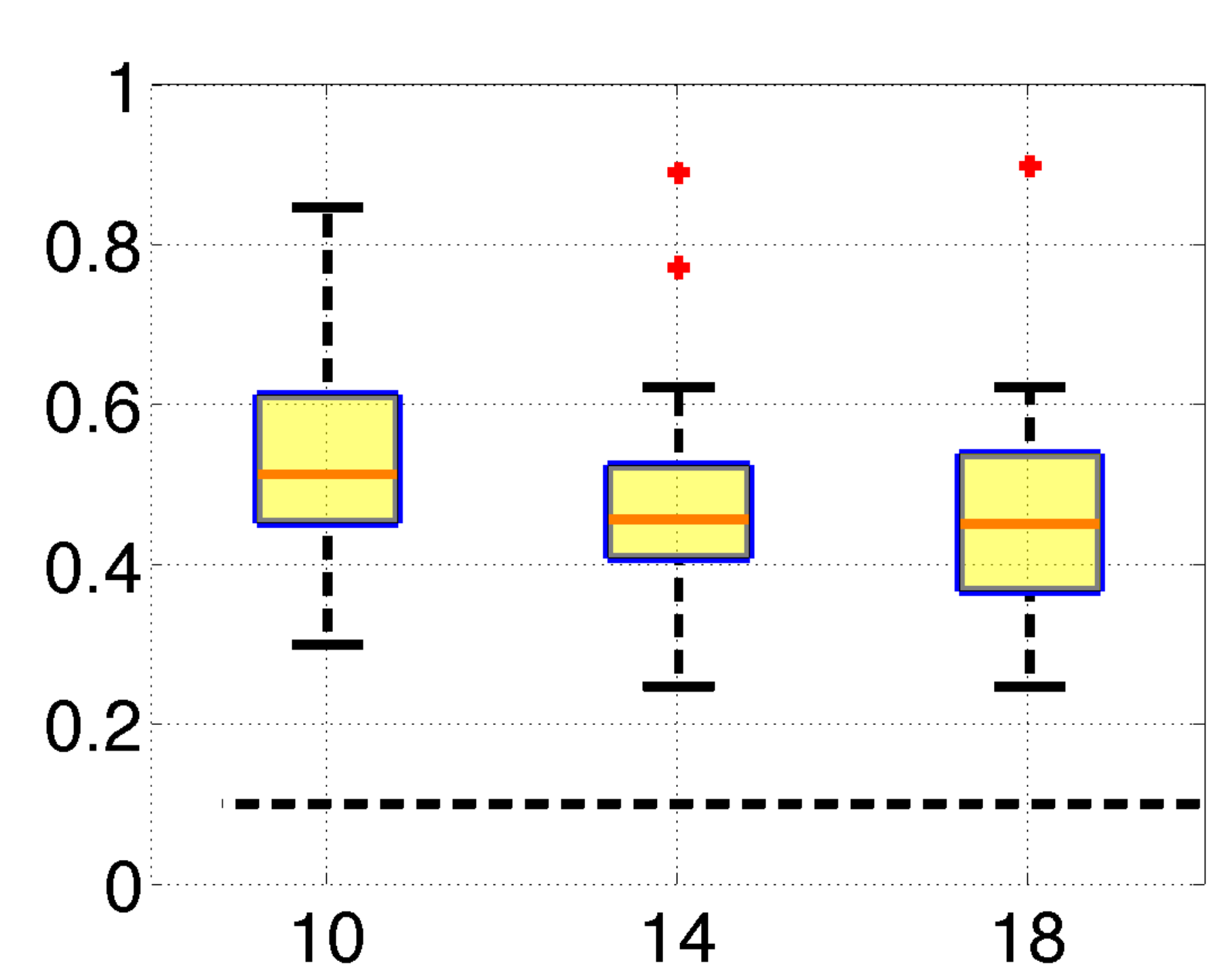}
  \includegraphics[scale=\lscale, clip=true, trim=\trleft cm \trbotbord cm \trright cm \trtop cm]{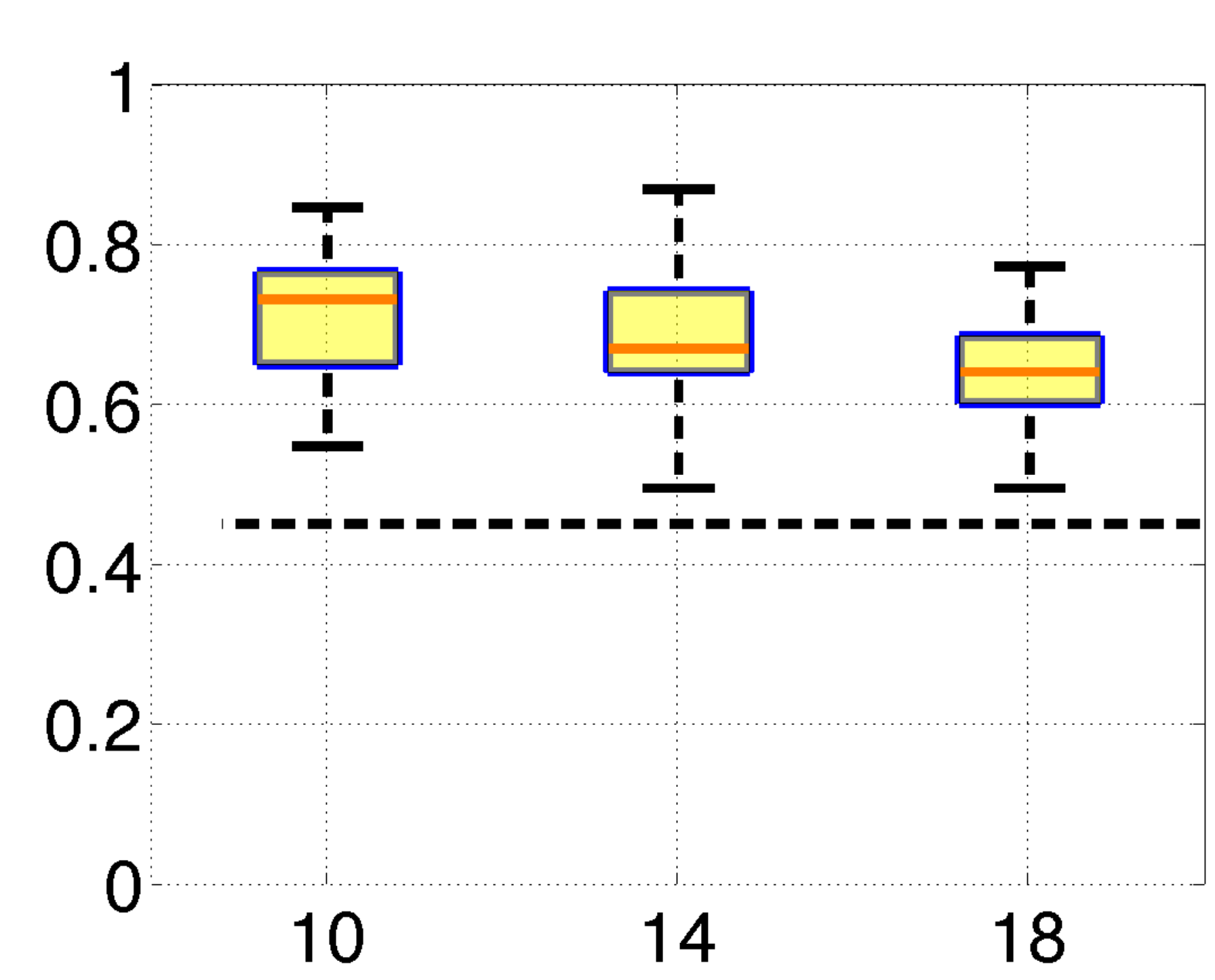}
 \includegraphics[scale=\lscale, clip=true, trim=\trleft cm \trbotbord cm \trright cm \trtop cm]{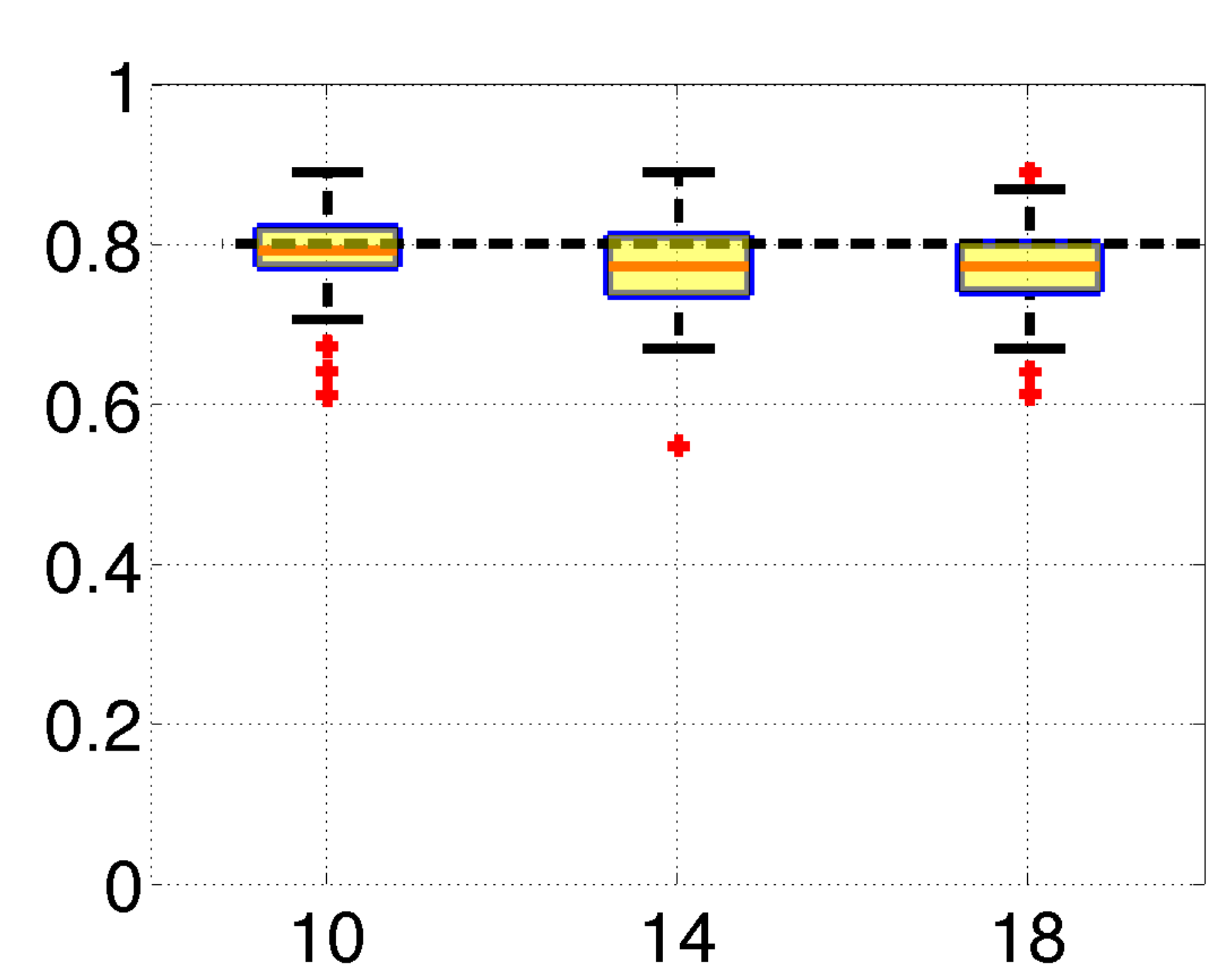}\\
 {\scriptsize $\hphantom{xxxx}\log_2 N$\hskip 1.9cm$\log_2 N$  \hskip 1.9cm$\log_2 N$}
   \vskip-2mm
 \caption{\textbf{Estimation performance of $\rho_{\rm x}$} as function of $\log_2 N$.
 \label{fig:perf_rho}}
 \end{figure}

The performance of $\hat \Theta^{M, BB}_N$, re-labelled $\widehat{\Theta}^{M}$ for simplicity, is compared against other existing estimation procedures.
The scaling exponents $(h_1, h_2)$ are estimated by means of the univariate wavelet based estimator for Hurst parameter, $(\widehat{h}^{\rm U}_1,\widehat{h}^{\rm U}_2)$, as described in \cite{Veitch1999} and applied to each component independently:
The univariate estimate of $h_1$ (resp. $h_2$) is obtained by taking the minimum (resp. maximum) between the linear regression coefficients of $\log_2 S_{11}(2^j)$ and $\log_2 S_{22}(2^j)$ versus $j\in J^{\rm fs}$ (resp.\ $j\in J^{\rm cs}$), with $J^{\rm fs}=\{j_1,\ldots,\lfloor\frac{j_1+j_2}{2}\rfloor\}$ (fine scales) and $J^{\rm cs}=\{ \lfloor\frac{j_1+j_2}{2}\rfloor+1,\ldots,j_2\}$ (coarse scales).
The parameters $(h_1, h_2, \beta)$ are also estimated using the multivariate semiparametric estimator $(\widehat{h}^{\rm W}_1,\widehat{h}^{\rm W}_2,\widehat{\beta}^{\rm W})$ proposed in \cite{AbryDidier2015}, which relies on the multiscale eigenstructure of $S(2^j)$. The statistics $\widehat{\Theta}^{\rm M}$, $(\widehat{h}^{\rm W}_1,\widehat{h}^{\rm W}_2,\widehat{\beta}^{\rm W})$ and $(\widehat{h}^{\rm U}_1,\widehat{h}^{\rm U}_2)$ are compared in Figs.~\ref{fig:perf_h1} to \ref{fig:perf_rho} in yellow, blue and magenta colors, respectively.

Though the full parametrization of Biv-OfBm requires a 7-dimensional vector parameter $$\Theta=(h_1,h_2,\rho_{\rm x},\sigma_{\rm x_1},\sigma_{\rm x_2},\beta,\gamma),$$ for ease of exposition we focus only on the $5$ most interesting parameters $(h_1,h_2,\rho_{\rm x},\beta,\gamma)$.
This follows the univariate literature that focuses on the estimation of $H$ for fBm, while neglecting the less interesting parameter $\sigma^2$.

   \begin{figure}[t]
      \captionsetup{justification=centering}
   \begin{center}
{\scriptsize \hskip .9cm$\rho_{\rm x}=0.1$\hskip 1.7cm$\rho_{\rm x}=0.45$  \hskip 1.4cm$\rho_{\rm x}=0.8$} \\
\rotatebox{90}{\hskip .2cm\scriptsize$\beta=\gamma=0$} \includegraphics[scale=\lscale, clip=true, trim=\trleftbord cm \trbot cm \trright cm \trtop cm]{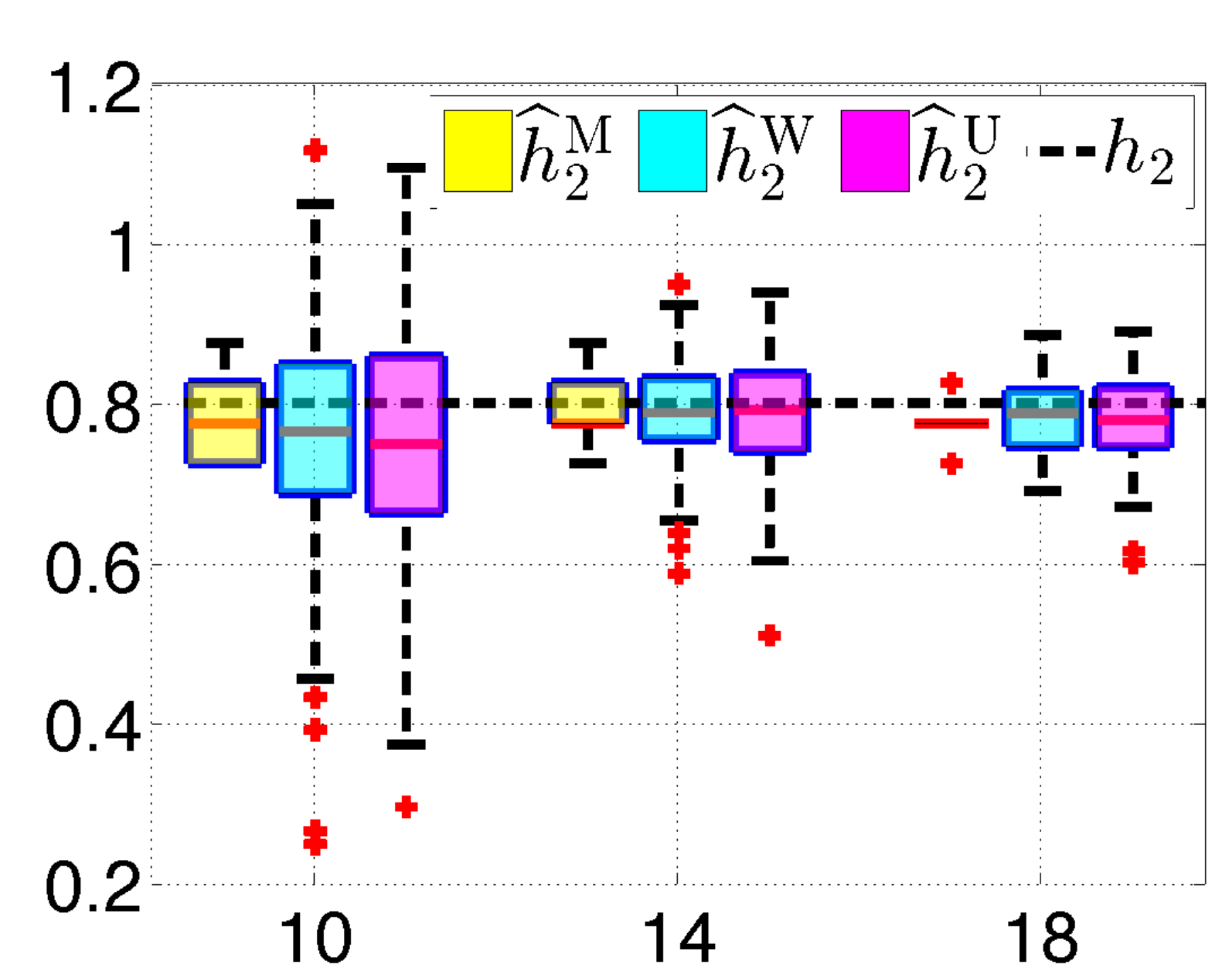}
  \includegraphics[scale=\lscale, clip=true, trim=\trleft cm \trbot cm \trright cm \trtop cm]{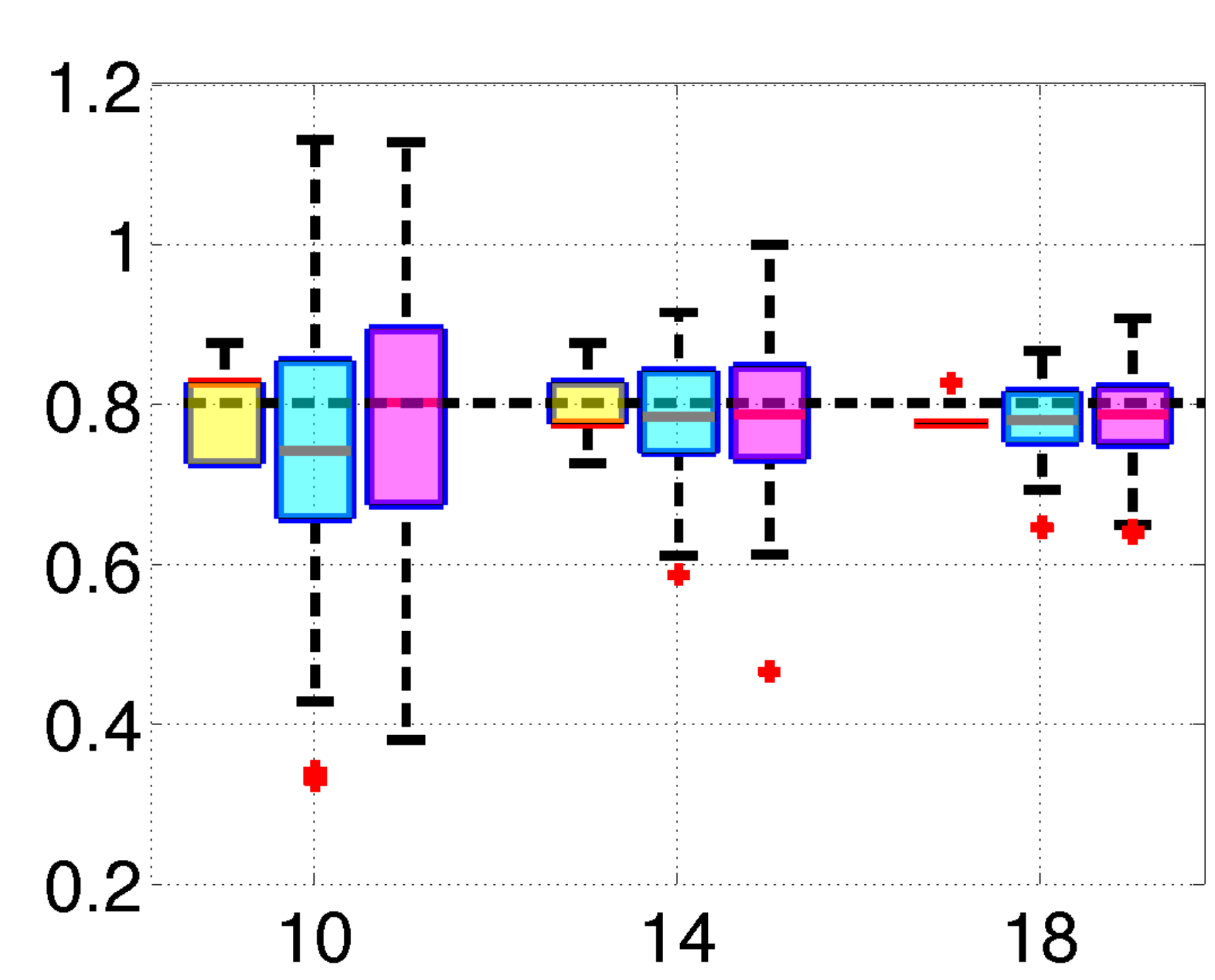}
 \includegraphics[scale=\lscale, clip=true, trim=\trleft cm \trbot cm \trright cm \trtop cm]{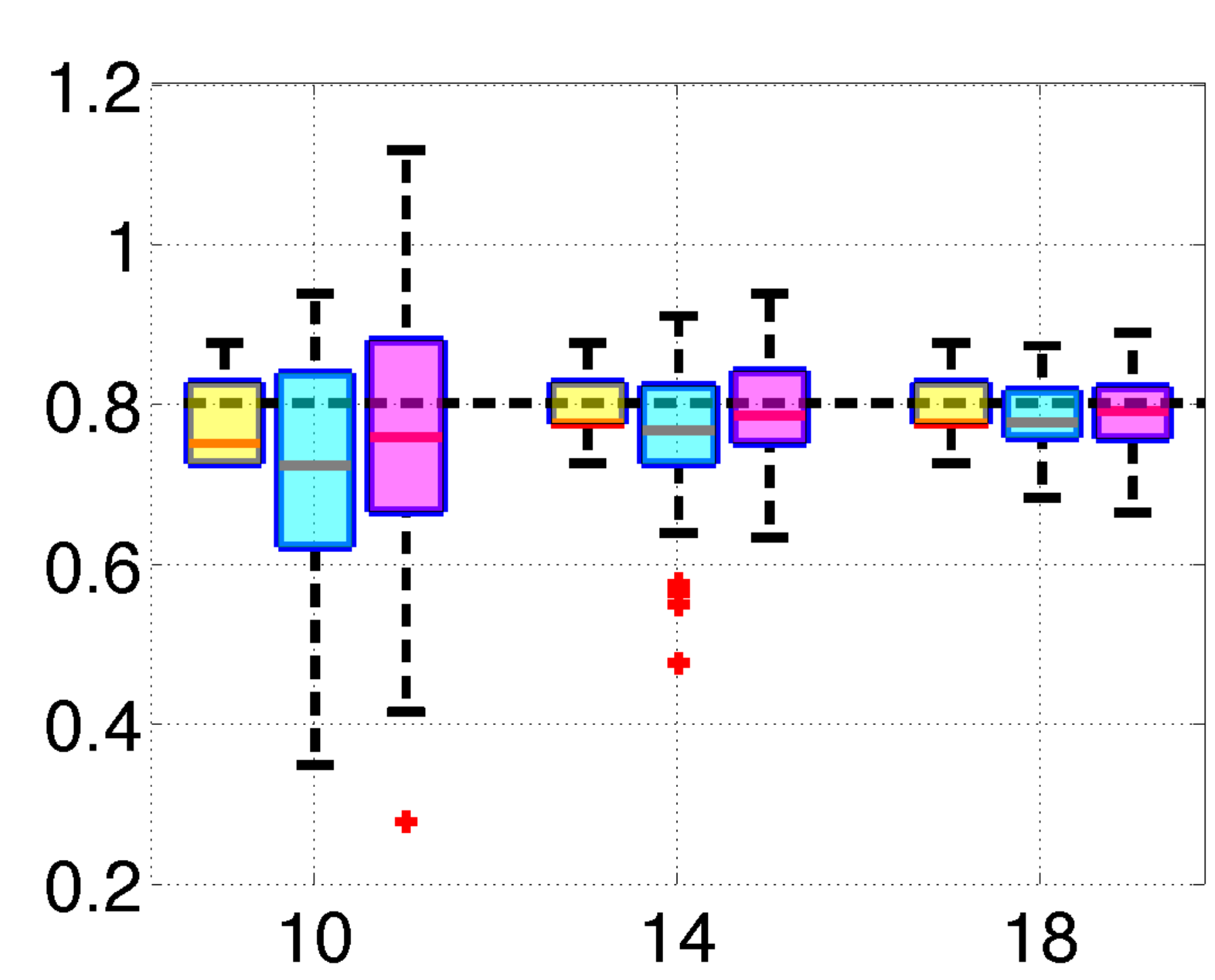}\\ \vskip -.3cm
 \rotatebox{90}{\hskip .2cm\scriptsize$\beta=\gamma=0.5$} \includegraphics[scale=\lscale, clip=true, trim=\trleftbord cm \trbot cm \trright cm \trtop cm]{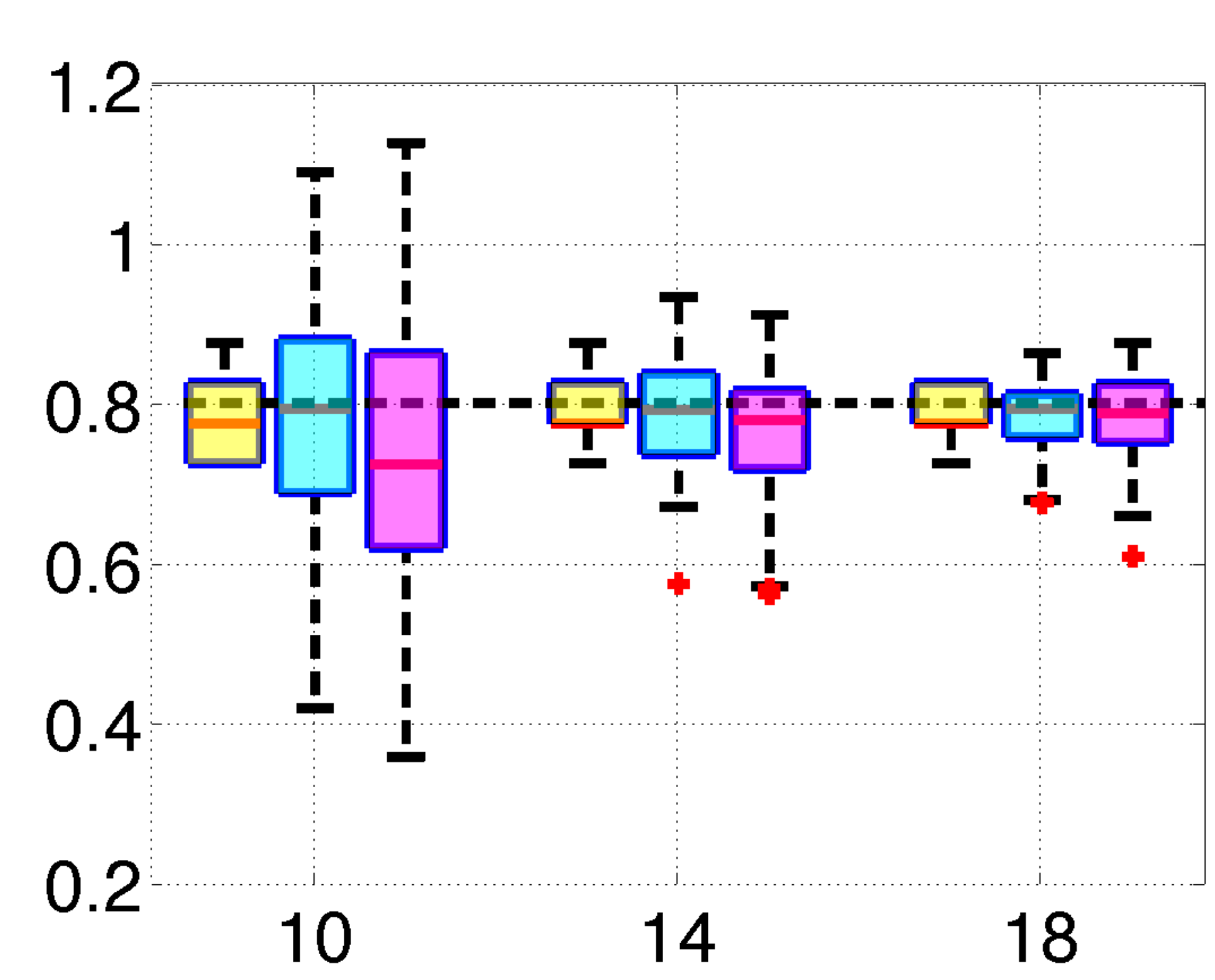}
  \includegraphics[scale=\lscale, clip=true, trim=\trleft cm \trbot cm \trright cm \trtop cm]{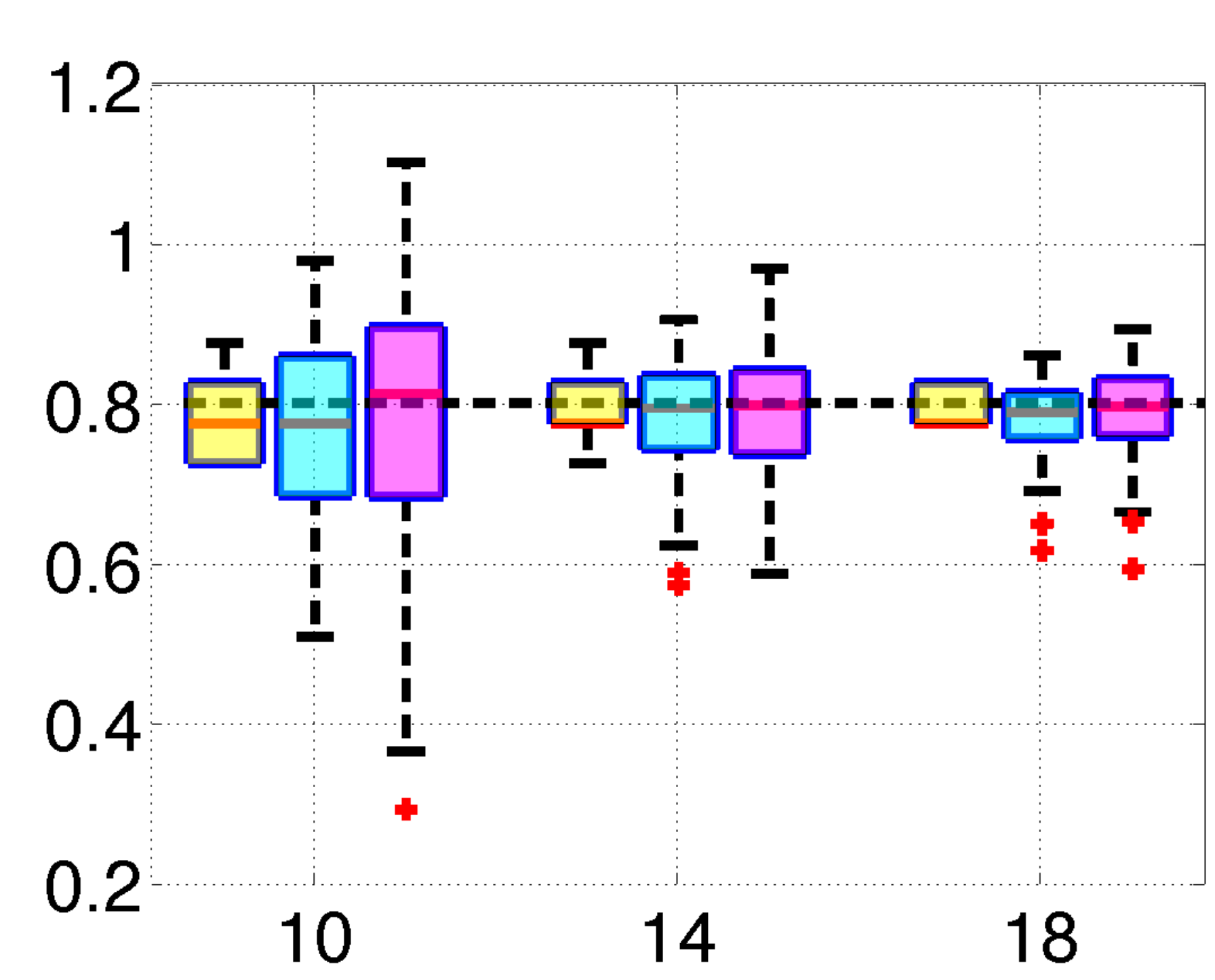}
 \includegraphics[scale=\lscale, clip=true, trim=\trleft cm \trbot cm \trright cm \trtop cm]{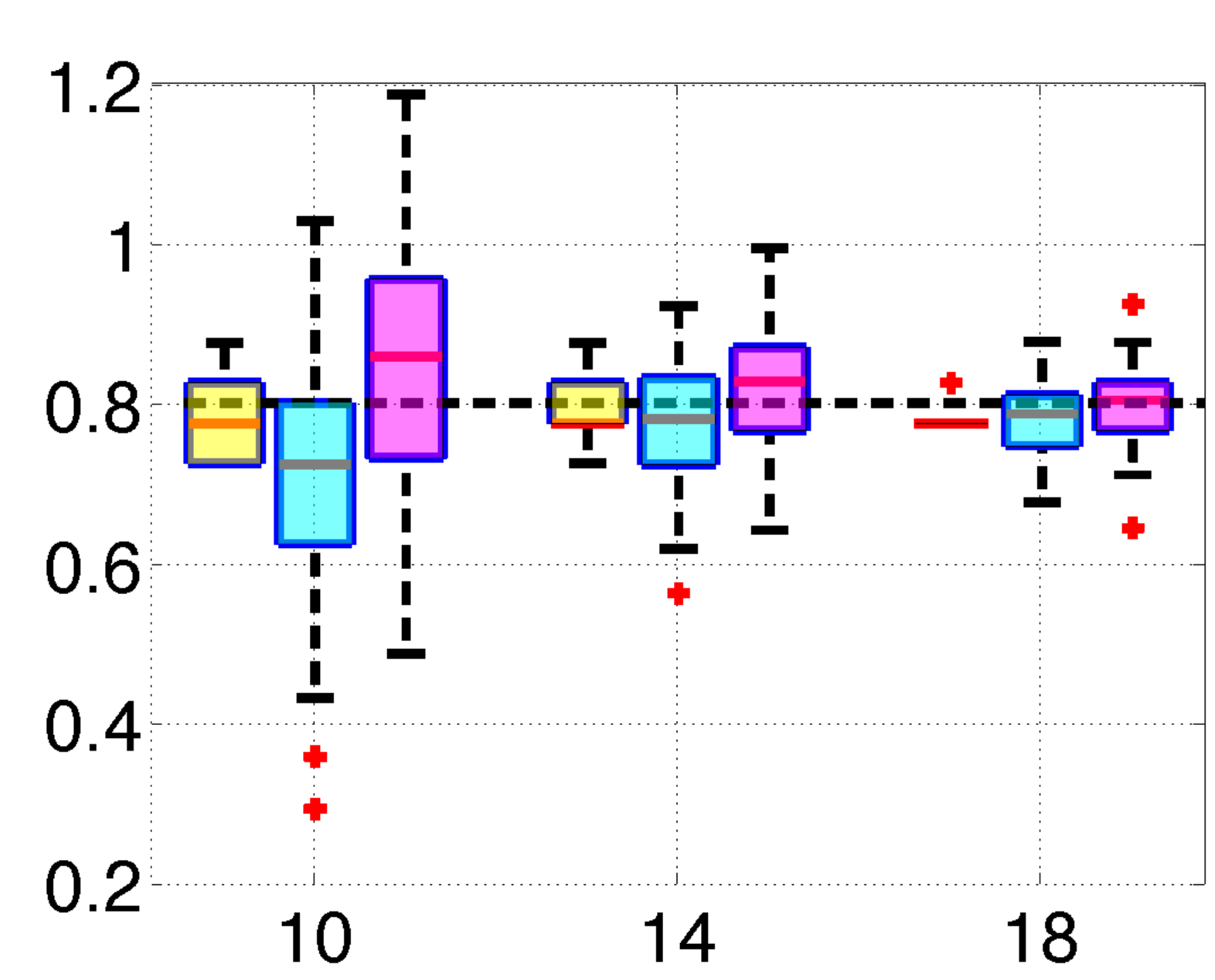}\\ \vskip -.3cm
 \rotatebox{90}{\hskip .1cm\scriptsize$\beta=-\gamma=0.5$} \includegraphics[scale=\lscale, clip=true, trim=\trleftbord cm \trbotbord cm \trright cm \trtop cm]{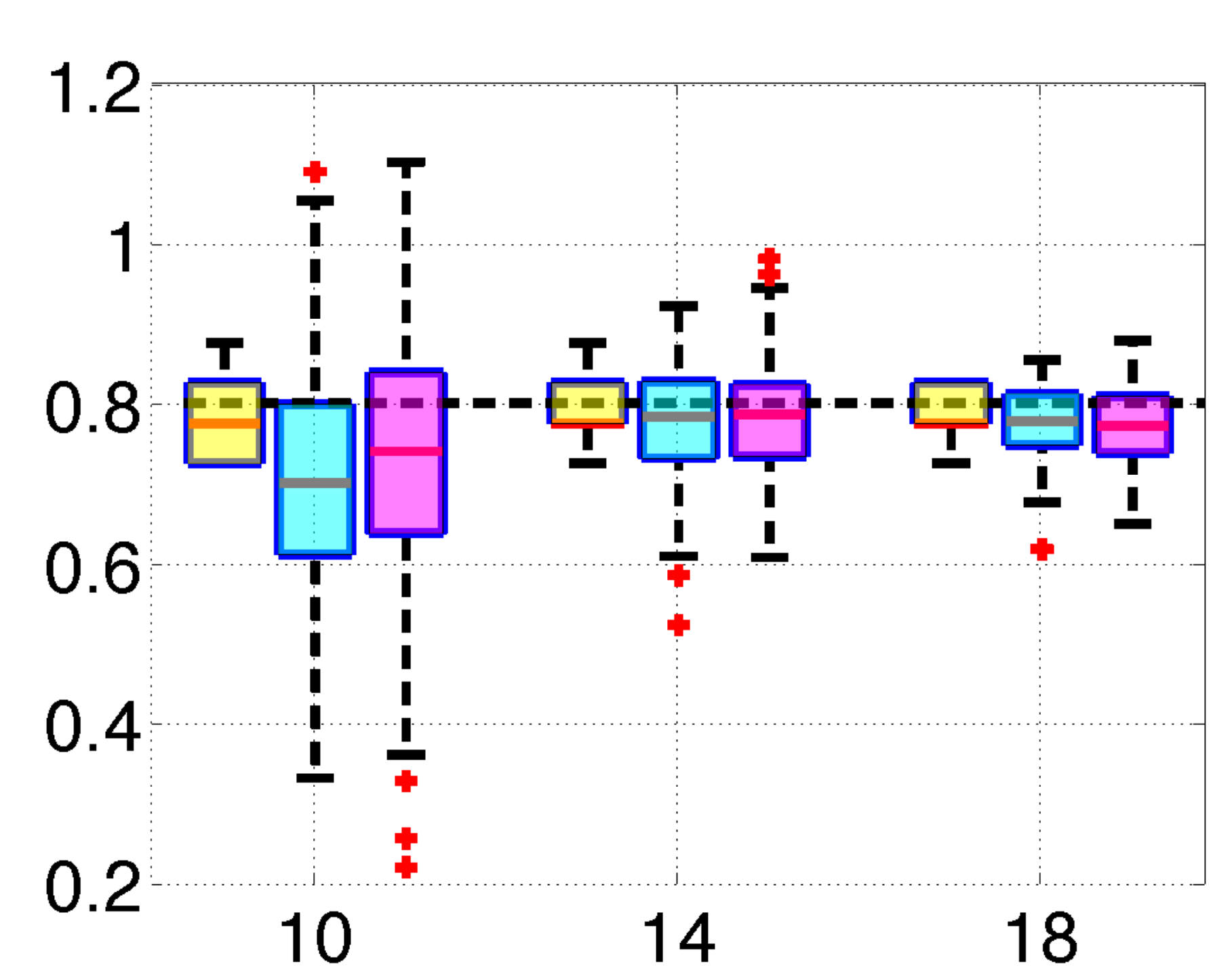}
  \includegraphics[scale=\lscale, clip=true, trim=\trleft cm \trbotbord cm \trright cm \trtop cm]{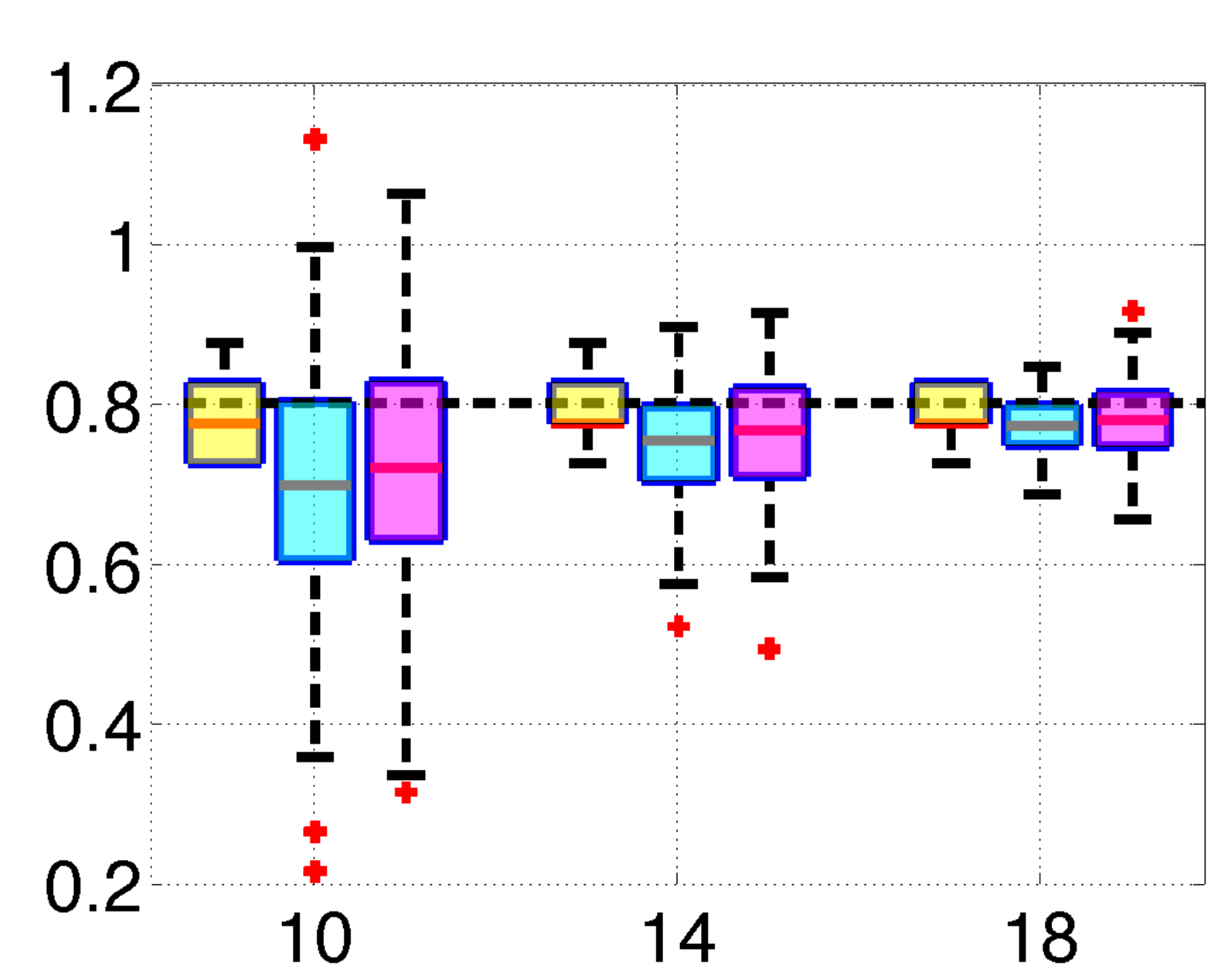}
 \includegraphics[scale=\lscale, clip=true, trim=\trleft cm \trbotbord cm \trright cm \trtop cm]{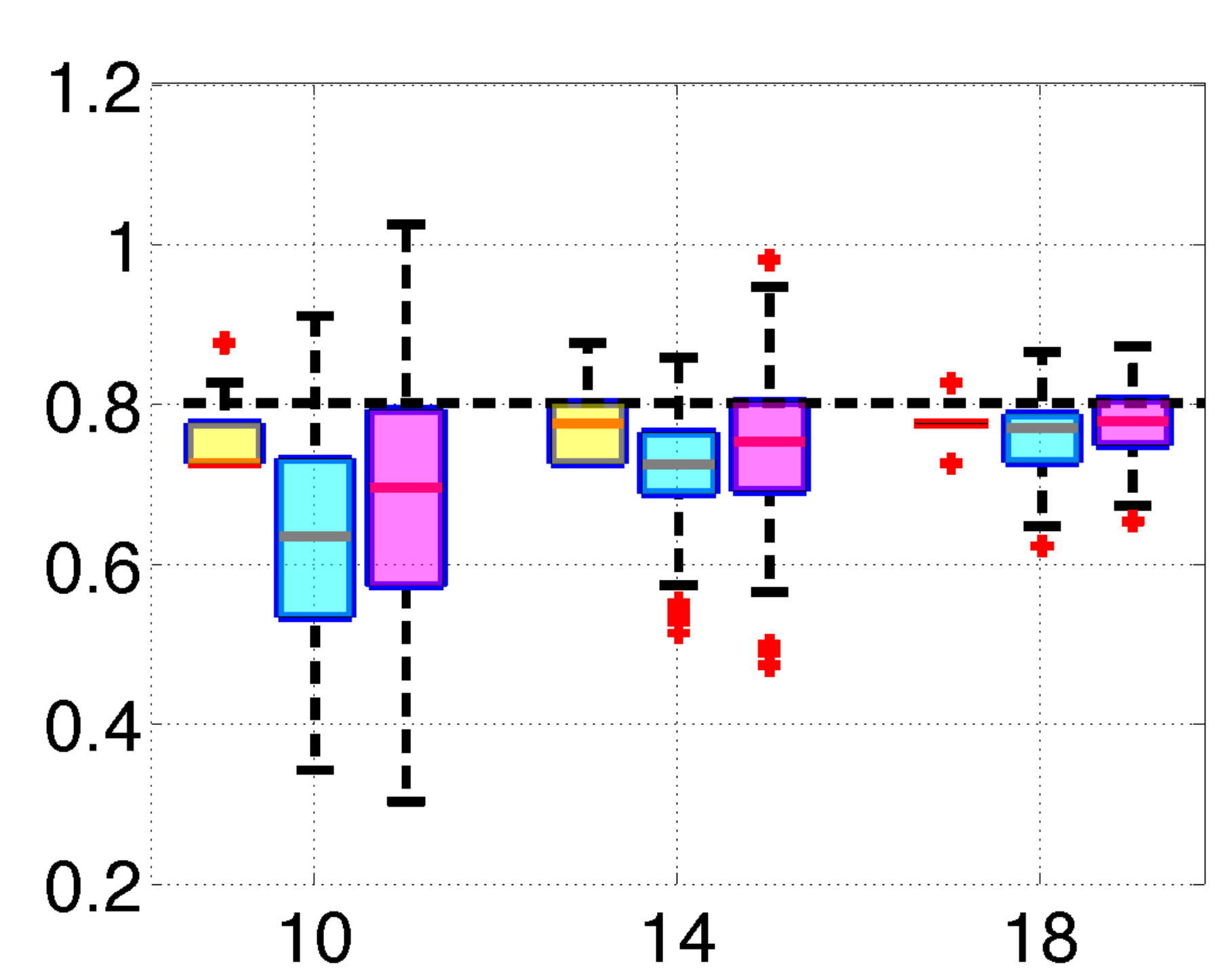}\\
 {\scriptsize $\hphantom{xxxx}\log_2 N$\hskip 1.9cm$\log_2 N$  \hskip 1.9cm$\log_2 N$}
   \vskip-2mm
 \caption{\textbf{Estimation performance of $h_2$} as function of $\log_2 N$.
  \label{fig:perf_h2}}
  \end{center}
 \end{figure}

 \begin{figure}[t]
  \centering
 {\scriptsize \hskip 1.3cm$\beta=\gamma=0$ \hskip 1.4cm $\beta=\gamma=0.5$ \hskip 1.1cm $\beta=-\gamma=0.5$}\\
\rotatebox{90}{\hskip .3cm\scriptsize Execution time}$\;$\includegraphics[scale=.155, clip=true,trim=0cm 0cm .3cm 0cm]{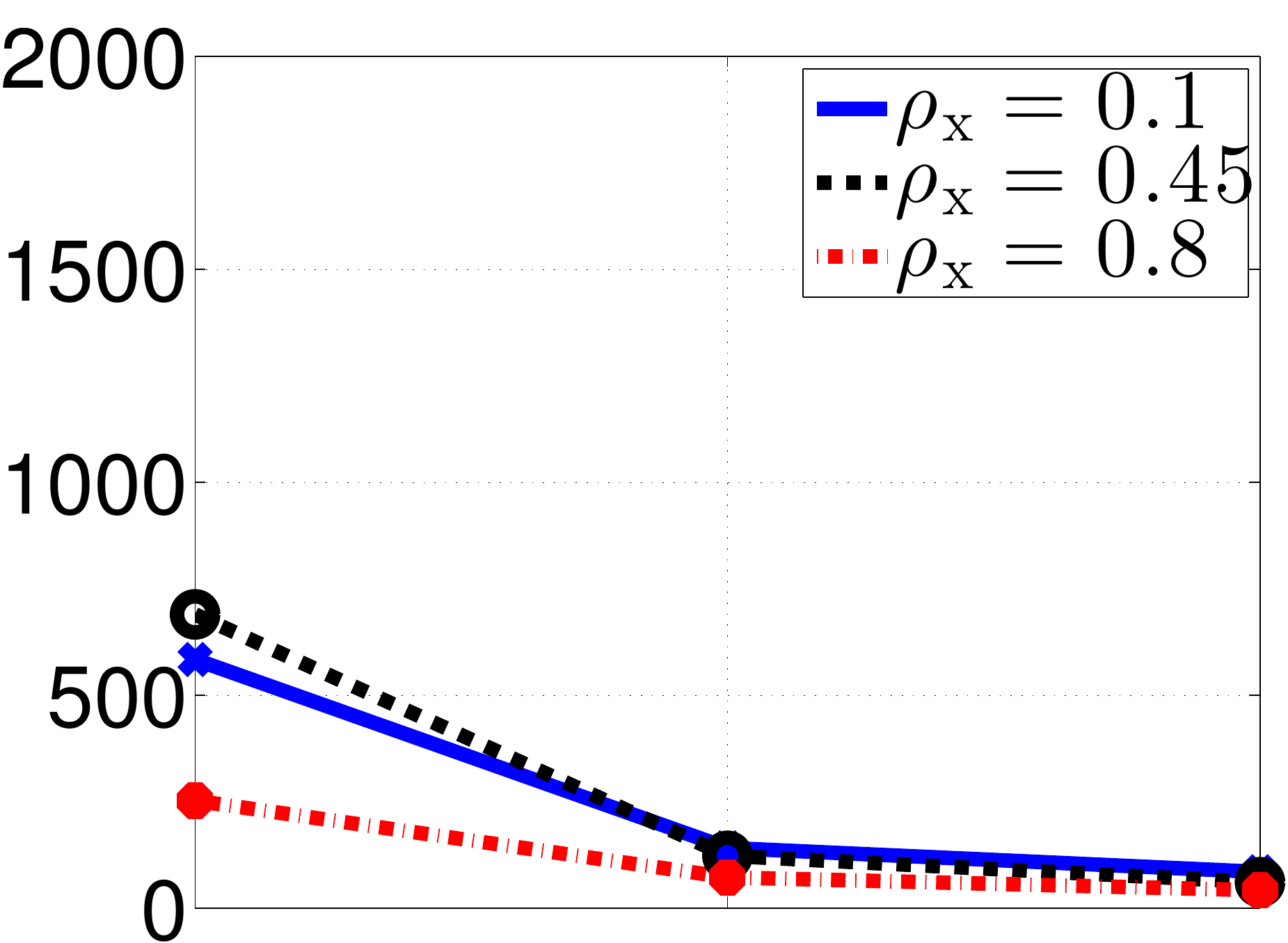}
\includegraphics[scale=.155, clip=true,trim=-1cm 0cm .3cm 0cm]{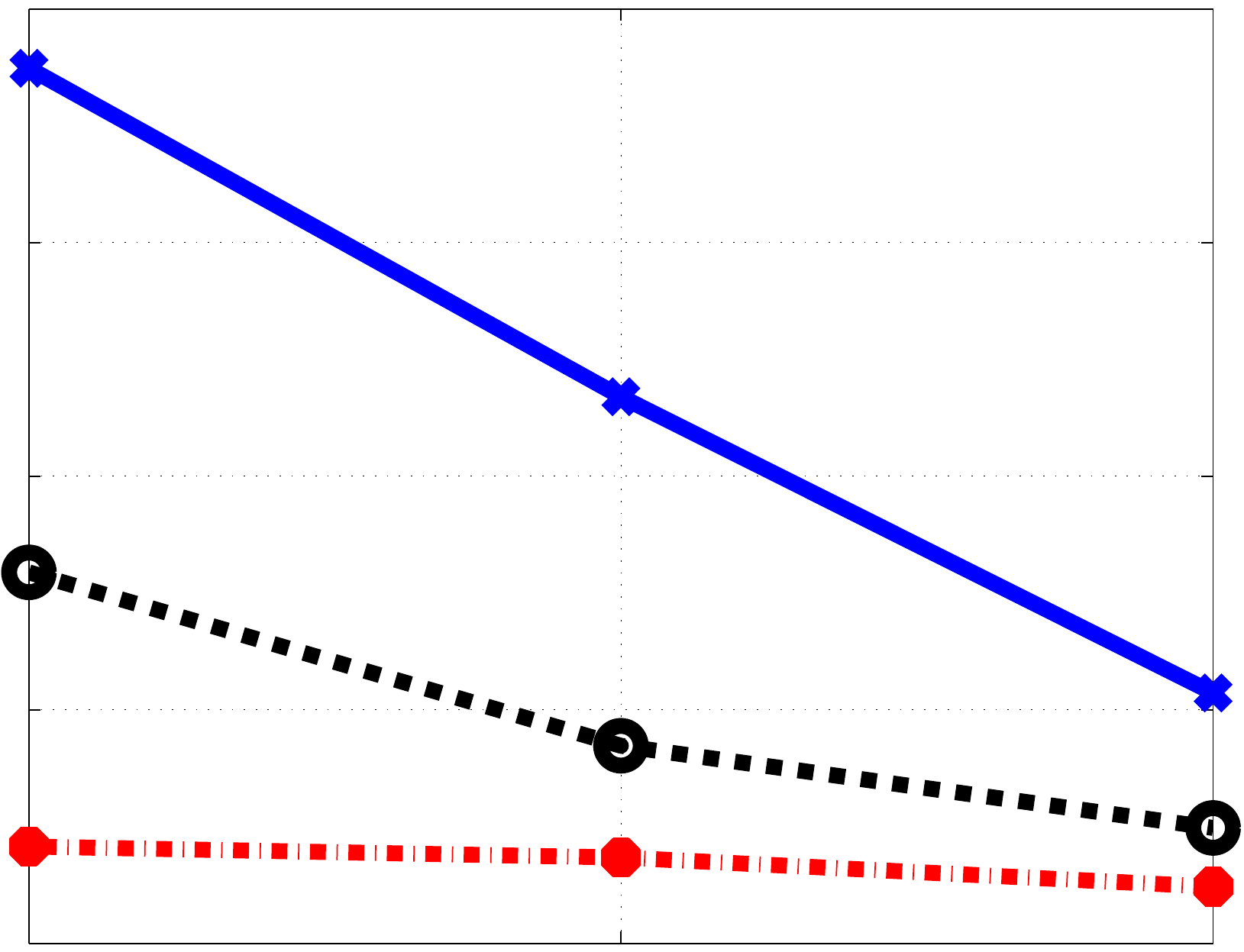}
\includegraphics[scale=.155, clip=true,trim=-1cm 0cm .3cm 0cm]{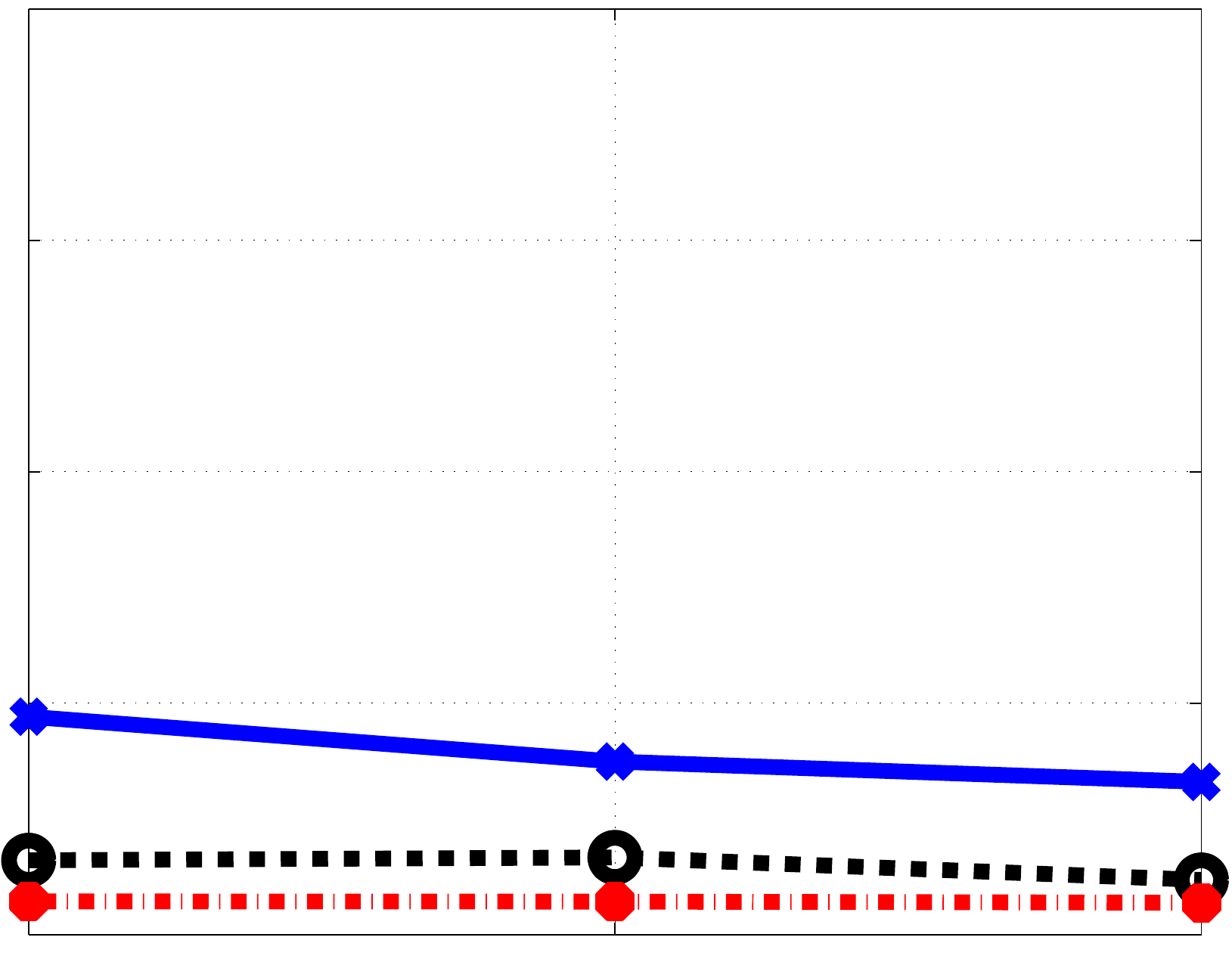}\\
\rotatebox{90}{\hskip -.1cm\scriptsize Percentage of iterations}$\;$\includegraphics[scale=.155, clip=true,trim=-.6cm 0cm .2cm 0cm]{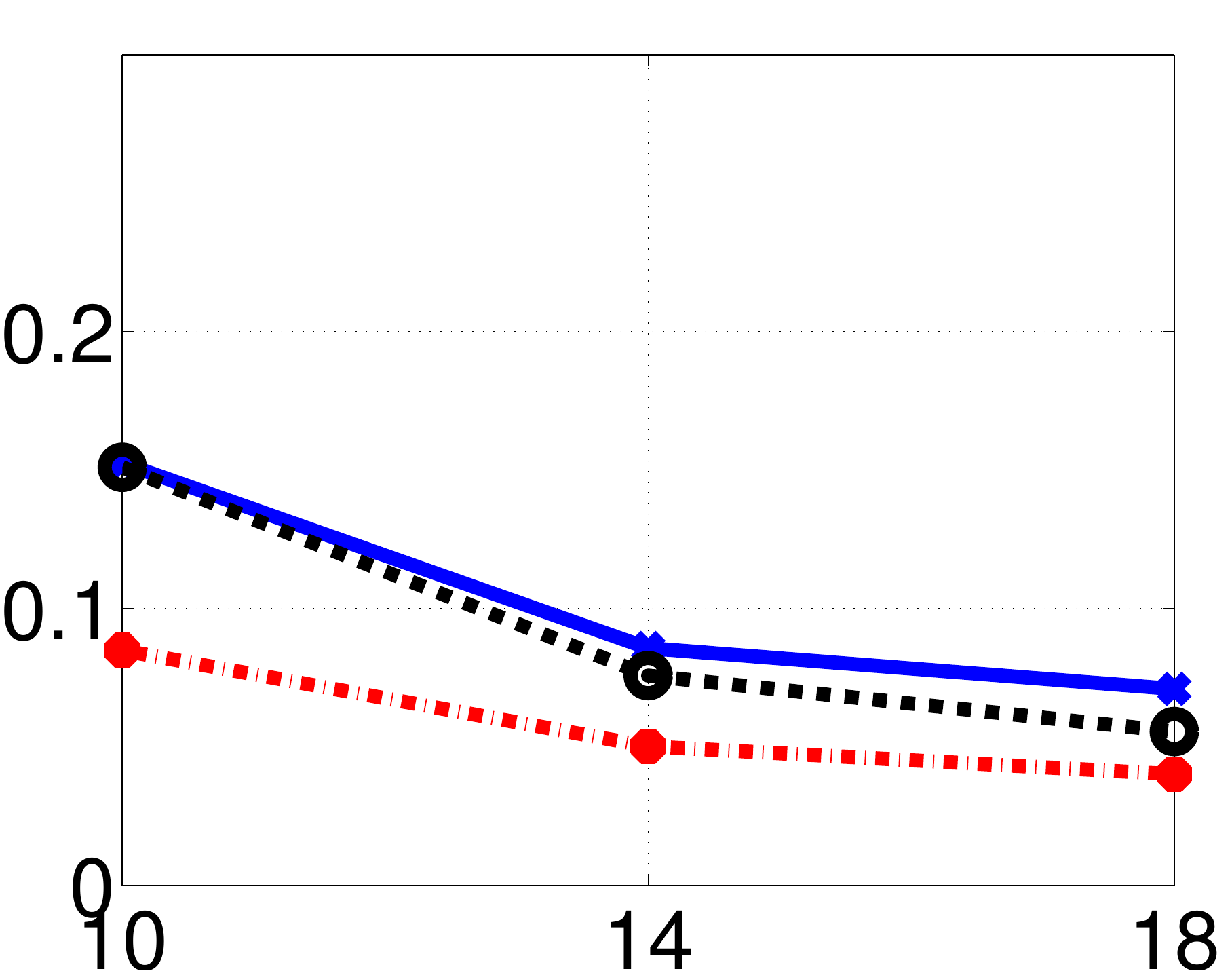}
\includegraphics[scale=.155, clip=true,trim=0cm 0cm .2cm 0cm]{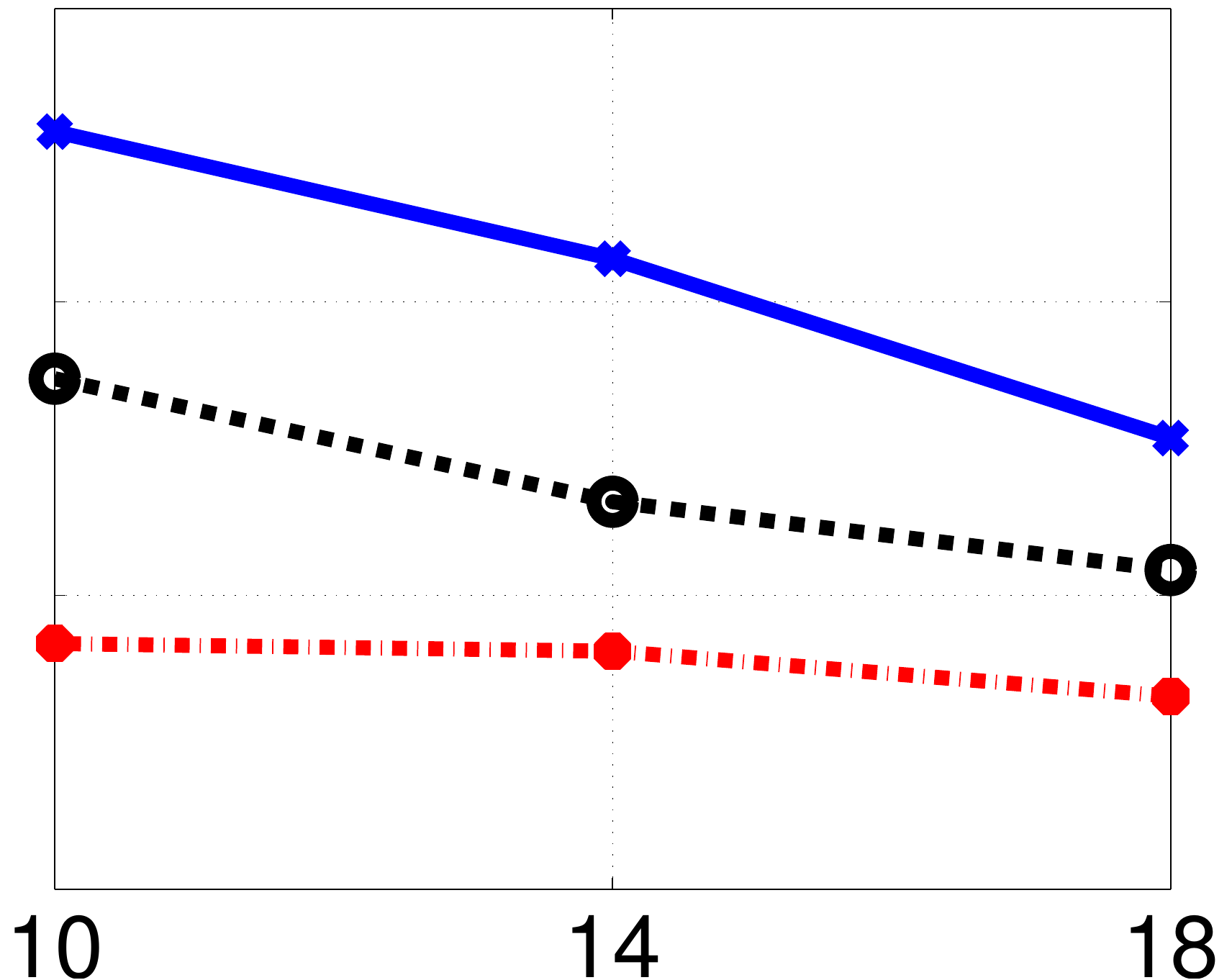}
\includegraphics[scale=.155, clip=true,trim=0cm 0cm .2cm 0cm]{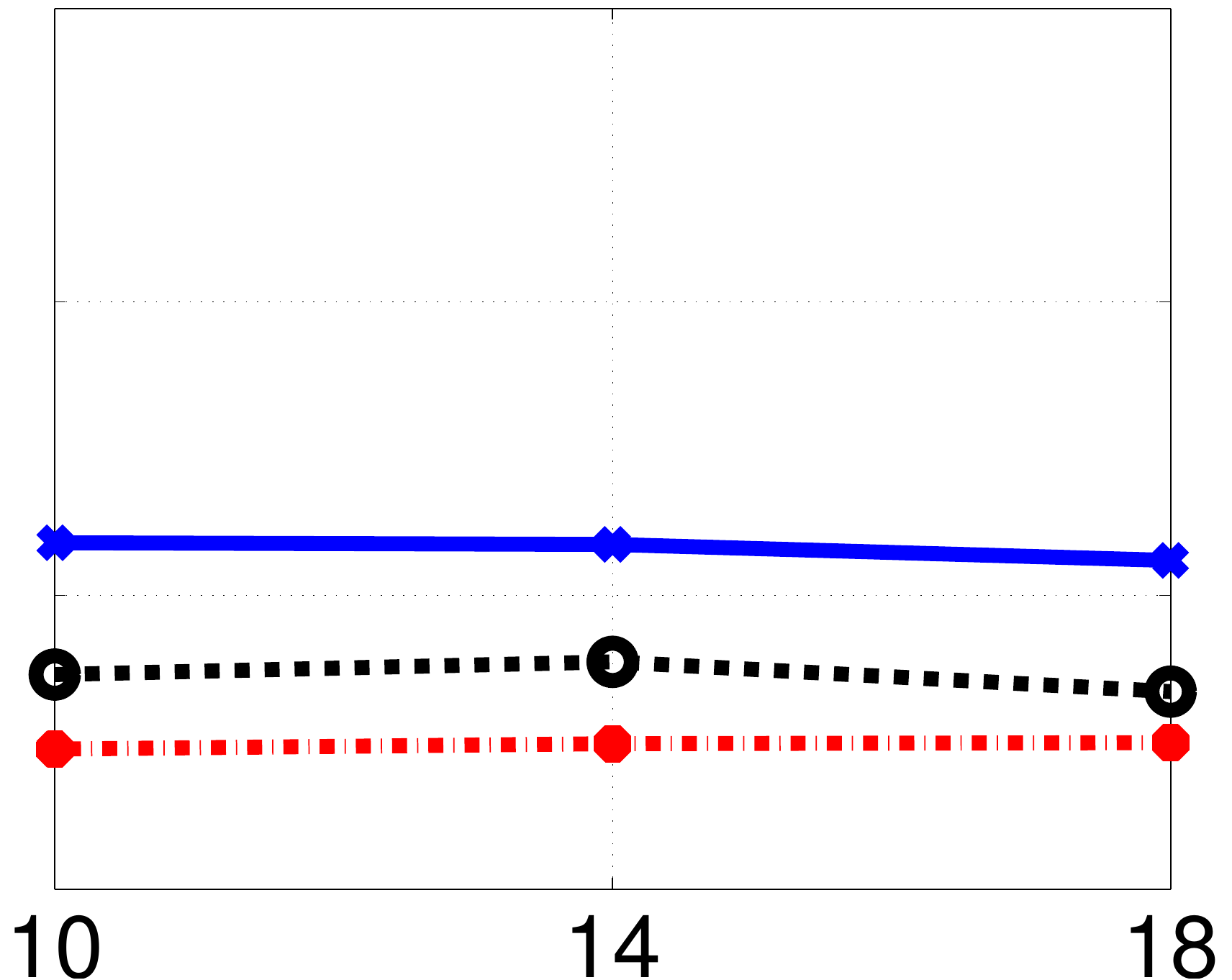}\\
 {\scriptsize $\hphantom{xxxx}\log_2 N$\hskip 1.9cm$\log_2 N$  \hskip 1.9cm$\log_2 N$}
  \vskip-2mm
\caption{\textbf{Computational costs as functions of $\log_2 N$.}  Top: Average execution time (in seconds). Bottom: percentage of performed iterations (compared to a greedy grid search algorithm with same accuracy).\label{fig:time}}
 \end{figure}

 \subsection{Estimation performance}

\noindent\textbf{Estimation of the dominant scaling exponent $h_2$ (Fig.~\ref{fig:perf_h2}).}
As expected, all procedures yield accurate estimates of the largest scaling exponent $h_2$. While all methods show comparable performance for large sample sizes, it is interesting that $\widehat{h}_2^{\rm M}$ displays better performance with lower bias and dispersion for small sample sizes by comparison to $\widehat{h}_2^{\rm W}$ and $\widehat{h}_2^{\rm U}$. The impact of the correlation $\rho_{\rm x}$ or the mixing parameters $(\beta, \gamma)$ on the performance is weak. \\

\noindent\textbf{Estimation of the non-dominant scaling exponent $h_1$ (Fig.~\ref{fig:perf_h1}).}
Estimating the lowest scaling exponent $h_1$ is intrinsically more difficult because the mixture of power laws masks the non-dominant Hurst eigenvalue.
As expected, univariate-type analysis fails to estimate correctly $h_1$ (except when there is no mixing $(\beta=\gamma=0)$).
While $\widehat{h}_1^{\rm M}$ and $\widehat{h}_2^{\rm W}$ show essentially the same performance for large sample size, it is interesting to note that $\widehat{h}_1^{\rm M}$ displays a far superior performance with lower bias and dispersion for small sample sizes.
However, the bias of $\widehat{h}_1^{\rm M}$ for $(\beta=-\gamma)$ and small correlation $\rho_{\rm x}=0.1$ are observed, showing that $\widehat{h}_1^{\rm M}$ is more strongly affected by the conjunction of low correlation amongst components and anti-orthogonal mixing than $\widehat{h}_2^{\rm M}$.\\

\noindent\textbf{Estimation of $\beta$ (Fig.~\ref{fig:perf_beta}).} A significant benefit of $\widehat{\beta}^{\rm M}$ consists of its being robust to small sample sizes, when $\widehat{\beta}^{\rm W}$ is not. While the mixing parameters seem not to impact the performance of $\widehat{\beta}^{\rm M}$, a low correlation value $\rho_{\rm x}$ hurts its performance. Other results not reported here for reasons of space also show that the performance of $\widehat{\beta}^{\rm M}$ is robust to a decrease of $h_2-h_1$, while that of $\widehat{\beta}^{\rm W}$ drastically deteriorates when $h_2-h_1 \rightarrow 0$.\\

\noindent\textbf{Estimation of $\gamma$ (Fig.~\ref{fig:perf_gamma}).} The performance of $\widehat{\gamma}^{\rm M}$ is very satisfactory, yet it is observed to be affected by low correlation. \\

\noindent\textbf{Estimation of $\rho_{\rm x}$ (Fig.~\ref{fig:perf_rho}).} The parameter $\rho_{\rm x}$ appears the most difficult to estimate, with significant bias for low correlation and anti-orthogonal mixing, a result in consistency with \cite{Amblard_P-0_2011_j-ieee-tsp_imfbm}.

 \subsection{Computational costs}

\begin{figure}[!h]
\centerline{
\includegraphics[width=0.32\linewidth]{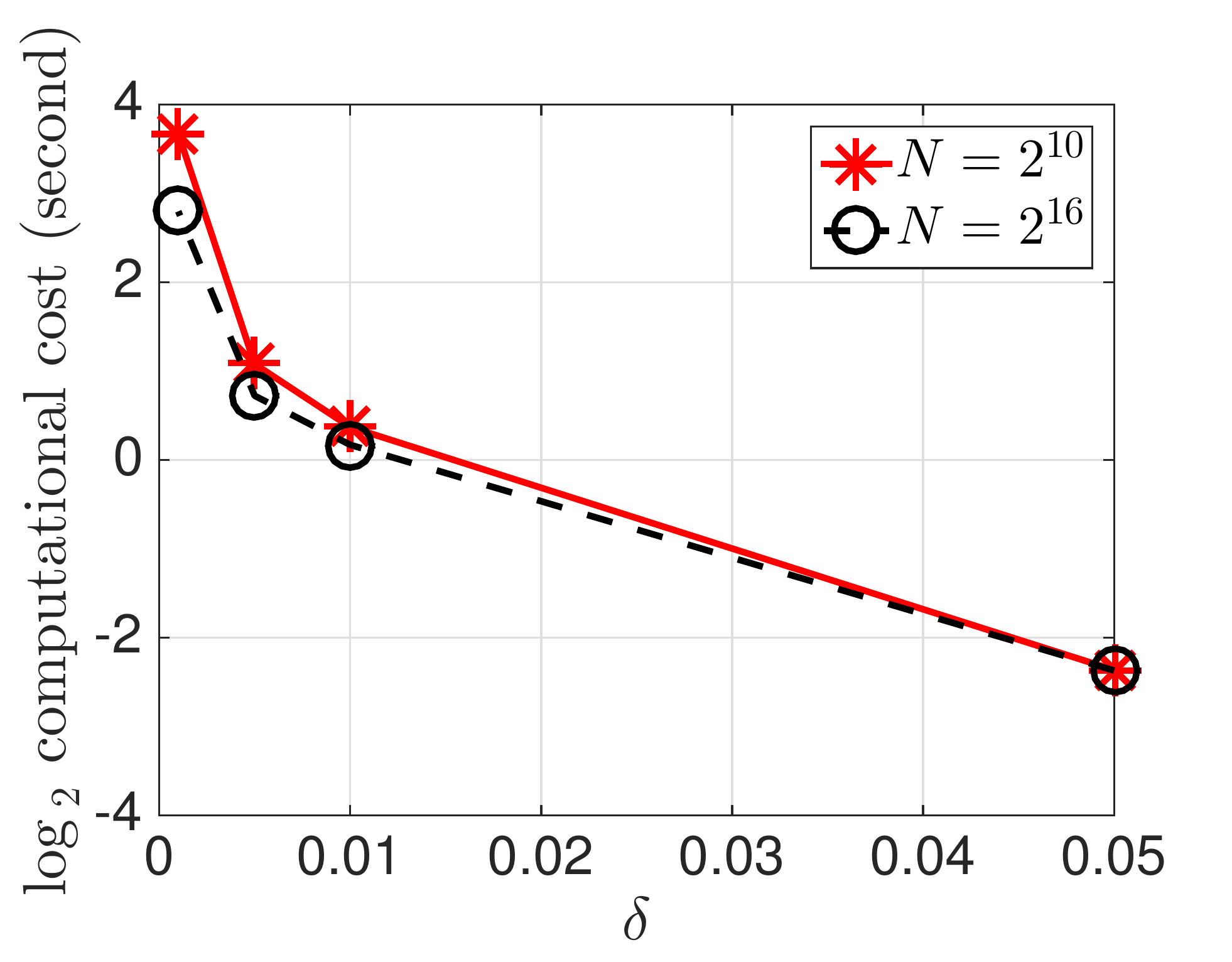}
\includegraphics[width=0.32\linewidth]{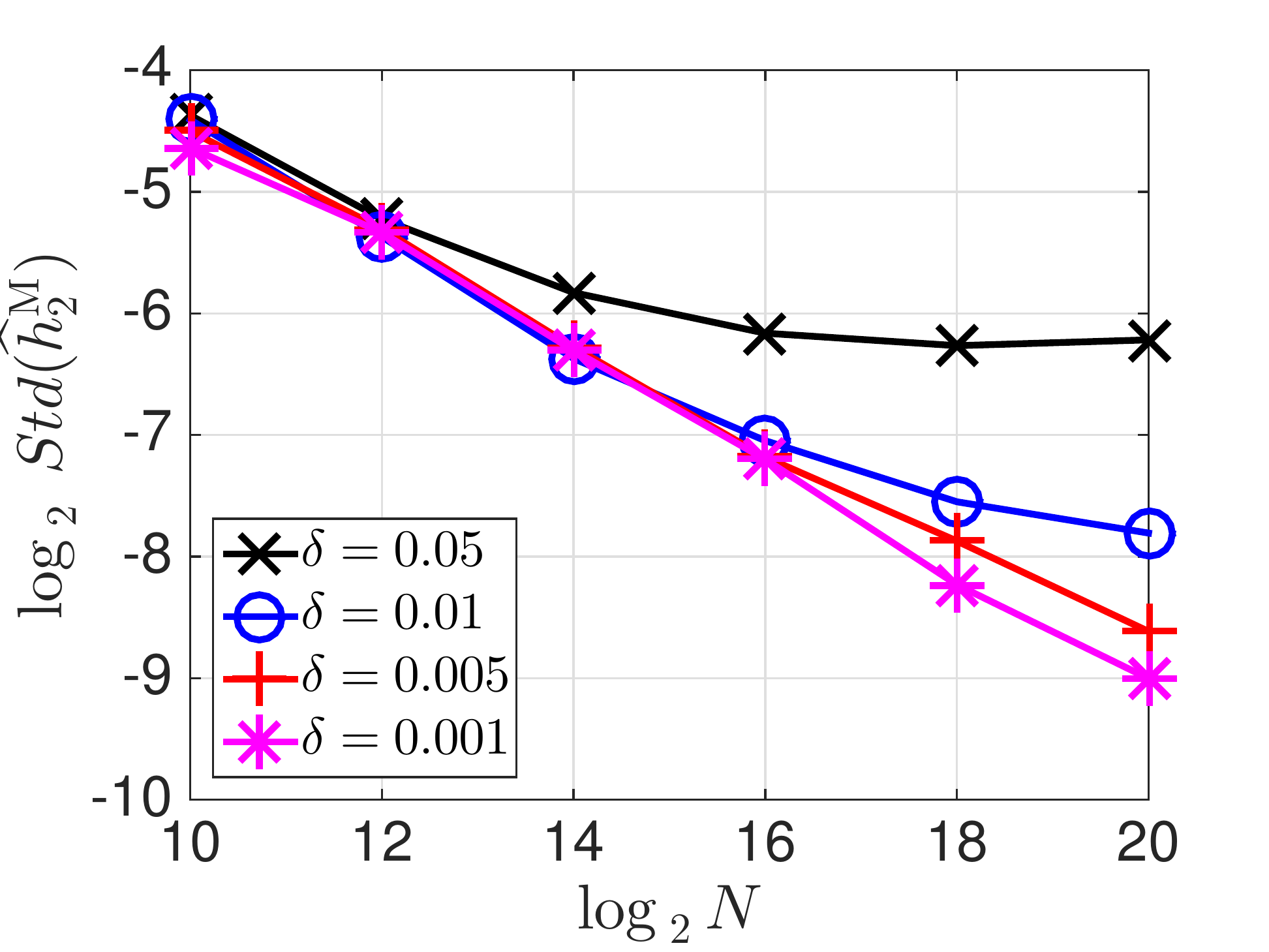}
}
\vskip-2mm
\caption{{\bf Precision.} Computational costs (left) and dispersion of the estimates, as functions of targeted precision $\delta$ \label{fig:preccosta} }
\end{figure}

The computational costs of the proposed identification algorithm described in Section~\ref{sec:algoBB} are reported in Fig.~\ref{fig:time} (top plots) as a function of $\log_2 N$ for each parameter setting. They are significantly smaller than those required by a systematic greedy grid search.
Fig.~\ref{fig:time} also clearly shows that the stronger the correlation amongst components, the easier the minimization of the functional $C_N$, thus indicating that the cross terms play a significant role in the identification of Biv-OfBm. Though surprising at first, the clear decrease of the computational costs with the increase of the sample size may be interpreted as the fact that it is obviously far easier to disentangle three different power laws when a large number of scales $2^j$ is available, which then requires large sample sizes. It is also worth noting that the orthogonal mixing, which may intuitively be thought of as the easiest, appears to be the most demanding in terms of iterations to minimize $C_N$. Unsurprisingly, the computational cost increases when the requested precision $\delta $ on the estimates is increased ($\delta \rightarrow 0$) Fig.~\ref{fig:preccosta} (left).

 \subsection{Sample size versus precision.}
 Figs.~\ref{fig:perf_h2} to \ref{fig:perf_rho} show that increasing the sample size $N$ improves the performance of $(\widehat{h}^{\rm U}_1,\widehat{h}^{\rm U}_2)$ and $(\widehat{h}^{\rm W}_1,\widehat{h}^{\rm W}_2,\widehat{\beta}^{\rm W})$, as both their median-bias and variance decrease with $N$.
Fig.~\ref{fig:preccosta} (right) indicates that the impact of $N$ is slightly more involved \GD{more involved?} for $\widehat{\Theta}^{\rm M}$.
As long as the dispersions of the estimates remain above the desired precision $\delta$ (targeted independently of $N$), one observes an expected decrease in the bias and dispersion when $N$ is increased. Because the minimization is stopped when the prescribed accuracy $\delta$ is reached, increasing $N$ without decreasing $\delta $ does not bring about any performance improvement, showing that the precision should be decreased when the sample size increases (empirically as $N^{-1/4}$) to improve the performance.
\vskip-4mm

\subsection{Asymptotic normality of $\hat \Theta$}

\begin{figure}[!h]
\centerline{
\includegraphics[width=0.32\linewidth]{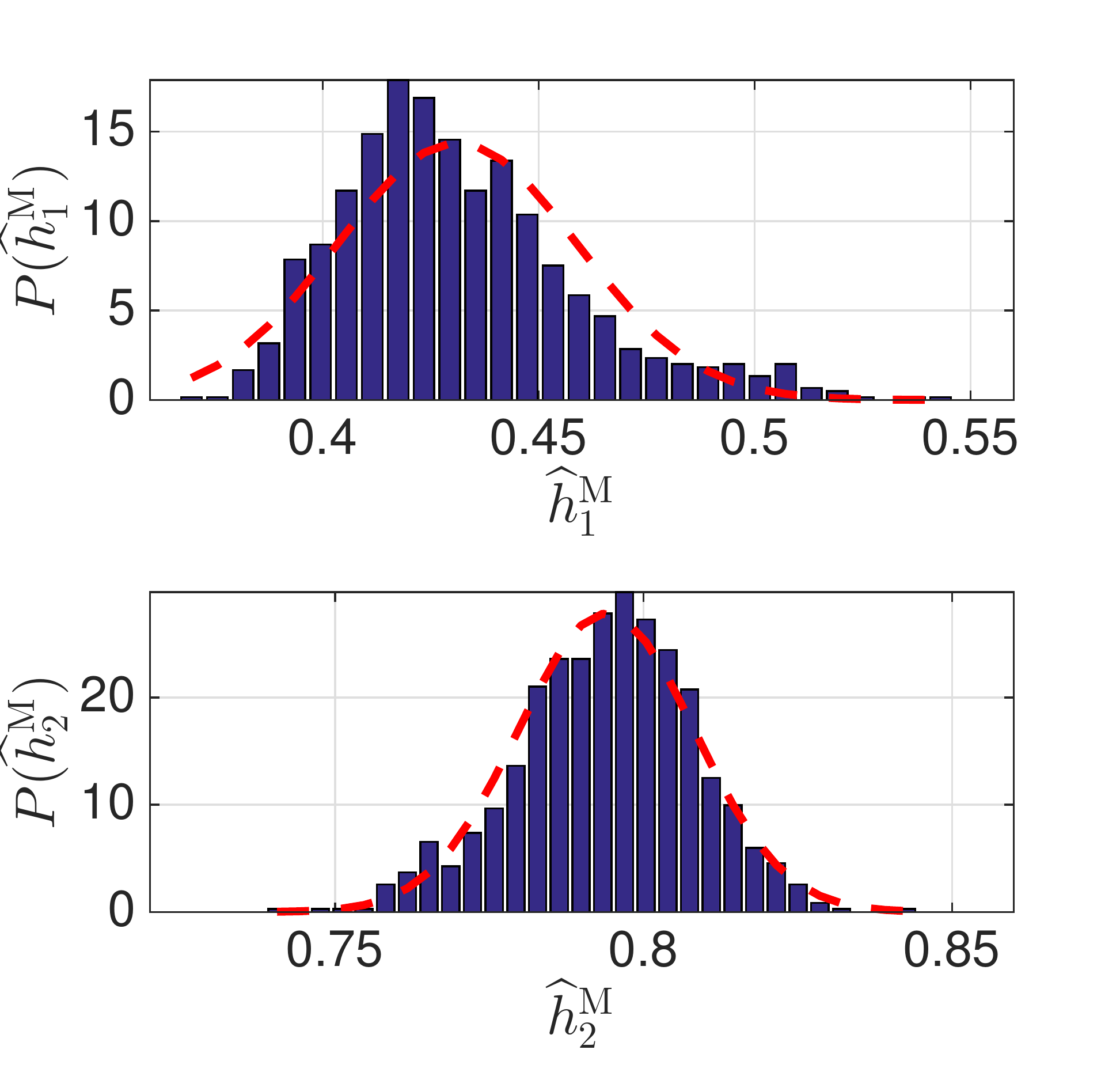}
\includegraphics[width=0.32\linewidth]{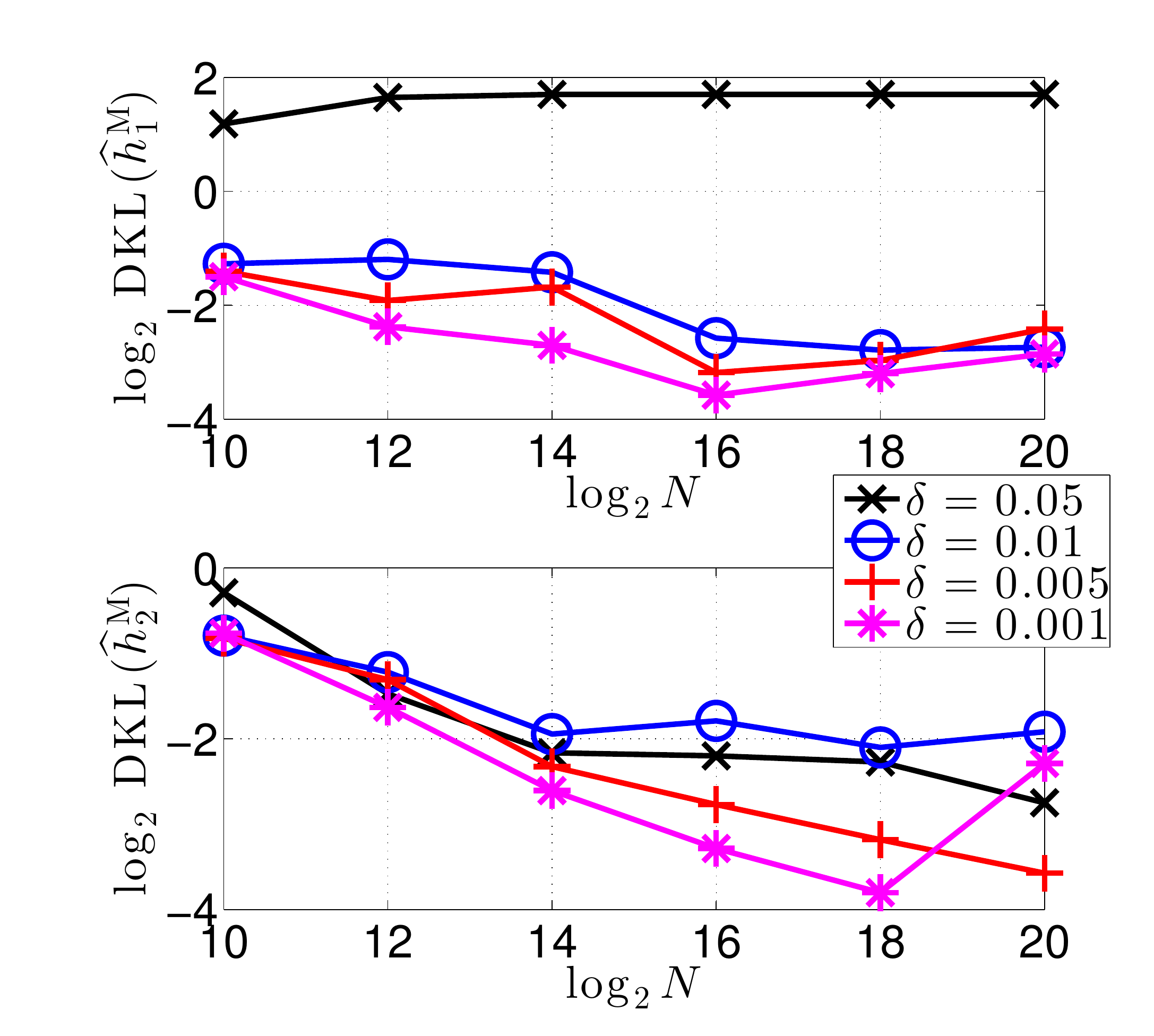}
}
\vskip-4mm
\caption{{\bf Asymptotic Normality.}  Left: empirical distributions of $\widehat{h}_1^{\rm M}$, $\widehat{h}_2^{\rm M}$ and best Gaussian fit, $N=2^{20}$, $\delta = 10^{-3}$.  Right: Kullback-Leibler divergence between the empirical distributions of $\widehat{h}_1^{\rm M}$, $\widehat{h}_2^{\rm M}$ and their best Gaussian fits, as functions of the sample size $N$. \label{fig:asna}}
\end{figure}

For a thorough study of normality,  $\Theta$ is restricted to (${h}_1$, ${h}_2$), while all other parameters are fixed \emph{a priori} and known.
Averages are obtained over 1000 independent realizations of Biv-OfBm for each sample size. Fig.~\ref{fig:asna} (left) visually compares the empirical distributions of $\widehat{h}_1^{\rm M}$, $\widehat{h}_2^{\rm M}$ to their best Gaussian fits. Fig.~\ref{fig:asna} (right) measures the Kullback-Leibler divergence between the empirical distributions of $\widehat{h}_1^{\rm M}$, $\widehat{h}_2^{\rm M}$ and their best Gaussian fits, as a function of the sample size $N$. Fig.~\ref{fig:asna} confirms the asymptotic normality of the estimates $\widehat{h}_1^{\rm M}$, $\widehat{h}_2^{\rm M}$ as theoretically predicted in Section~\ref{sec.Co}. It further shows that the higher the requested precision $\delta$, the faster the convergence to normality. Moreover, asymptotic normality is reached for far smaller sample sizes for the largest Hurst exponent estimate $\widehat{h}_2^{\rm M}$ than for that of the smallest one $\widehat{h}_1^{\rm M}$.

\section{Conclusion}
\label{sec.conc}

To the best of our knowledge, this contribution proposes the first full identification procedure for Biv-OfBm.
Its originality is to formulate identification as a non-linear regression as well as to  propose a Branch \& Bound procedure to provide efficient and elegant solutions to the corresponding non-convex optimization problem.

Consistence and asymptotic normality of the estimates are shown theoretically in a general multivariate setting.
The estimation performance is assessed for finite sample sizes by Monte Carlo simulations and found globally satisfactory  for all parameters.
Estimation of parameters $ \gamma $ (mixing) and $\rho_{\rm x}$ (correlation amongst components)
remain the most difficult parameters to estimate, though no other estimation procedure has yet been proposed in the literature.
However, including estimation of $ \gamma $ and $\rho_{\rm x}$ into the non-linear regression formulation permits to outperform the state-of-the art method, $(\widehat{h}^{\rm W}_1,\widehat{h}^{\rm W}_2,\widehat{\beta}^{\rm W})$, for estimating $h_1$, $h_2$ and $\beta$, at the price of massively increasing computational costs.
The proposed Branch \& Bound procedure is yet still shown to have a significantly lower computational cost compared to the infeasible greedy grid search strategy.
The estimation performance is satisfactory and controlled enough so that actual use of real world data can be investigated.

This estimation procedure, together with its performance assessments are paving the road for hypothesis testing, where, e.g., testing the absence of mixing (i.e., $W$ is diagonal) or of correlation amongst components (i.e., $\rho_{\rm x} \equiv 0$) are obviously interesting issues in practice.

Routines permitting both the identification and synthesis of OfBm will be made publicly available at time of publication.

This work has been explored on real Internet traffic data~\cite{Frecon_J_2016_p-icassp_nlrbssi-aaditjsapb}.

 \appendix

 \section{Interval arithmetic}
\label{app:IA}
\begin{figure}[t]
\centerline{
\includegraphics[width=0.32\linewidth]{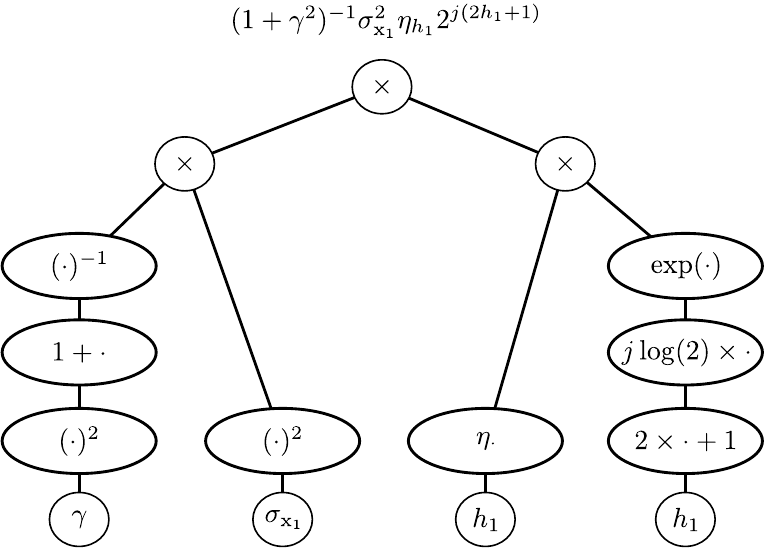}\qquad
\includegraphics[width=0.32\linewidth]{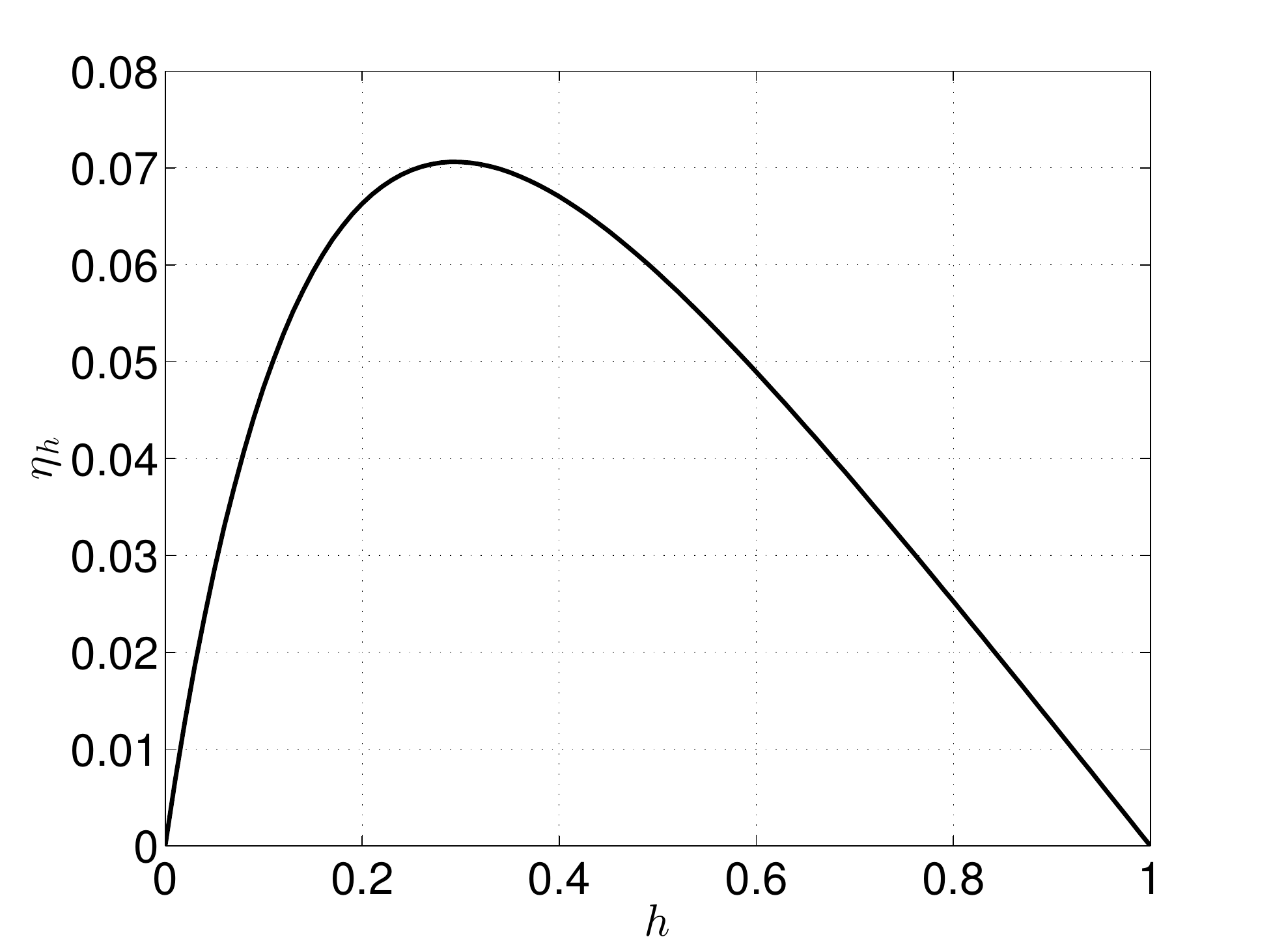}
}
\caption{{\bf Interval arithmetic and function $\eta_{{h}}$.}  Left: calculus tree associated with the first term in Eq.\eqref{eq:E11}\label{fig:calculusTreeP1}.
Right: function $\eta_{{h}}.$}
\end{figure}

Interval arithmetic is classically used to compute lower bounds for an objective criterion $C_N$ on a convex set $\mathcal{R}$ \cite{Moore_RE_1966_book_ia,Jaulin_J_2001_book_aia, Moore_R_2009_book_iia}.
It relies on the explicit decomposition of $C_N$ into several elementary functions such as sum, product, inverse, square, logarithm, exponential,\ldots, referred to as the calculus tree. For the sake of readability, we do not detail the calculus tree associated to the full $C_N$, but only its first term \eqref{eq:E11}, sketched in Fig.~\ref{fig:calculusTreeP1} (right). The leaves of the tree, i.e., the bottom line, consist of occurrences of the variables involved in $\Theta$.
Each node of the graph corresponds to an elementary function applied to its children. The intervals are composed and propagated from bottom to top to obtain a bound for $C$ within $\mathcal{R}$.
The literature on interval arithmetic provides boundaries for most of the elementary functions on intervals:
\begin{itemize}
\item $[\underline{x},\overline{x}] + [\underline{y},\overline{y}] = [\underline{x}+\underline{y}, \overline{x}+\overline{y}]$
\item $[\underline{x},\overline{x}] - [\underline{y},\overline{y}] = [\underline{x}-\overline{y}, \overline{x}-\underline{y}]$.
\item $[\underline{x},\overline{x}] \times  [\underline{y},\overline{y}] = [\min\{ \underline{x}\times\underline{y}, \underline{x}\times\overline{y}, \overline{x}\times\underline{y},\overline{x}\times\overline{y} \}, \max\{ \underline{x}\times\underline{y}, \underline{x}\times\overline{y}, \overline{x}\times\underline{y},\overline{x}\times\overline{y} \}]$
\item $[\underline{x},\overline{x}] \div  [\underline{y},\overline{y}] = [\underline{x},\overline{x}] \times [\frac{1}{\overline{y}},\frac{1}{\underline{y}}]$, if $0\not\in  [\underline{y},\overline{y}]$
\item $\log([\underline{x},\overline{x}]) = [\log(\underline{x}),\log(\overline{x})]$, if $0\not\in  [\underline{x},\overline{x}]$
\item $\exp([\underline{x},\overline{x}]) = [\exp(\underline{x}),\exp(\overline{x})]$, if $0\not\in  [\underline{x},\overline{x}]$
\end{itemize}
The only non-elementary function involved in $C_N$ is $h\mapsto\eta_h$, defined in \eqref{eq:etah} and illustrated in Fig.~\ref{fig:calculusTreeP1} (right), for a least asymmetric orthogonal Daubechies wavelet $\psi_0$ with $N_\psi = 2$. From the study of its monotonicity, we devise the following empirical bounding scheme:
\begin{equation}
(\forall[\underline{h},\overline{h}]\subseteq[0,1])\quad \eta_{[\underline{h},\overline{h}]} = \begin{cases} [\eta_{\underline{h}}, \eta_{\overline{h}}], \text{if } \overline{h}<0.3, \\
[\eta_{\overline{h}}, \eta_{\underline{h}}], \text{if } \underline{h}\geq0.3, \\
[\min(\eta_{\overline{h}},  \eta_{\underline{h}}), 0.071], \text{otherwise}.
 \end{cases}
\end{equation}
\vskip-4mm

\section{Proof of Theorem \ref{t:consistency}}
\label{app:AA}

Consider any sequence $\{\Theta_N\}_{N \in \bbN} \in \Xi$ (i.e., not necessarily composed of minima). We claim that
\vspace{-0.1cm}
\begin{equation}\label{e:Theta_nu->Theta_0}
C_N(\Theta_N) \stackrel{P}\rightarrow 0 \Rightarrow \Theta_N \stackrel{P}\rightarrow \Theta_0.
\end{equation}
By contradiction, assume that we can choose a subsequence $\{\Theta_{N(r)}\}_{r \in \bbN}$ such that, with positive probability, $C_{N(r)}(\Theta_{N(r)}) < r^{-1}$ and $\|\Theta_{N(r)} - \Theta_0\| \geq C_0 > 0$. Then, the conditions \eqref{e:EW(2j)_is_identifiable} and \eqref{e:smoothness} imply that there are indices $j^*$, $i^*_1$ and $i^*_2$, a constant $\delta > 0$
and a sequence of sets $E_{\delta,N(r)} = \{\omega: |\log_2 |E_{i^*_1,i^*_2}(2^{j^*},\Theta_{N(r)})| - \log_2 |E_{i^*_1,i^*_2}(2^{j^*},\Theta_0)||\geq \delta \}$ such that $P(E_{\delta,N(r)}) \geq C_1 >0$ for some $C_1 > 0$. So, choose $\varepsilon \in (0, \delta)$. By Theorem \ref{t:asymptotic_normality_wavecoef_fixed_scales}, for some $C_2$,
\begin{multline*}
\frac{1}{r} > \sum^{j_2}_{j = j_1}\sum_{1 \leq i_1 \leq i_2 \leq P} \{\log_2 |S_{i_1,i_2}(2^{j})| - \log_2 |E_{i_1,i_2}(2^{j},\Theta_0)| +
\log_2 |E_{i_1,i_2}(2^{j},\Theta_0)| - \log_2 |E_{i_1,i_2}(2^{j},\Theta_{N(r)})|\}^2
\\
\geq \Big\{ |\log_2 |E_{i^*_1,i^*_2}(2^{j^*},\Theta_0)| - \log_2 |E_{i^*_1,i^*_2}(2^{j^*},\Theta_{N(r)})|| - \varepsilon \Big\}^2  \geq C_2 > 0
\end{multline*}
with non-vanishing probability (contradiction). Thus, \eqref{e:Theta_nu->Theta_0} holds. Now consider the sequence \eqref{e:minima}, and note that
\begin{equation}
0 \leq  C_N(\widehat{\Theta}^M_{N}) = \inf_{\Theta \hspace{0.5mm}\in \hspace{0.5mm}\Xi}C_N(\Theta) \leq C_N(\Theta_0) \stackrel{P}\rightarrow 0,
\end{equation}
by Theorem \ref{t:asymptotic_normality_wavecoef_fixed_scales}. Therefore, by \eqref{e:Theta_nu->Theta_0}, the limit \eqref{e:consistency} holds. $\Box$

\section{Proof of Theorem \ref{t:asymptotic_normality}}
\label{app:AB}

Rewrite
$$
C_N (\Theta) = \sum^{j_2}_{j = j_1}\sum_{1 \leq i_1 \leq i_2 \leq P} (f_{N})_{i_1,i_2,j}(\Theta) \in \bbR,
$$
where
$$
(f_{N})_{i_1,i_2,j}(\Theta) = \{\log_2 |S_{i_1,i_2}(2^j)| - \log_2 |E_{i_1,i_2}(2^j,\Theta)|\}^2.
$$
It is clear that for $k\in\mathbb{N}^*$,
\begin{equation}\label{e:deriv_C=sum_deriv_f}
C^{(k)}_N (\Theta) = \sum^{j_2}_{j = j_1}\sum_{1 \leq i_1 \leq i_2 \leq P} (f^{(k)}_{N})_{i_1,i_2,j}(\Theta).
\end{equation}
Fix a triple $(i_1,i_2,j)$. By \eqref{e:parametrization_Theta0_intXi}, \eqref{e:EW(2j,Theta)_entry_not_zero} and \eqref{e:smoothness}, the first three derivatives of $(f_{N})_{i_1,i_2,j}(\Theta)$ with respect to $\Theta$ are well-defined in $\textnormal{int}\hspace{0.5mm}\Xi$. The first two derivatives at $\Theta$ can be expressed as
\begin{align}\label{e:f'N}
(f'_{N})_{i_1,i_2,j}(\Theta) &= \nabla_{\Theta}(f_{N})_{i_1,i_2,j}(\Theta)^*\\ &= -2 \{\log_2 |S_{i_1,i_2}(2^j)| - \log_2 |E_{i_1,i_2}(2^j,\Theta)|\} \Lambda_{i_1,i_2}(2^j,\Theta),
\end{align}
\begin{align}\label{e:f''N}
(f''_{N})_{i_1,i_2,j}(\Theta) &= \nabla_{\Theta}[\nabla_{\Theta}(f_{N})_{i_1,i_2,j}(\Theta)^*] \\ &= 2 \Big\{\Lambda_{i_1,i_2}(2^j,\Theta)\Lambda_{i_1,i_2}(2^j,\Theta)^*
- [\log_2 |S_{i_1,i_2}(2^j)| - \log_2 |E_{i_1,i_2}(2^j,\Theta)|] \nabla_{\Theta}\Lambda_{i_1,i_2}(2^j,\Theta)\}
\end{align}
Similarly, $(f'''_{N})_{i_1,i_2,j}(\Theta)$ consists of sums and products of $\log_2 |S(2^j)| - \log_2 |E(2^j,\Theta)|$ and derivatives of $\Lambda_{i_1,i_2}(2^j,\Theta)$.

Rewrite $C'_{N}(\Theta) = \{(C'_{N})_{l}(\Theta)\}_{l=1,\hdots,\textnormal{dim}\hspace{0.5mm}\Xi}$. By a second order Taylor expansion of $(C'_{N})_l(\Theta)$ with Lagrange remainder,
$$
\bbR \ni (C'_{N})_l(\widehat{\Theta}_{N}) - (C'_{N})_l(\Theta_0) = \nabla_{\Theta}(C'_{N})_l(\Theta_0) (\widehat{\Theta}_{N} - \Theta_0)
$$
\begin{equation}\label{e:Taylor_expansion_entry-l}
+ (\widehat{\Theta}_{N} - \Theta_0)^* \frac{\nabla_{\Theta}[\nabla_{\Theta}(C'_{N})_l((\widetilde{\Theta}_{N})_l)^*]}{2!}(\widehat{\Theta}_{N} - \Theta_0),
\end{equation}
where each entry $(\widetilde{\Theta}_{N})_{l}$, $l = 1,\hdots,\textnormal{dim}\hspace{0.5mm}\Xi$, is a parameter value lying in a segment between $\Theta_0$ and $\widehat{\Theta}_N$. Thus,
$$
\bbR^{\textnormal{dim}\hspace{0.5mm}\Xi} \ni C'_{N}(\widehat{\Theta}_{N}) - C'_{N}(\Theta_0) = C''_{N}(\Theta_0)(\widehat{\Theta}_{N} - \Theta_0)
$$
\begin{equation}\label{e:Taylor_expansion}
+ \Big((\widehat{\Theta}_{N} - \Theta_0)^* \frac{\nabla_{\Theta}[\nabla_{\Theta}(C'_{N})_l((\widetilde{\Theta}_{N})_{l})^*]}{2!}\Big)_{l=1,\hdots,\textnormal{dim}\hspace{0.5mm}\Xi}(\widehat{\Theta}_{N} - \Theta_0),
\end{equation}
where each entry $(\widehat{\Theta}_{N} - \Theta_0)^*\nabla_{\Theta}[\nabla_{\Theta}(C'_{N})_l((\widetilde{\Theta}_{N})_{l})^*]$, $l=1,\hdots,\textnormal{dim}\hspace{0.5mm}\Xi$, is a row vector. By \eqref{e:parametrization_Theta0_intXi}, $\widehat{\Theta}_{N} \in \textnormal{int}\hspace{0.5mm} \Xi$ for large $N$ with probability going to 1. Thus, $C'_{N}(\widehat{\Theta}_{N}) = 0$. Solving \eqref{e:Taylor_expansion} for $\widehat{\Theta}_N - \Theta_0$ yields
\begin{equation}\label{e:solving_Taylor_expansion_for_Theta-Theta0}
\sqrt{N}(\widehat{\Theta}_{N} - \Theta_0) = -\Big\{C''_{N}(\Theta_0) + \Big((\widehat{\Theta}_{N} - \Theta_0)^*\frac{\nabla_{\Theta}[\nabla_{\Theta}C'_{N}((\widetilde{\Theta}_{N})_l)^*]}{2!}\Big)_{l=1,\hdots,\textnormal{dim}\hspace{0.5mm}\Xi} \Big\}^{-1}\sqrt{N}C'_{N}(\Theta_0).
\end{equation}
Under \eqref{e:smoothness}, by the consistency of $\widehat{\Theta}_{N}$ for $\Theta_0$ and the expression for $(f_N''')_{i_1,i_2,j}(\Theta)$, we have for $l = 1,\hdots,\textnormal{dim}\hspace{0.5mm}\Xi$,
\begin{equation}\label{e:f'''(Theta-tilde)_consistency}
\nabla_{\Theta}[\nabla_{\Theta}C'_{N}((\widetilde{\Theta}_{N})_l)^*] \stackrel{P}\rightarrow \nabla_{\Theta}[\nabla_{\Theta}C'_{N}(\Theta_{0})^*].
\end{equation}
The invertibility of the matrix between braces in \eqref{e:solving_Taylor_expansion_for_Theta-Theta0} is ensured with probability going to 1 by the condition \eqref{e:parametrization_Theta0_intXi_deriv} and Theorem \ref{t:asymptotic_normality_wavecoef_fixed_scales}, since the expression \eqref{e:f''N} entails
\begin{equation}\label{e:C''_limit}
C''_{N}(\Theta_0) \stackrel{P}\rightarrow 2\sum^{j_2}_{j=j_1}\sum_{1 \leq i_1 \leq i_2 \leq P}\Lambda_{i_1,i_2}(2^j,\Theta_0) \Lambda_{i_1,i_2}(2^j,\Theta_0)^*.
\end{equation}
Let $\|\cdot\|_{l^1}$ be the entry-wise $l^1$ matrix norm. By the relations \eqref{e:deriv_C=sum_deriv_f} and \eqref{e:f'N}, as well as a first order Taylor expansion of $\log_2|\cdot|$ around $E_{i_1,i_2}(2^j,\Theta_0)$ under the condition \eqref{e:EW(2j,Theta)_entry_not_zero}, we can reexpress $\sqrt{N}C'_{N}(\Theta_0)$ as
$$
-2 \sum^{j_2}_{j=j_1} 2^{j/2}{\hskip -.2cm} \sqrt{K_{j}}{\hskip -.5cm}\sum_{1 \leq i_1 \leq i_2 \leq P}{\hskip -.5cm} \frac{ S_{i_1,i_2}(2^j) - E_{i_1,i_2}(2^j,\Theta_0)}{ (\log 2) E_{i_1,i_2}(2^j,\Theta_0)} \Lambda_{i_1,i_2}(2^j,\Theta_0)
$$
\begin{equation}\label{e:sqrt(nu)C'=first_order_Taylor}
+ o \Big(\sum^{j_2}_{j=j_1}\sqrt{K_j}\|S(2^j) - E(2^j,\Theta_0)\|_{l^1} \Big) \in \bbR^{\textnormal{dim}\hspace{0.5mm}\Xi}.
\end{equation}
Then, the weak limit \eqref{e:asymptotic_normality} is a consequence of \eqref{e:solving_Taylor_expansion_for_Theta-Theta0}, \eqref{e:f'''(Theta-tilde)_consistency}, \eqref{e:C''_limit}, \eqref{e:sqrt(nu)C'=first_order_Taylor}, Theorem \ref{t:asymptotic_normality_wavecoef_fixed_scales} and the Cram\'{e}r-Wold device. $\Box$

\section{Discussion of remark \ref{rem:technical}}
\label{app:AC}

The condition \eqref{e:assumption_BB} amounts to requiring the full rank of a sum of rank 1 terms. To fix ideas, consider a sum of the form $vv^* + ww^*$, where $v, w \in \bbR^2 \backslash \{0\}$. Then, the sum has deficient rank 1 if and only if $v$ and $w$ are collinear. Indeed, assume that there is a vector $u \neq 0$ such that $u^*\{vv^* + ww^*\}u = 0$. Thus, $u^*vv^*u  = 0 = u^*ww^*u$, whence collinearity follows.

So, for notational simplicity, we assume that we can rewrite the wavelet spectrum as
$$
|E_{i_1,i_2}(2^j,\Theta)| = a_{i_1,i_2} 2^{j2h_1}+ b_{i_1,i_2} 2^{j(h_1 + h_2)}+ c_{i_1,i_2} 2^{j2h_2},
$$
where $i_1,i_2=1,\hdots,P$ (c.f.\ \cite{AbryDidier2015}, Lemma 4.2). Further assume that the only parameters to be estimated are $h_1<h_2$. Then, for a fixed $j$ and a pair of indices $(i_1,i_2)$,
\begin{equation*}
\frac{\partial}{\partial h_1}\log_2|E_{i_1,i_2}(2^j,\Theta)| = \frac{1}{ \log 2 }\frac{a_{i_1,i_2} \log(2^{2j})2^{j2h_1}+ b_{i_1,i_2} \log(2^j) 2^{j(h_1 + h_2)}}{a_{i_1,i_2} 2^{j2h_1}+ b_{i_1,i_2} 2^{j(h_1 + h_2)}+ c_{i_1,i_2}2^{j2h_2}},
\end{equation*}
\begin{equation*}
\frac{\partial}{\partial h_2}\log_2|E_{i_1,i_2}(2^j,\Theta)| = \frac{1}{ \log 2 }\frac{b_{i_1,i_2} \log(2^{2j})2^{j(h_1+h_2)}+ c_{i_1,i_2} \log(2^j) 2^{j2h_2}}{a_{i_1,i_2} 2^{j2h_1}+ b_{i_1,i_2} 2^{j(h_1 + h_2)}+ c_{i_1,i_2} 2^{j2h_2}}.
\end{equation*}
This suggests that for at least two triplets $(i_1,i_2,j)$, the vectors
$$
\Big( \frac{\partial}{\partial h_1}\log_2|E_{i_1,i_2}(2^j,\Theta)|, \frac{\partial}{\partial h_2}\log_2|E_{i_1,i_2}(2^j,\Theta)| \Big)^*
$$
will not be collinear for most parametrizations in practice.

\bibliographystyle{abbrv}
\bibliography{abbr,ofbm}

\begin{thebibliography}{10}

\bibitem{AbryDidier2015}
P.~Abry and G.~Didier.
\newblock Wavelet estimation for operator fractional {B}rownian motion.
\newblock {\em Bernoulli}, page to appear, 2015.

\bibitem{Achard2008}
S.~Achard, D.~Bassett, A.~Meyer-Lindenberg, and E.~Bullmore.
\newblock Fractal connectivity of long-memory networks.
\newblock {\em Phys. Rev. E}, 77(3):036104, 2008.

\bibitem{Amblard_P-0_2011_j-ieee-tsp_imfbm}
P.-O. Amblard and J.-F. Coeurjolly.
\newblock {Identification of the Multivariate Fractional Brownian Motion}.
\newblock {\em IEEE Trans. Signal Process.}, 59(11):5152--5168, Nov. 2011.

\bibitem{Ambard_P-O_2013_bsmfbm}
P.-O. Amblard, J.-F. Coeurjolly, F.~Lavancier, and A.~Philippe.
\newblock Basic properties of the multivariate fractional {Brownian} motion.
\newblock {\em Bulletin de la Soci\'et\'e Math\'ematique de France,
  S\'eminaires et Congr\`es}, 28:65--87, 2012.

\bibitem{Bardet2003a}
J.-M. Bardet, G.~Lang, G.~Oppenheim, A.~Philippe, S.~Stoev, and M.~Taqqu.
\newblock Semi-parametric estimation of the long-range dependence parameter: a
  survey.
\newblock In P.~Doukhan, G.~Oppenheim, and M.~S. Taqqu, editors, {\em Theory
  and applications of Long-range dependence}, pages 557--577, Boston, 2003.
  Birkh\"auser.

\bibitem{Ciuciu14}
P.~Ciuciu, P.~Abry, and B.~J. He.
\newblock Interplay between functional connectivity and scale-free dynamics in
  intrinsic {fMRI} networks.
\newblock {\em NeuroImage}, 95:248--263, 2014.

\bibitem{Coeurjolly_J-F_2013_ESAIM_wamfbm}
J.-F. Coeurjolly, P.-O. Amblard, and S.~Achard.
\newblock {Wavelet analysis of the multivariate fractional Brownian motion}.
\newblock {\em {ESAIM: Probability and Statistics}}, 17:592--604, Aug. 2013.

\bibitem{Helgason:Didier:Abry:2015}
G.~Didier, H.~Helgason, and P.~Abry.
\newblock Demixing multivariate-operator self-similar processes.
\newblock In {\em Proc. Int. Conf. Acoust., Speech Signal Process.}, Brisbane,
  Australia, 2015.

\bibitem{Didier_G_2011_bernouilli_irpofbm}
G.~Didier and V.~Pipiras.
\newblock Integral representations and properties of operator fractional
  {B}rownian motions.
\newblock {\em Bernoulli}, 17(1):1--33, 2011.

\bibitem{didier:pipiras:2012}
G.~Didier and V.~Pipiras.
\newblock Exponents, symmetry groups and classification of operator fractional
  {B}rownian motions.
\newblock {\em Journal of Theoretical Probability}, 25(2):353--395, 2012.

\bibitem{Frecon_J_2016_p-icassp_nlrbssi-aaditjsapb}
J.~Frecon, R.~Fontugne, G.~Didier, N.~Pustelnik, K.~Fukuda, and P.~Abry.
\newblock Non-linear regression for bivariate self-similarity identification -
  application to anomaly detection in {I}nternet traffic based on a joint
  scaling analysis of packet and byte counts.
\newblock In {\em Proc. Int. Conf. Acoust., Speech Signal Process.}, Shanghai,
  China, March 2016.

\bibitem{Helgason_H_2011_j-sp_fessmgtsuce}
H.~Helgason, V.~Pipiras, and P.~Abry.
\newblock Fast and exact synthesis of stationary multivariate {G}aussian time
  series using circulant embedding.
\newblock {\em Signal Process.}, 91(5):1123 -- 1133, 2011.

\bibitem{Helgason_H_2011_j-sp_smsspmdccme}
H.~Helgason, V.~Pipiras, and P.~Abry.
\newblock Synthesis of multivariate stationary series with prescribed marginal
  distributions and covariance using circulant matrix embedding.
\newblock {\em Signal Process.}, 91(8):1741 -- 1758, 2011.

\bibitem{Ichida_K_1979_computing_aiamgo}
K.~Ichida and Y.~Fujii.
\newblock An interval arithmetic method for global optimization.
\newblock {\em Computing}, 23(1):85--97, 1979.

\bibitem{Jaulin_J_2001_book_aia}
L.~Jaulin, M.~Kieffer, and O.~Didrit.
\newblock {\em Applied interval analysis : with examples in parameter and state
  estimation, robust control and robotics}.
\newblock Springer, London, 2001.

\bibitem{Kearfott_RB_1992_j-global-opt_aibbabcop}
R.~B. Kearfott.
\newblock An interval branch and bound algorithm for bound constrained
  optimization problems.
\newblock {\em J. Global Optim.}, 2(3):259--280, 1992.

\bibitem{Mallat_S_2008_book_awtspsw}
S.~Mallat.
\newblock {\em A Wavelet Tour of Signal Processing, Third Edition: The Sparse
  Way}.
\newblock Academic Press, 3rd edition, 2008.

\bibitem{Moore_R_1992_book_rmgo}
R.~Moore, E.~Hansen, and A.~Leclerc.
\newblock Recent advances in global optimization.
\newblock chapter Rigorous Methods for Global Optimization, pages 321--342.
  Princeton University Press, Princeton, NJ, USA, 1992.

\bibitem{Moore_RE_1966_book_ia}
R.~E. Moore.
\newblock {\em Interval analysis}.
\newblock Prentice-Hall Inc., Englewood Cliffs, N.J., 1966.

\bibitem{Moore_R_2009_book_iia}
R.~E. Moore, R.~B. Kearfott, and M.~J. Cloud.
\newblock {\em Introduction to Interval Analysis}.
\newblock Society for Industrial and Applied Mathematics, 2009.

\bibitem{Ratschek_H_1988_book_ncmgo}
H.~Ratschek and J.~Rokne.
\newblock {\em New Computer Methods for Global Optimization}.
\newblock Halsted Press, New York, NY, USA, 1988.

\bibitem{Samorodnitsky1994}
G.~Samorodnitsky and M.~Taqqu.
\newblock {\em Stable non-{G}aussian random processes}.
\newblock Chapman and Hall, New York, 1994.

\bibitem{vandervaart:2000}
A.~Van~der Vaart.
\newblock {\em Asymptotic Statistics}, volume~3.
\newblock Cambridge University Press, 2000.

\bibitem{Veitch1999}
D.~Veitch and P.~Abry.
\newblock A wavelet-based joint estimator of the parameters of long-range
  dependence.
\newblock {\em IEEE Trans. Inform. Theory}, 45(3):878--897, 1999.

\bibitem{WENDT:2007:E}
H.~Wendt, P.~Abry, and S.~Jaffard.
\newblock Bootstrap for empirical multifractal analysis.
\newblock {\em IEEE Signal Process. Mag.}, 24(4):38--48, 2007.

\end{thebibliography}

 \end{document}